# FLORENTIN SMARANDACHE

## A UNIFYING FIELD IN LOGICS: NEUTROSOPHIC LOGIC. NEUTROSOPHY, NEUTROSOPHIC SET, NEUTROSOPHIC PROBABILITY AND STATISTICS

### (fourth edition)

$$NL(A_1 \veebar A_2) =$$
$$( T_1 \odot (\{1\} \ominus T_2) \oplus T_2 \odot (\{1\} \ominus T_1) \ominus T_1 \odot T_2 \odot (\{1\} \ominus T_1) \odot (\{1\} \ominus T_2),$$
$$I_1 \odot (\{1\} \ominus I_2) \oplus I_2 \odot (\{1\} \ominus I_1) \ominus I_1 \odot I_2 \odot (\{1\} \ominus I_1) \odot (\{1\} \ominus I_2),$$
$$F_1 \odot (\{1\} \ominus F_2) \oplus F_2 \odot (\{1\} \ominus F_1) \ominus F_1 \odot F_2 \odot (\{1\} \ominus F_1) \odot (\{1\} \ominus F_2) ).$$

$$NL(A_1 \mapsto A_2) = ( \{1\} \ominus T_1 \oplus T_1 \odot T_2, \{1\} \ominus I_1 \oplus I_1 \odot I_2, \{1\} \ominus F_1 \oplus F_1 \odot F_2 ).$$

$$NL(A_1 \leftrightarrow A_2) = ( (\{1\} \ominus T_1 \oplus T_1 \odot T_2) \odot (\{1\} \ominus T_2 \oplus T_1 \odot T_2),$$
$$(\{1\} \ominus I_1 \oplus I_1 \odot I_2) \odot (\{1\} \ominus I_2 \oplus I_1 \odot I_2),$$
$$(\{1\} \ominus F_1 \oplus F_1 \odot F_2) \odot (\{1\} \ominus F_2 \oplus F_1 \odot F_2) ).$$

# FLORENTIN SMARANDACHE

# A UNIFYING FIELD IN LOGICS: NEUTROSOPHIC LOGIC. NEUTROSOPHY, NEUTROSOPHIC SET, NEUTROSOPHIC PROBABILITY AND STATISTICS

## (fourth edition)

$$NL(A_1 \veebar A_2) =$$
$$( T_1 \odot (\{1\} \ominus T_2) \oplus T_2 \odot (\{1\} \ominus T_1) \ominus T_1 \odot T_2 \odot (\{1\} \ominus T_1) \odot (\{1\} \ominus T_2),$$
$$I_1 \odot (\{1\} \ominus I_2) \oplus I_2 \odot (\{1\} \ominus I_1) \ominus I_1 \odot I_2 \odot (\{1\} \ominus I_1) \odot (\{1\} \ominus I_2),$$
$$F_1 \odot (\{1\} \ominus F_2) \oplus F_2 \odot (\{1\} \ominus F_1) \ominus F_1 \odot F_2 \odot (\{1\} \ominus F_1) \odot (\{1\} \ominus F_2) ).$$

$$NL(A_1 \mapsto A_2) = ( \{1\} \ominus T_1 \oplus T_1 \odot T_2, \{1\} \ominus I_1 \oplus I_1 \odot I_2, \{1\} \ominus F_1 \oplus F_1 \odot F_2 ).$$

$$NL(A_1 \leftrightarrow A_2) = ( (\{1\} \ominus T_1 \oplus T_1 \odot T_2) \odot (\{1\} \ominus T_2 \oplus T_1 \odot T_2),$$
$$(\{1\} \ominus I_1 \oplus I_1 \odot I_2) \odot (\{1\} \ominus I_2 \oplus I_1 \odot I_2),$$
$$(\{1\} \ominus F_1 \oplus F_1 \odot F_2) \odot (\{1\} \ominus F_2 \oplus F_1 \odot F_2) ).$$

**Contents:**

**Preface by C. Le:  3**

## Preface to Neutrosophy and Neutrosophic Logic

### by C. Le

**1 Introduction.**

It was a surprise for me when in 1995 I received a manuscript from the mathematician, experimental writer and innovative painter Florentin Smarandache, especially because the treated subject was of philosophy - revealing paradoxes - and logics.

He had generalized the fuzzy logic, and introduced two new concepts:

a) "*neutrosophy*" – study of neutralities as an extension of dialectics;

b) and its derivative "*neutrosophic*", such as "neutrosophic logic", "neutrosophic set", "neutrosophic probability", and "neutrosophic statistics" and thus opening new ways of research in four fields: philosophy, logics, set theory, and probability/statistics.

It was known to me his setting up in 1980's of a new literary and artistic avant-garde movement that he called "*paradoxism*", because I received some books and papers dealing with it in order to review them for the German journal "Zentralblatt für Mathematik". It was an inspired connection he made between literature/arts and science, philosophy.

We started a long correspondence with questions and answers.

Because paradoxism supposes multiple value sentences and procedures in creation, anti-sense and non-sense, paradoxes and contradictions, and it's tight with neutrosophic logic, I would like to make a small presentation.

## 2. pARadOXisM, the Last Avant-Garde of the Second Millennium.

**2.1. Definition:**

PARADOXISM is an avant-garde movement in literature, art, philosophy, science, based on excessive used of antitheses, antinomies, contradictions, parables, odds, paradoxes in creations.

It was set up and led by the writer Florentin Smarandache since 1980's, who said: "The goal is to enlargement of the artistic sphere through non-artistic elements. But especially the counter-time, counter-sense creation. Also, to experiment."

**2.2. Etymology:**

Paradoxism = paradox+ism, means the theory and school of using paradoxes in literary, artistic, philosophical, scientific creations.

**2.3. History:**

"Paradoxism started as an anti-totalitarian protest against a closed society, Romania of 1980's, where the whole culture was manipulated by a small group. Only their ideas and their publications counted. We couldn't publish almost anything.

Then, I said: Let's do literature... without doing literature! Let's write... without actually writing anything. How? Simply: object literature! 'The flying of a bird', for example, represents a "natural poem", that is not necessary to write down, being more palpable and perceptible in any language that some signs laid on the paper, which, in fact, represent an "artificial poem": deformed, resulted from a translation by the observant of the observed, and by translation one falsifies. 'The cars jingling on the street' was a "city poem", 'peasants mowing' a "disseminations poem", 'the dream with open eyes' a "surrealist poem", 'foolishly speaking' a "dadaist poem", 'the conversation in Chinese for an ignorant of this language' a "leftist poem", 'alternating discussions of travelers, in a train station, on different themes' a "post-modern poem" (inter-textualism).

Do you want a vertically classification? "Visual poem", "sonorous poem", "olfactory poem", "taste poem", "tactile poem".

Another classification in diagonal: "poem-phenomenon", "poem-(soul) status", "poem-thing".

In painting, sculpture similarly - all existed in nature, already fabricated.

Therefore, a mute protest we did!

Later, I based it on contradictions. Why? Because we lived in that society a double life: an official one - propagated by the political system, and another one real. In mass-media it was promulgated that 'our life is wonderful', but in reality 'our life was miserable'. The paradox flourishing! And then we took the creation in derision, in inverse sense, in a syncretic way. Thus the paradoxism was born. The folk jokes, at great fashion in Ceausescu's 'Epoch', as an intellectual breathing, were superb springs.

The "No" and "Anti" from my paradoxist manifestos had a creative character, not at all nihilistic (C. M. Popa). The passing from paradoxes to paradoxism was documentarily described by Titu Popescu in his classical book concerning the movement: "Paradoxism's Aesthetics" (1994). While I. Soare, I. Rotaru, M. Barbu, Gh. Niculescu studied paradoxism in my literary work. N. Manolescu asserted, about one of my manuscripts of non-poems, that they are against-the-hair.

I didn't have any forerunner to influence me, but I was inspired from the 'upside-down situation' that existed in the country. I started from politic, social, and immediately got to literature, art, and philosophy, even science.

Through experiments one brings new literary, artistic, philosophical or scientific terms, new procedures, methods or even algorithms of creation. In one of my manifestos I proposed the sense of embezzling, changes from figurative to proper sense, upside-down interpretation of linguistic expressions.

In 1993 I did a paradoxist tour to literary associations and universities in Brazil.

Within 2o years of existence, 25 books and over 200 commentaries (articles, reviews) have been published, plus 3 national and international anthologies."

(Florentin Smarandache)

## 2.4. Features of Paradoxism (by Florentin Smarandache)

# *Basic Thesis of Paradoxism:*

everything has a meaning and a non-meaning in a harmony with each other.

# *Essence of Paradoxism:*
    a) sense has a non-sense, and reciprocally:
    b) non-sense has a sense.

# *Motto of Paradoxism:*
"All is possible, the impossible too!"

# *Symbol of Paradoxism:*
(a spiral -- optic illusion, or vicious circle).

# *Delimitation from Other Avant-Gardes:*
- paradoxism has a significance, while dadaism, lettrism, the absurd movement do not;
- paradoxism especially reveals the contradictions, the antinomies, anti-theses, anti-phrases, antagonism, non-conformism, the paradoxes in other words of anything (in literature, art, science), while futurism, cubism, surrealism, abstractism and all other avant-gardes do not focus on them.

*2.5. Directions for Paradoxism:*
- to use science methods (especially algorithms) for generating (and studying also) contradictory literary and artistic works;
- to create contradictory literary and artistic works in scientific spaces (using scientific: symbols, meta-language, matrices, theorems, lemmas, etc.).

*2.6. Third Paradoxist Manifesto:*
Therefore, don't enforce any literary rules on me! Or, if you do, I'll certainly encroach upon them. I'm not a poet, that's why I write poetry.
I'm an anti-poet or non-poet.
I thus came to America to re-build the Statue of Liberty of the Verse, delivered from the tyranny of the classic and its dogma.
I allowed any boldness:
- anti-literature and its literature;
- flexible forms fixed, or the alive face of the death!
- style of the non-style;
- poems without verse
(because poems don't mean words)- dumb poems with loud voice;
- poems without poems (because the notion of "poem" doesn't match any definition found in dictionaries or encyclopedias) - poems which exist by their absence;
- after-war literature: pages and pages bombed by filthiness, triteness, and non-poeticality;
- paralinguistic verse (only!): graphics, lyrical portraits, drawings, drafts...
- non-words and non-sentence poems;
- very upset free verse and trivial hermetic verse;
- intelligible unintelligible language;
- unsolved and open problems of mathematics like very nice poems of the spirit - we

must scientificize the art in this technical century;
- impersonal texts personalized;
- electrical shock;
- translation from the impossible into the possible, or transformation of the abnormal to the normal;
- pro Non-Art Art;
- make literature from everything, make literature from nothing!

The poet is not a prince of ducks! The notion of "poetry" and its derivatives have become old-fashioned in this century, and people laugh at them in disregard. I'm ashamed to affirm that I create lyrical texts, I hide them. People neither read nor listen to lyrical texts anymore, but they will read this volume because it's nothing to read!

However, the Paradoxist Movement is neither nihilism, nor disparity.

The book of the non-poems is a protest against art's marketing.

Do you writers sell your feelings? Do you create only for money??

Only books about crimes, sex, horror are published. Where is the true Art?

In begging... .

You may find in this book of uncollected poems everything you don't need and don't like: poems not to be read, not to be heard, and not to be written at all!

Enjoy them. Only after nuisance you really know what pleasure means.

They provide a mirror of everybody's infinite soul. Art, generally speaking, is pushed up to its last possible frontiers toward non-art, and even more...

Better a book of blank pages, than one that says nothing.

A very abstract and symbolic language is further used, but very concrete at the same time: non-restrictive verse from any form or content. It takes advantage of cliche against itself.

EVERYTHING IS POSSIBLE, THEREFORE: THE IMPOSSIBLE TOO! Hence don't wonder about this anti-book! If you don't understand it, that means you understand all. That is the goal of the manifesto. Because Art is not for the mind, but for feelings. Because Art is also for the mind.

Try to interpret the uninterpretable! Your imagination may flourish as a cactus in a desert.

But, The American Manifesto of the PARADOXISM is especially a revolt of the emigrant to the United States who doesn't speak English, against the language - an anti-language book written in more than a broken English (the American speech of Tomorrow?)...

[From the book: NonPoems, by Florentin Smarandache, Xiquan Publishing House, Phoenix, Chicago, 1991, 1992, 1993;

the volume contains very experimental so called , such as:
- poems without verse;
- poems without poems;
- poem-drafts;
- drawn-poems;
- poems in Pirissanorench (language spoken in the South-West of the United States by a single person!);
- super-poems;

- graphic poems;
- upset-poems.]

## 3. "Dictionary of Computing", by Denis Howe.

A well-documented and large dictionary, dealing with terms needed in Computer Science,
*The Free Online Dictionary of Computing*, is edited by Denis Howe from England
<*dbh@doc.ic.ac.uk*>. With his accordance, I cite from his dictionary the definitions regarding this issue of MVJ:

### 3.1. Neutrosophy:

<*philosophy*> (From Latin "neuter" - neutral, Greek "sophia" - skill/wisdom) A branch of philosophy, introduced by Florentin Smarandache in 1980, which studies the origin, nature, and scope of neutralities, as well as their interactions with different ideational spectra.
Neutrosophy considers a proposition, theory, event, concept, or entity, "A" in relation to its opposite, "Anti-A" and that which is not A, "Non-A", and that which is neither "A" nor "Anti-A", denoted by "Neut-A". Neutrosophy is the basis of neutrosophic logic, neutrosophic probability, neutrosophic set, and neutrosophic statistics.
*Home*.
["Neutrosophy / Neutrosophic Probability, Set, and Logic", Florentin Smarandache, American Research Press, 1998].
(1999-07-29)

### 3.2. Neutrosophic Logic:

<*logic*> (Or "**Smarandache logic**") A generalization of fuzzy logic based on Neutrosophy. A proposition is t true, i indeterminate, and f false, where t, i, and f are real values from the ranges T, I, F, with no restriction on T, I, F, or the sum n=t+i+f. Neutrosophic logic thus generalizes:
- intuitionistic logic, which supports incomplete theories (for 0<n<100, 0<=t,i,f<=100);
- fuzzy logic (for n=100 and i=0, and 0<=t,i,f<=100);
- Boolean logic (for n=100 and i=0, with t,f either 0 or 100);
- multi-valued logic (for 0<=t,i,f<=100);
- paraconsistent logic (for n>100, with both t,f<100);
- dialetheism, which says that some contradictions are true (for t=f=100 and i=0; some paradoxes can be denoted this way).
Compared with all other logics, neutrosophic logic introduces a percentage of "indeterminacy" - due to unexpected parameters hidden in some propositions. It also allows each component t,i,f to "boil over" 100 or "freeze" under 0. For example, in some tautologies t>100, called "overtrue".
*Home*.
["Neutrosophy / Neutrosophic probability, set, and logic", F. Smarandache, American Research Press, 1998].
(1999-10-04)

### 3.3. Neutrosophic Set:

<*logic*> A generalization of the intuitionistic set, classical set, fuzzy set, paraconsistent set, dialetheist set, paradoxist set, tautological set based on Neutrosophy. An element x(T, I, F) belongs to the set in the following way: it is t true in the set, i indeterminate in the set, and f false, where t, i, and f are real numbers taken from the sets T, I, and F with no restriction on T, I, F, nor on their sum n=t+i+f.
The neutrosophic set generalizes:

- the intuitionistic set, which supports incomplete set theories (for 0<n<100, 0<=t,i,f<=100);
- the fuzzy set (for n=100 and i=0, and 0<=t,i,f<=100);
- the classical set (for n=100 and i=0, with t,f either 0 or 100);
- the paraconsistent set (for n>100, with both t,f<100);
- the dialetheist set, which says that the intersection of some disjoint sets is not empty (for t=f=100 and i=0; some paradoxist sets can be denoted this way).
*Home*.
["Neutrosophy / Neutrosophic Probability, Set, and Logic", Florentin Smarandache, American Research Press, 1998].
(1999-12-14)

### 3.4. Neutrosophic Probability:

<*logic*> An extended form of probability based on Neutrosophy, in which a statement is held to be t true, i indeterminate, and f false, where t, i, f are real values from the ranges T, I, F in ]-0, 1+[, with no restriction on T, I, F or the sum n=t+i+f.
*Home*.
["Neutrosophy / Neutrosophic Probability, Set, and Logic", Florentin Smarandache, American Research Press, 1998].
(1999-10-04)

### 3.5. Neutrosophic Statistics:

<*statistics*> Analysis of events described by neutrosophic probability.
["Neutrosophy / Neutrosophic Probability, Set, and Logic", Florentin Smarandache, American Research Press, 1998].
(1999-07-05)

These definitions receive a more general form in Smarandache's next two papers, and are explained, exemplified, giving them the last versions.  Why was it necessary to extend the *fuzzy logic*?
A)  Because a paradox, as proposition, can not be described in fuzzy logic;
B)  and because the neutrosophic logic helps make a distinction between a 'relative truth' and an 'absolute truth', while fuzzy logic does not.
Due to the fact that a paradox is a proposition which is true and false in the same time, the neutrosophic logic value NL(paradox) = (1, i, 1), but this notation cannot be used in determining the fuzzy logic value FL(paradox), because if FL(paradox) = 1 (the truth) then automatically the fuzzy component of falsity is 0.  That's why it's interesting to study the neutrosophics.

0.Introduction:
# The Non-Standard Real Unit Interval

*Abstract*: In this paper one defines the non-standard real unit interval $\Vert^{-}0, 1^{+}\Vert$ (a simpler notation is also used: ]-0, 1+[ ) as a support for neutrosophy, neutrosophic logic, neutrosophic set, neutrosophic probability and statistics encountered in the next papers.

*Keywords and Phrases*: Non-standard analysis, hyper-real number, infinitesimal, monad, non-standard real unit interval, operations with sets

## 0.1. A small introduction to Non-Standard Analysis.

In 1960s Abraham Robinson has developed the *non-standard analysis*, a formalization of analysis and a branch of mathematical logic, which rigorously defines the infinitesimals. Informally, an infinitesimal is an infinitely small number. Formally, x is said to be infinitesimal if and only if for all positive integers n one has $|x| < 1/n$. Let $\varepsilon >= 0$ be a such infinitesimal number. The *hyper-real number set* is an extension of the real number set, which includes classes of infinite numbers and classes of infinitesimal numbers. Let's consider the non-standard finite numbers $(1^{+}) = 1 + \varepsilon$, where "1" is its standard part and "$\varepsilon$" its non-standard part, and $(^{-}0) = 0 - \varepsilon$, where "0" is its standard part and "$\varepsilon$" its non-standard part.

Then, we call $\Vert^{-}0, 1^{+}\Vert$ a non-standard unit interval. Obviously, 0 and 1, and analogously non-standard numbers infinitely small but less than 0 or infinitely small but greater than 1, belong to the non-standard unit interval. Actually, by "$^{-}$a" one signifies a hyper monad, i.e. a set of hyper-real numbers in non-standard analysis:

$(^{-}a) = \{a-x: x \in \mathbf{R}^{*}, x \text{ is infinitesimal}\}$,

and similarly "$b^{+}$" is a hyper monad:

$(b^{+}) = \{b+x: x \in \mathbf{R}^{*}, x \text{ is infinitesimal}\}$.

Generally, the left and right borders of a non-standard interval $\Vert^{-}a, b^{+}\Vert$ are vague, imprecise, themselves being non-standard (sub)sets $(^{-}a)$ and $(b^{+})$ as defined above.

Combining the two before mentioned definitions we now introduce the hyper-real binad of $(c^{+})$:

$(^{-}c^{+}) = \{c-x: x \in \mathbf{R}^{*}, x \text{ is infinitesimal}\} \cup \{c+x: x \in \mathbf{R}^{*}, x \text{ is infinitesimal}\}$, which is a collection of open punctured neighborhoods (balls) of c.

Of course, $(^{-}a) < a$ and $(b^{+}) > b$. No order between $(^{-}c^{+})$ and c.

By extension let $\inf \Vert^{-}a, b^{+}\Vert = (^{-}a)$ and $\sup \Vert^{-}a, b^{+}\Vert = (b^{+})$.

**Operations with hyper-real monads and binads**:

Addition of non-standard finite numbers with themselves or with real numbers:

$(^-a) + b = (^-(a + b))$, since $(^-a) + b = (a-\varepsilon)+b = (a+b)-\varepsilon = (^-(a + b))$, where $\varepsilon$ is an infinitesimal number.

$a + (b^+) = ((a + b)^+)$, since $a + (b^+) = a+(b+\varepsilon) = (a+b)+\varepsilon = ((a + b)^+)$.

$(^-a) + (b^+) = (^-(a + b)^+)$, since $(^-a) + (b^+) = (a-\varepsilon_1)+(b+\varepsilon_2)=(a+b)+(-\varepsilon_1+\varepsilon_2)$, but $-\varepsilon_1+\varepsilon_2$ can be both $> 0$ or $< 0$ since $\varepsilon_1, \varepsilon_2$ can be any infinitesimals, so we have both $(^-(a + b))$, and $((a + b)^+)$. So, the results is a hyper-binad.

$(^+a) + (b^-) = (^-(a + b)^+)$

$(^-a) + (^-b) = (^-(a + b))$ (the left hyper monads absorb themselves)

$(a^+) + (b^+) = ((a + b)^+)$ (analogously, the right hyper monads absorb themselves)

For the hyper-binads, we just introduced, it is similar:

$a + (^-b^+) = (^-(a + b)^+)$, since $a + (^-b^+) = a + (b-\varepsilon_1+\varepsilon_2) = (a+b)+(-\varepsilon_1+\varepsilon_2) = (^-(a + b)^+)$.

$(^-a^+) + b = (^-(a + b)^+)$.

$(^-a) + (^-b^+) = (^-(a + b)^+)$, since $(^-a) + (^-b^+) = (a-\varepsilon_1)+(b-\varepsilon_2+\varepsilon_3) = (a+b)+(-\varepsilon_1-\varepsilon_2+\varepsilon_3) = (^-(a + b)^+)$ because it is possible to have both $-\varepsilon_1-\varepsilon_2+\varepsilon_3 > 0$ and $< 0$.

$(^-a^+) + (^-b^+) = (^-(a + b)^+)$.

$(^-a^+) + (b^+) = (^-a^+) + (^-b) = (^-(a + b)^+)$.

Similarly for subtraction, multiplication, division, roots, and powers of non-standard finite numbers with themselves or with real numbers.

We can consider $(^-a)$ equals to the open interval $(a-\varepsilon, a)$, where $\varepsilon$ is a positive infinitesimal number. Thus:

$(^-a) = (a-\varepsilon, a)$

$(a^+) = (b, b+\varepsilon)$

$(^-a^+) = (a-\varepsilon_1, a) \cup (a, a+\varepsilon_2)$, where all $\varepsilon, \varepsilon_1, \varepsilon_2$ are positive infinitesimal numbers.

And, in consequence, all operations on hyper-monads and hyper-binads are equivalent to operations on their corresponding open intervals.

For example, $(\bar{}a^+)^3 = (\bar{}(a^3)^+)$, $\sqrt{}(b^+) = ((\sqrt{}b)^+)$, $(\bar{}a^+)/2 = (\bar{}(a/2)^+)$.

Let $i = \sqrt{}(-1)$ be the imaginary unit.

We now extend the hyper-real monads and hyper-real binads to the set of complex numbers a+bi, and we get:

**Hyper-complex number monads**:

(a)+(b)i, where either (a) or (b) or both are hyper-real monads, and none is a hyper-real binad.

H**yper-complex number binads**:

(a)+(b)i, where either (a) or (b) or both are hyper-real binads, and none is a hyper-real monad.

H**yper-complex number mixed monad-binad**:

(a)+(b)i, where one of (a), (b) is a hyper-real monad and the other is a hyper-real binad.

The operations on these hyper-complex number monads/binads/mixed are reduced to operations with hyper-real number monads and binads.

If we note by **($\bar{}$R)** the <u>set of left hyper-real number monads of all real umbers</u>, i.e. ($\bar{}$R) = {($\bar{}$a), a $\in$ R}, then {($\bar{}$R), +, ×} is a commutative ordered field, where ($\bar{}$0) is the unit element for the addition and ($\bar{}$1) is the unit element for the multiplication.

The relation of order is simply defined as ($\bar{}$a) < ($\bar{}$b) iff a < b.

{($\bar{}$R), +,×, ·} is also a linear algebra on the field of real numbers R.

Similarly, if we note by **(R$^+$)** <u>set of right hyper-real number monads of all real umbers</u>, i.e. (R$^+$) = {(a$^+$), a $\in$ R}, then {(R$^+$), +, ×} is a commutative ordered field, where (0$^+$) is the unit element for the addition and (1$^+$) is the unit element for the multiplication.

The relation of order is simply defined as (a$^+$) < (b$^+$) iff a < b.

{(R$^+$), +,×, ·} is also a linear algebra on the field of real numbers R.

And also, if we note by **($\bar{}$R)** the <u>set of hyper-real number binads of all real umbers</u>, i.e.

($\bar{}$R$^+$) = {($\bar{}$a$^+$), a $\in$ R}, then {($\bar{}$R$^+$), +, ×} is a commutative ordered field, where ($\bar{}$0$^+$) is the unit element for the addition and ($\bar{}$1$^+$) is the unit element for the multiplication.

The relation of order is simply defined as ($\bar{}$a$^+$) < ($\bar{}$b$^+$) iff a < b.

{($\bar{}$R$^+$), +,×, ·} is also a linear algebra on the field of real numbers R.

**0.2. Definition of neutrosophic components.**

Let T, I, F be standard or non-standard real subsets of $]^-0, 1^+[$,

with    sup T = t_sup, inf T = t_inf,
        sup I = i_sup,  inf I = i_inf,
        sup F = f_sup,  inf F = f_inf,
and    n_sup = t_sup+i_sup+f_sup,
        n_inf = t_inf+i_inf+f_inf.

The sets T, I, F are not necessarily intervals, but may be any real sub-unitary subsets: discrete or continuous; single-element, finite, or (countable or uncountable) infinite; union or intersection of various subsets; etc.

They may also overlap. The real subsets could represent the relative errors in determining t, i, f (in the case when the subsets T, I, F are reduced to points).

**Statically T, I, F are subsets, but dynamically the components T, I, F are set-valued vector functions/operators depending on many parameters**, such as: time, space, etc. (some of them are hidden parameters, i.e. unknown parameters): T(t, s, …), I(t, s, …), F(t, s, …), where t=time, s=space, etc., that's why the neutrosophic logic can be used in quantum physics. The Dynamic Neutrosophic Calculus can be used in psychology. Neutrosophics try to reflect the dynamics of things and ideas.

See an example:

The proposition "Tomorrow it will be raining" does not mean a fixed-valued components structure; this proposition may be say 40% true, 50% indeterminate, and 45% false at time $t_1$; but at time $t_2$ may change at 50% true, 49% indeterminate, and 30% false (according with new evidences, sources, etc.); and tomorrow at say time $t_{145}$ the same proposition may be 100%, 0% indeterminate, and 0% false (if tomorrow it will indeed rain). This is the dynamics: the truth value changes from a time to another time.

In other examples: the truth value of a proposition may change from a place to another place, for example: the proposition "It is raining" is 0% true, 0% indeterminate, and 100% false in Albuquerque (New Mexico), but moving to Las Cruces (New Mexico) the truth value changes and it may be (1, 0, 0).

Also, the truth value depends/changes with respect to the observer (subjectivity is another parameter of the functions/operators T, I, F). For example: "John is smart" can be (.35, .67, .60) according to his boss, but (.80, .25, .10) according to himself, or (.50, .20, .30) according to his secretary, etc.

In the next papers, T, I, F, called *neutrosophic components*, will represent the truth value, indeterminacy value, and falsehood value respectively referring to neutrosophy, neutrosophic logic, neutrosophic set, neutrosophic probability, neutrosophic statistics.

This representation is closer to the human mind reasoning. It characterizes/catches the *imprecision* of knowledge or linguistic inexactitude received by various observers (that's why T, I, F are subsets - not necessarily single-elements), *uncertainty* due to incomplete knowledge or acquisition errors or stochasticity (that's why the subset I exists), and *vagueness* due to lack of clear contours or limits (that's why T, I, F are subsets and I exists; in particular for the appurtenance to the neutrosophic sets).

One has to specify the superior (x_sup) and inferior (x_inf) limits of the subsets because in many problems arises the necessity to compute them.

### 0.3. Operations with sets.

Let $S_1$ and $S_2$ be two (unidimensional) real standard or non-standard subsets, then one defines:

### 0.3.1. Addition of sets:

$S_1 \oplus S_2 = \{x \mid x = s_1 + s_2, \text{ where } s_1 \in S_1 \text{ and } s_2 \in S_2\}$,

with inf $S_1 \oplus S_2 = \inf S_1 + \inf S_2$, sup $S_1 \oplus S_2 = \sup S_1 + \sup S_2$;

and, as some particular cases, we have

$\{a\} \oplus S_2 = \{x \mid x = a + s_2, \text{ where } s_2 \in S_2\}$

with inf $\{a\} \oplus S_2 = a + \inf S_2$, sup $\{a\} \oplus S_2 = a + \sup S_2$;

also $\{1^+\} \oplus S_2 = \{x \mid x = 1^+ + s_2, \text{ where } s_2 \in S_2\}$

with inf $\{1^+\} \oplus S_2 = 1^+ + \inf S_2$, sup $\{1^+\} \oplus S_2 = 1^+ + \sup S_2$.

### 0.3.2. Subtraction of sets:

$S_1 \ominus S_2 = \{x \mid x = s_1 - s_2, \text{ where } s_1 \in S_1 \text{ and } s_2 \in S_2\}$.

For real positive subsets (most of the cases will fall in this range) one gets inf $S_1 \ominus S_2 = \inf S_1 - \sup S_2$, sup $S_1 \ominus S_2 = \sup S_1 - \inf S_2$;

and, as some particular cases, we have

$\{a\} \ominus S_2 = \{x \mid x = a - s_2, \text{ where } s_2 \in S_2\}$,

with inf $\{a\} \ominus S_2 = a - \sup S_2$, sup $\{a\} \ominus S_2 = a - \inf S_2$;

also $\{1^+\} \ominus S_2 = \{x \mid x = 1^+ - s_2, \text{ where } s_2 \in S_2\}$,

with inf $\{1^+\} \ominus S_2 = 1^+ - \sup S_2$, sup $\{1^+\} \ominus S_2 = 1^+ - \inf S_2$.

### 0.3.3. Multiplication of sets:

$S_1 \odot S_2 = \{x \mid x = s_1 \cdot s_2, \text{ where } s_1 \in S_1 \text{ and } s_2 \in S_2\}$.

For real positive subsets (most of the cases will fall in this range) one gets inf $S_1 \odot S_2 = \inf S_1 \cdot \inf S_2$, sup $S_1 \odot S_2 = \sup S_1 \cdot \sup S_2$;

and, as some particular cases, we have

$\{a\} \odot S_2 = \{x \mid x = a \cdot s_2, \text{ where } s_2 \in S_2\}$,

with inf $\{a\} \odot S_2 = a * \inf S_2$, sup $\{a\} \odot S_2 = a \cdot \sup S_2$;

also $\{1^+\} \odot S_2 = \{x \mid x = 1^+ s_2, \text{ where } s_2 \in S_2\}$,

with inf $\{1^+\} \odot S_2 = 1^+ \cdot \inf S_2$, sup $\{1^+\} \odot S_2 = 1^+ \cdot \sup S_2$.

### 0.3.4. Division of a set by a number:

Let $k \in R^*$, then $S_1 \oslash k = \{x \mid x = s_1 / k, \text{ where } s_1 \in S_1\}$.

**Acknowledgements**:

The author would like to thank Drs. C. Le and Ivan Stojmenovic for encouragement and invitation to write this and the following papers.

Kluwer Academic Publishers. xiv, 311 p. (2000). [ISBN 0-7923-6340-X/hbk; ISSN 0921-3791]

# Neutrosophy, a New Branch of Philosophy

*Abstract*:

In this paper is presented a new branch of philosophy, called *neutrosphy*, which studies the origin, nature, and scope of neutralities, as well as their interactions with different ideational spectra.

The Fundamental Thesis: Any idea <A> is T% true, I% indeterminate, and F% false, where T, I, F are standard or non-standard subsets included in $\rbrack^-0, 1^+\lbrack$.

The Fundamental Theory: Every idea <A> tends to be neutralized, diminished, balanced by <Non-A> ideas (not only <Anti-A>, as Hegel asserted) - as a state of equilibrium.

Neutrosophy is the base of *neutrosophic logic*, a multiple value logic that generalizes the fuzzy logic, of *neutrosophic set* that generalizes the fuzzy set, and of *neutrosphic probability*
and *neutrosophic statistics*, which generalize the classical and imprecise probability and statistics respectively.

*Keywords and Phrases*: Non-standard analysis, hyper-real number, infinitesimal, monad, non-standard real unit interval, operations with sets

*1991 MSC:* 00A30, 03-02, 03B50

## 1.1. Foreword.

Because world is full of indeterminacy, a more precise imprecision is required. That's why, in this study, one introduces a new viewpoint in philosophy, which helps to the generalization of classical 'probability theory', 'fuzzy set' and 'fuzzy logic' to <neutrosophic probability>, <neutrosophic set> and <neutrosophic logic> respectively. They are useful in artificial intelligence, neural networks, evolutionary programming, neutrosophic dynamic systems, and quantum mechanics.

Especially in quantum theory there is an uncertainty about the energy and the momentum of particles, and, because the particles in the subatomic world don't have exact positions, we better calculate their neutrosophic probabilities (i.e. involving a percent of incertitude, doubtfulness, indetermination as well - behind the percentages of truth and falsity respectively) of being at some particular points than their classical probabilities.

Besides Mathematics and Philosophy interrelationship, one searches Mathematics in connection with Psychology, Sociology, Economics, and Literature.

This is a foundation study of the NEUTROSOPHIC PHILOSOPHY because, I think, a whole collective of researchers should pass through all philosophical schools/movements/theses/ideas and
extract their positive, negative, and neuter features.
Philosophy is subject to interpretation.
This is a *propédeutique* (Fr.), and a first attempt of such treatise.
(An exhaustive (if possible) neutrosophic philosophy should be a synthesis of all-times philosophies inside of a neutrosophic system.)

This article is a collection of concise fragments, short observations, remarks, various citations, aphorisms, some of them in a poetical form. (Main references are listed after several individual fragments.) It also introduces and explores new terms within the framework of avant-garde and experimental philosophical methods under multiple values logics.

The research is a part of a National Science Foundation grant proposal for Interdisciplinary Logical Sciences.

## 1.2. Neutrosophy, a New Branch of Philosophy

**A) Etymology:**
Neutro-sophy [French *neutre* < Latin *neuter*, neutral, and Greek *sophia*, skill/wisdom] means knowledge of neutral thought.

**B) Definition:**
Neutrosophy is a new branch of philosophy, which studies the origin, nature, and scope of neutralities, as well as their interactions with different ideational spectra.

**C) Characteristics:**
This mode of thinking:
- proposes new philosophical theses, principles, laws, methods, formulas, movements;
- reveals that world is full of indeterminacy;
- interprets the uninterpretable;
- regards, from many different angles, old concepts, systems:
showing that an idea, which is true in a given referential system, may be false in another one, and vice versa;
- attempts to make peace in the war of ideas,
and to make war in the peaceful ideas;
- measures the stability of unstable systems,
and instability of stable systems.

**D) Methods of Neutrosophic Study:**
mathematization (neutrosophic logic, neutrosophic probability and statistics, duality), generalization, complementarity, contradiction,  paradox, tautology, analogy, reinterpretation, combination, interference, aphoristic, linguistic, transdisciplinarity.

**E) Formalization:**

Let's note by <A> an idea, or proposition, theory, event, concept, entity, by <Non-A> what is not <A>, and by <Anti-A> the opposite of <A>. Also, <Neut-A> means what is neither <A> nor <Anti-A>, i.e. neutrality in between the two extremes. And <A'> a version of <A>.

<Non-A> is different from <Anti-A>.

For example:

If <A> = white, then <Anti-A> = black (antonym),

but <Non-A> = green, red, blue, yellow, black, etc. (any color, except white),

while <Neut-A> = green, red, blue, yellow, etc. (any color, except white and black),

and <A'> = dark white, etc. (any shade of white).

In a classical way:

<Neut-A> ≡ <Neut-(Anti-A)>, i.e. neutralities of <A> are identical with neutralities of <Anti-A>, also <Non-A> ⊃ <Anti-A>, and <Non-A> ⊃ <Neut-A> as well,

<A> ∩ <Anti-A> = Ø, <A> ∩ <Non-A> = Ø.

<A>, <Neut-A>, and <Anti-A> are disjoint two by two.

<Non-A> is the completeness of <A> with respect to the universal set.

But, since in many cases the borders between notions are vague, imprecise, it is possible that <A>, <Neut-A>, <Anti-A> (and <Non-A> of course) have common parts two by two.

**F) Main Principle:**

Between an idea <A> and its opposite <Anti-A>, there is a continuum-power spectrum of neutralities <Neut-A>.

**G) Fundamental Thesis:**

Any idea <A> is T% true, I% indeterminate, and F% false, where T, I, F ⊂ $\| {}^{-}0, 1^{+} \|$.

**H) Main Laws:**

Let <α> be an attribute, and (T, I, F) ⊂ $\| {}^{-}0, 1^{+} \|^{3}$. Then:

- There is a proposition <P> and a referential system {R}, such that <P> is T% <α>, I% indeterminate or <Neut-α>, and F% <Anti-α>.

- For any proposition <P>, there is a referential system {R}, such that <P> is T% <α>, I% indeterminate or <Neut-α>, and F% <Anti-α>.

- <α> is at some degree <Anti-α>, while <Anti-α> is at some degree <α>.

Therefore:

For each proposition <P> there are referential systems {$R_1$}, {$R_2$}, ..., so that <P> looks differently in each of them - getting all possible states from <P> to <Non-P> until <Anti-P>.

And, as a consequence, for any two propositions <M> and <N>, there exist two referential systems {$R_M$} and {$R_N$} respectively, such that <M> and <N> look the same. The referential systems are like mirrors of various curvatures reflecting the propositions.

**I) Mottos:**

- All is possible, the impossible too!
- Nothing is perfect, not even the perfect!

### J) Fundamental Theory:

Every idea <A> tends to be neutralized, diminished, balanced by <Non-A> ideas (not only <Anti-A> as Hegel asserted) - as a state of equilibrium.  In between <A> and <Anti-A> there are infinitely many <Neut-A> ideas, which may balance <A> without necessarily <Anti-A> versions.

To neuter an idea one must discover all its three sides:  of sense (truth), of nonsense (falsity), and of undecidability (indeterminacy) - then reverse/combine them.  Afterwards, the idea will be classified as neutrality.

### K) Delimitation from Other Philosophical Concepts and Theories:

1. Neutrosophy is based not only on analysis of oppositional propositions, as *dialectic* does, but on analysis of neutralities in between them as well.

2. While *epistemology* studies the limits of knowledge and justification, neutrosophy passes these limits and takes under magnifying glass not only defining features and substantive conditions of an entity <E> - but the whole <E'> derivative spectrum in connection with <Neut-E>.

Epistemology studies philosophical contraries, e.g. <E> versus <Anti-E>, neutrosophy studies <Neut-E> versus <E> and versus <Anti-E> which means logic based on neutralities.

3-4. *Neutral monism* asserts that ultimate reality is neither physical nor mental. Neutrosophy considers a more than pluralistic viewpoint:  infinitely many separate and ultimate substances making up the world.

5. *Hermeneutics* is the art or science of interpretation, while neutrosophy also creates new ideas and analyzes a wide range ideational field by balancing instable systems and unbalancing stable systems.

6. *Philosophia Perennis* tells the common truth of contradictory viewpoints, neutrosophy combines with the truth of neutral ones as well.

7. *Fallibilism* attributes uncertainty to every class of beliefs or propositions, while neutrosophy accepts 100% true assertions, and 100% false assertions as well - moreover, checks in what referential systems the percent of uncertainty approaches zero or 100.

### L) Philosophy's Limits:

The whole philosophy is a *tautologism*:  true in virtue of form, because any idea when first launched is proved true by its initiator(s).  Therefore, philosophy is empty or uninformative, and a priori knowledge.

One can ejaculate:  All is true, even the false!

And yet, the whole philosophy is a *nihilism*:  because any idea, first proved true, is later proved false by followers.  It is a contradiction: false in virtue of form.  Therefore, now philosophy is overinformative, and a posteriori knowledge.

Now, one can ejaculate:  All is false, even the truth!

All not-yet-contradicted philosophical ideas will be sooner or later contradicted because every philosopher attempts to find a breach in the old systems.  Even this new theory (that I am sure it is not pretty sure!) will be inverted... And, later, others will reinstall it back...

Consequently, philosophy is logically necessary and logically impossible. Agostoni Steuco of Gubbio was right, the differences between philosophers are undifferentiable. Leibniz's expression <true in all possible world> is superfluous, derogatory, for our mind may construct impossible world as well, which become possible in our imagination (F.Smarandache, "Inconsistent Systems of Axioms", 1995).

- In this theory one can prove anything!
- In this theory one can deny anything!

Philosophism = Tautologism + Nihilism.

    **M) Classification of Ideas:**
a) easily accepted, quickly forgotten;
b) easily accepted, heavily forgotten;
c) heavily accepted, quickly forgotten;
d) heavily accepted, heavily forgotten.
And various versions in between any two categories.

    **N) Evolution of an Idea** <A> in the world is not cyclic (as Marx said), but discontinuous, knotted, boundless:
<Neut-A> = existing ideational background, before arising <A>;
<Pre-A> = a pre-idea, a forerunner of <A>;
<Pre-A'> = spectrum of <Pre-A> versions;
<A> = the idea itself, which implicitly gives birth to
<Non-A> = what is outer <A>;
<A'> = spectrum of <A> versions after (miss)interpretations
    (miss)understanding by different people, schools,
    cultures;
<A/Neut-A> = spectrum of <A> derivatives/deviations, because <A>
     partially mixes/melts first with neuter ideas;
<Anti-A> = the straight opposite of <A>, developed inside of
     <Non-A>;
<Anti-A'> = spectrum of <Anti-A> versions after
    (miss)interpretations (miss)understanding by different
    people, schools, cultures;
<Anti-A/Neut-A> = spectrum of <Anti-A> derivatives/deviations,
     which means partial <Anti-A> and partial
     <Neut-A> combined in various percentage;
<A'/Anti-A'> = spectrum of derivatives/deviations after mixing
     <A'> and <Anti-A'> spectra;
<Post-A> = after <A>, a post-idea, a conclusiveness;
<Post-A'> = spectrum of <Post-A> versions;
<Neo-A> = <A> retaken in a new way, at a different level, in new
    conditions, as in a non-regular curve with inflection
    points, in evolute and involute periods, in a
    recurrent mode; the life of <A> restarts.

Marx's 'spiral' of evolution is replaced by a more complex differential curve with ups-and-downs, with knots - because evolution means cycles of involution too.

This is *dynaphilosophy* = the study of infinite road of an idea.

<Neo-A> has a larger sphere (including, besides parts of old <A>, parts of <Neut-A> resulted from previous combinations), more characteristics, is more heterogeneous (after combinations with various <Non-A> ideas). But, <Neo-A>, as a whole in itself, has the tendency to homogenize its content, and then to de-homogenize by mixture with other ideas.

And so on, until the previous <A> gets to a point where it paradoxically incorporates the entire <Non-A>, being indistinct of the whole. And this is the point where the idea dies, can not be distinguished from others. The Whole breaks down, because the motion is characteristic to it, in a plurality of new ideas (some of them containing grains of the original <A>), which begin their life in a similar way. As a multi-national empire.

It is not possible to pass from an idea to its opposite without crossing over a spectrum of idea's versions, deviations, or neutral ideas in between.

Thus, in time, <A> gets to mix with <Neut-A> and <Anti-A>.

We wouldn't say that "extremes attract each other", but <A> and <Non-A> (i.e., inner, outer, and neutron of an idea).

Therefore, Hegel was incomplete when he resumed that: a thesis is replaced by another, called anti-thesis; contradiction between thesis and anti-thesis is surpassed and thus solved by a synthesis. So Socrates in the beginning, or Marx and Engels (dialectic materialism).

There is not a triadic scheme:
 - thesis, antithesis, synthesis (Hegelians);
or
 - assertion, negation, negation of negation (Marxists);
but a pluradic pyramidal scheme, as seen above.

Hegel's and Marx's antithesis <Anti-T> does not simply arise from thesis <T> only. <T> appears on a background of preexistent ideas, and mixes with them in its evolution. <Anti-T> is built on a similar ideational background, not on an empty field, and uses in its construction not only opposite elements to <T>, but elements of <Neut-T> as well, and even elements of <T>.

For, a thesis <T> is replaced not only by an antithesis <Anti-T>, but also by various versions of neutralities <Neut-T>.

We would resume this at: neuter-thesis (ideational background before thesis), pre-thesis, thesis, pro-thesis, non-thesis (different, but not opposite), anti-thesis, post-thesis, neo-thesis.

Hegel's scheme was purist, theoretic, idealistic. It had to be generalized. From simplism to organicism.

**O) Philosophical Formulas:**

Why are there so many distinct (even contrary) philosophical Schools?

Why, concomitantly with the introduction of a notion <A>, its reverse <Non-A> is resulting?

Now, one presents philosophical formulas just because in the spiritual field it is really difficult to obtain (exact) formulas.

### a) Law of Equilibrium:

The more <A> increases, the more <Anti-A> decreases. One has the following relationship:

<A>·<Anti-A> = k·<Neut-A>,

where k is a constant depending on <A>, and <Neut-A> is a supporting point for balancing the two extremes.

If the supporting point is the neutralities' centroid, then the above formula is simplified to:

<A>·<Anti-A> = k,

where k is a constant depending on <A>.

Interesting particular cases:
Industrialization × Spiritualization = constant, for any society.
The more industrialized a society is, the less spiritual level its citizens have.
Science × Religion = constant.
White × Black = constant.
Plus × Minus = constant.
Pushing to the limits, in other words calculating in the absolute space, one gets:
Everything × Nothing = universal constant,
or $\infty \times 0$ (= $0 \times \infty$) = universal constant.

We are directing towards a mathematization of philosophy, but not in a Platonian sense.

Graph 5. O.a.1:
Materialism × Idealism = constant, for any society.

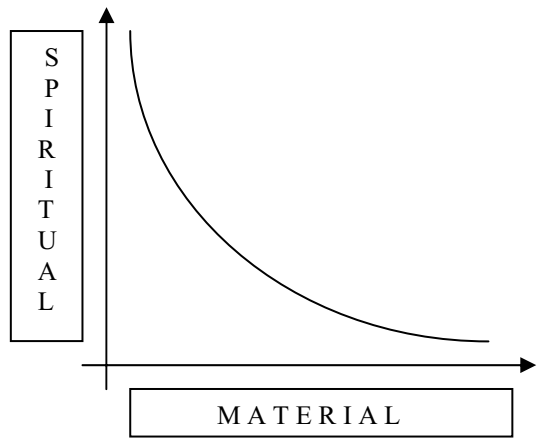

The vertical and horizontal Cartesian axes are asymptotes for the curve M·I = k.

### b) Law of Anti-Reflexivity:

<A> in the mirror of <A> gradually vanishes itself.
Or <A> of <A> may transform into a distorted <A>.

Examples:
Marriage between relatives gives birth to vapid (often handicapped) descendants.
That's why crossing the species of plants (and sometimes races of animals and humans as well) we get hybrids with better qualities and/or quantities. Biological theory of mixing species.
That's why emigration is benign for bringing new blood in a static population.
Nihilism, spread out after Turgeniev's "Parents and children" novel in 1862 as an absolute negation, denies everything, therefore itself too!
Dadaism of the dadaism vanishes either.

### c) Law of Complementarity:

<A> feels like completing with <Non-A> in order to form a whole.

Examples:
Persons, who are different, feel like completing each other and associate. (Man with woman.)
Complementary colors (that, combined in the right intensities, produce white).

### d) Law of Inverse Effect:

When trying to convert someone to an idea, belief, or faith by boring repetitions or by force, he ends up to hating it.

Examples:
The more you ask someone to do something, the less he would.
Doubling the rule, brings to halving.
What's much, it's not good...
(inversely proportional).
When you are sure, don't be!

When pressing someone to do something, he would do a different (but not necessary the opposite, as Newton's third motion laws' axiom stated) at-various-slopes reaction:

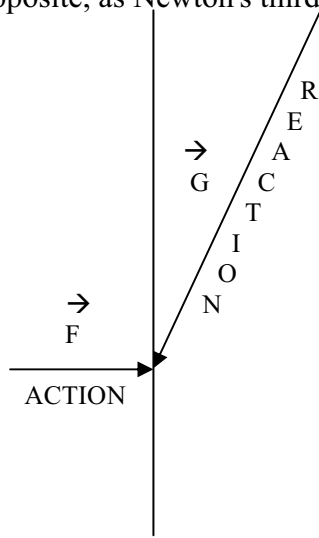

***e) Law of Reverse Identification:***
<Non-A> is a better <A> than <A>.

Example:
Poetry is more philosophical than philosophy.

***f) Law of Joined Disjointedness:***
<A> and <Non-A> have elements in common.

Examples:
There is little distinction between "good" and "bad".
Rational and irrational work together unseparately.
Consciousness and unconsciousness similarly.
"Come, my soul said, let's write poems for my body, for we are One" (Walt Whitman).
Finite is infinite [see the microinfinity].

***g) Law of Identities' Disjointedness:***
The permanent fight between <A> and <A'> (different shades of <A>).

Examples:
The permanent fight between absolute truth and relative truth.
The distinction between crisp false and neutrosophic false (the second one means a combination of falsity, indeterminacy, and truth degrees).

***h) Law of Compensation:***
If <A> now, then <Non-A> later.

Examples:
Any loss / Has its gain
[meaning later it will be better, because you learned from the loss].
There is no success without failure
[patience guys!].

***i) Law of Prescribed Condition:***
One cannot jump out of own limits.
(One spins inside own circle.)

***j) Law of Particular Ideational Gravitation:***
Every idea <A> attracts and rejects other idea <B> with a force directly proportional with the product of their neutrosophic measures and the exponential of their distance.

(By opposition to the modern restatement of Newton's law of gravitation of particles of matter, the distance influences directly - not indirectly - proportional: the more opposite (distanced) ideas, the stronger attraction)

### k) Law of Universal Ideational Gravitation:

<A> tends towards <Non-A> (not <Anti-A> as Hegel said), and reciprocally.
There are forces which act on <A>, directing it towards <Non-A>,
until a critical point is attained, and then <A> turns back.
<A> and <Non-A> are in a continuous motion, and their frontiers changing accordingly.

Examples:
Perfection leads to imperfection.
Ignorance is pleased.

Particular Case:
Everybody tends to approach his specific level of incompetence!
This is not a joke, but very truly:
X gets a job at level say L1;
if he is good, he's promoted to level L2;
if, in the new position he's good, he's promoted further to L3;
and so on... until he's not good anymore, therefore not promoted;
thus, he got to his level of incompetence.
<A> tends towards <Non-A>.
Therefore, everybody's ideal is to tend towards what he/she is
not able to do.

But the movement is nonlinear.
<Non-A> has a large range (power of continuum) of "what is not
<A>" (outer <A>) versions, let index them in the set $\{<Non-A>_i\}_i$.

(All $\{<Anti-A>_i\}_i$ versions are included in <Non-A>.)
Hence, infinitely many <Non-A>$_i$ versions gravitate, as planets around a star, on orbits of
<A>. And, between each <Non-A>$_i$ version and the centroid "star" <A>, there are
attraction and rejection forces. They approach each other until arriving to certain
minimum critical points: $P_{m(i)}$ for <A>, and $Q_{m(i)}$ for <Non-A>$_i$, and then again they go
far from each other until touching certain
maximum points: $P_{M(i)}$ for <A>, and $Q_{M(i)}$ for <Non-A>$_i$.
Through differential equations we may calculate the minimum and maximum (spiritual)
distances between <A> and <Non-A>$_i$, the Cartesian coordinates of the critical points,
and the status quo of each version.
We would say that <A> and a <Non-A>$_i$ version meet in an absolute/infinite point.
When all <Non-A>$_i$ versions fall into <A> we have a catastrophe!

### P) Neutrosophic Studies and Interpretations of Known Theories, Modes, Views, Processes of Reason, Acts, Concepts in Philosophy.

This section, which is a *neutrosophic epistemology*, has a structure alike Wittgenstein's
tractatus: short (from 1-2 lines to maximum 10-15 lines) independent philosophical reflections, metaphysical and metaphorical comments – which are separated by blank rows. It is an *analytical study*, and it is related to multiple-valued logic because in almost each small paragraph one shows that a statement <A> was proved true by a philosopher X whereas latter another philosopher Y proved the opposite statement <Anti-A> was true. Therefore, both <A> and <Anti-A> were true. {Whence one can deduce that both <A> and <Anti-A> could be false.} Even more, using a neutrosophic interpretation, one could say that other ideas in between <A> and <Anti-A> and related to them, noted by <Neut-A>, could be true as well. This relates to dialetheism, which says that some contradictions are true, to paraconsistent logic, to intuitionistic logic, till neutrosophic logic (where <A>, <Anti-A>, and ideas in between them belonging to <Neut-A> could all be true or partially true).

Many paradoxes are treated here, and one knows that a paradox is a proposition true and false in the same time – i.e. connected to multiple-valued logic as well, and not many logics approached the paradoxes. Other reflections show that a *subject* may be characterized by an *attribute* and *its opposite* simultaneously (also related to multiple-valued logic in a philosophical way).

To any launched idea there are pro- and contra- reactions, but also neuter (indifference, neutrality) as well. Hegel's <dialectic> [Gr. *dialektikě* < *dia* with, *legein* to speak] doesn't work, it consequently has to be extended to a somehow improper term *trialectic*, and even more to a *pluralectic* because there are various degrees of positive, and of negative, and of indifference as well - all of them interpenetrated. Going to a continuum-power *transalectic* (∞-alectic).

"+" not only asks for "-" for equilibrium, as Hegel said, but for "0" as well as a support point for the thinking lever.

Hegel's self-development of an idea <A> is not determined on its internal contraries only, but on its neutralities as well - because they all fare and interfere. Self-development of an idea is also determined by external (pro, contra, neuter) factors (*Comparative Philosophy*, as comparative literature).

Between particular and universal there are P% particular, I% indeterminate (neutral), and U% universal things, with P, I, U ⊂ $\|^-0, 1^+\|$.

The atom's structure holds in the history of any idea. The reasoning is based upon the analysis of positive, negative, and neutral propositions.
This should be called *Quantum Philosophy*.

In nuclear fission a free neutron strongly interacts with nuclei and is readily absorbed, then it decays into a proton, an electron, and a 'neutrino' (Enrico Fermi) with a half-time of near 12 minutes.

Neutrosophy equally encompasses a philosophical viewpoint, and mode of reflection, and concept, and method in itself, and action, and movement, and general theory, and process of reasoning.

This approach differs from *neutrosophism*, which is a point of view that neutrosophy is a fundamental science to study the world from that perspective.

Neutrosophy studies not only an idea's conditions of possibility, but of impossibility as well.  And focuses on its historical development (past and present interpretation – by using classical analysis, and future interpretation - by using neutrosophic probability and statistics).

In economics Keynes chose for the concept of "unstable equilibrium" (<The General Theory>), whereas Anghel M. Rugină passed to that of "stable disequilibrium" (<Truth in the Abstract (Analytical) versus Truth in the Concrete (Empirical)>).
A self-regulating and self-unregulating mechanism is functioning in each system, moving from equilibrium to dis-equilibrium back and forth.
A unstable-made stability, and stable-made instability.  Or equilibrium in disequilibrium, and disequilibrium in equilibrium.

We mean, a very dynamic system by rapid small changes, characterized by a derivative.
The static system is dead.
Leon Walras was right:  monopolies reduce the competitions, and thus the progress.

My opinion is that some philosophers grope, stumble.  They don't have clear ideas or systems, or even precise directions on a subject.  Paroles, paroles...
What one asserts today, another will deny tomorrow.
Many times they talk too much for saying nothing.  Some have points contrary to experience and evidence, others have an inadequate reason.
That's why a mathematization (even more, an axiomatization but not in stricto sensu) of all knowledge fields would be required, especially in philosophy (alike Mendeleev's Table of Chemical Elements).
The mathematization is required because it is not possible!

Philosophy is semiscientific and semiempiric.  It is less scientific than psychology, but more scientific than poetry.

Human is dependent and independent in the same time.

I understand spirit as quality, and material as quantity.
Of course, they melt.

I see truth like a body, an object with a form.
I see material as a dense/condensed spirit, a viscous idea.

The structure of ideas reflects the structure of objects.
And reciprocally.

In the mind-body problem:
The mental phenomenon is of physical nature, and the physical phenomenon of mental as well.

"(...) it feels sometimes the economy was propelled on the symmetry principle, which demands that every new theory always be exactly the reverse of the old one" [Mark Blaug, <Economics Theory in Retrospective>].

Neohegelians:
Reconciliation of contraries (Bradley), or irreconciliation of contraries (Wahl)?
Both!

Neutrosophy:
- has the aim of unifying field in humanistic (as Einstein tried to find in science);
- explores the differences between:
  . thinkers,
  . philosophical schools, movements, theories, principles and proves they are minimum;
- reveals that no thought school is better than another, and no philosopher is greater than another;
- is an attempt to reconcile reluctant viewpoints with inoffensive others;
- the truth may not be separated from false;
- if the philosopher X enunciated a proposition P, try to contrarily think and to compare with <Neut-P> too.

Ignorantism:
Power countries deliberately ignore the arts, literature, science, culture, traditions of third world countries. More, they even boycott, scorn them...
Third world countries creators and inventors are also handicapped by language, poor living conditions, less technology for doing research.
In histories of arts, literature, science you see only westerners:
- rare exceptions of other people being in there confirm the rule!
A minor poet, for example, who wrote in English or French or German is better known than a genius like Eminescu who wrote in a not-international Romanian language.

Negativity (Heraclitus, Spinoza, Kant, Hegel) passes through diverse phases:
from assertion to a spectrum of partial negativity, and eventually to a higher degree negativity.

Not com-plementarity (used by Bohr and Heisenberg in physics philosophy), but *tri-plementarity* (negativity, positivity, and neutrality - corresponding to 0, 1, and 1/2 respectively), even *n-plementarity* (which means: n disjoint elements together forming a whole), or generalized to *∞-plementarity* (with power of continuum), for there are complex mixed versions of them.
Beyond indefinitely many states between 0 and 1, the midpoint ½ represents neither negative, nor positive - or represents both of them (which cancel each other).

Hermeneutics of Philosophical Hermeneutics:
If prejudgments can not be eliminated in the judgements,
why do we need the science of interpretation?

Arguing with Plekhanov (historical development is not managed at will), one says that it is at some degree managed and at another degree not managed by the will.

Ab'lard's conceptualism, which states that universalia post rem (general is besides things), i.e. general is not in things, is partially true, because general persists in each individual,
that's why it is possible to form classes of individuals with similar particular characteristics.

*Philosophy of Philosophy*:
- why do we need philosophy today?
- why don't we need philosophy today?
- what direction is philosophy going to?
- what direction isn't philosophy going to?
One feels that philosophy is for people who have nothing else to do, like puzzles or rebus!

Neutrosophy means/encompasses:
- philosophy seen by a mathematician and poet;
- study of History of Philosophy;
- controversial themes of philosophy
(to explore the offensiveness and inoffensiveness);
- evolution of an idea from <A> to <Non-A> and then to <Anti-A>;
- how to get patterns where they do not look to be,
i.e. to find common characteristics at "+", "-", and "0" attributes;
- how an idea appears from different viewpoints,
from all viewpoints;
- to find the vanishing point of all philosophical ideas.

Neutrosophy can also be seen as:
- new approach to philosophy;
- philosophy of philosophies;
- non-philosophy;
- super-philosophy;
- neophilosophy;
- God and Devil of the philosophy;
- meta-philosophy, macro-philosophy;
- New World Order in philosophy;
- paradox of philosophy and philosophy of the paradox;
- thought of thought;
- showing the philosophy's perfection and imperfection
simultaneously;
- paradox within/from paradox:  there are infinitely many;
- world's enigma;
- nature's essence;
- enigma of the world;

- any substance ultimately has a neutrosophic attribute;
- life without paradox would be monotonous and boring, linear;
- paradoxist intuition is a high level of awareness;
- postmodernist;
- an algebraic, physical and chemical philosophy;
- consistent with its inconsistence.

Transcendentalism (Emerson especially, and Kant, Hegel, Fichte), which proposes to discover the nature of reality by investigating the process of thought, is combined with pragmatism (Williams James), which first "tries to interpret each notion or theory by tracing its respective practical consequences".
We mean to know reality through thought, and thought through reality.

In India's VIII-th - IX-th centuries one promulgated the Non-Duality (Advaita) through the non-differentiation between Individual Being (Atman) and Supreme Being (Brahman). The philosopher Sańkaracharya (782-814 A.C.) was then considered the savior of Hinduism, just in the moment when the Buddhism and the Jainism were in a severe turmoil and India was in a spiritual crisis.
Non-Duality means elimination of ego, in order to blend yourself with the Supreme Being (to reach the happiness).
Or, arriving to the Supreme was done by Prayer (Bhakti) or Cognition (Jnana). It is a part of Sańkaracharya's huge merit (charya means teacher) the originality of interpreting and synthesizing the Source of Cognition (Vedas, IV th century B.C.),
the Epic (with many stories), and the Upanishads (principles of Hindu philosophy) concluding in Non-Duality.
Then Special Duality (Visishta Advaita) follows, which asserts that Individual Being and Supreme Being are different in the beginning, but end to blend themselves (Rāmānujacharya, XI-th century).
And later, to see that the neutrosophic scheme perfectly functions, Duality (Dvaita) ensues, through whom the Individual Being and Supreme Being were differentiated (Madhvacharya, XIII-th - XIV-th centuries).
Thus: Non-Duality converged to Duality.
<Non-A> converges to <A>.

Know yourself to know the others.
Study the others to understand yourself.

In conclusion, I want to be what I don't want to be:
a Philosopher. That's why I am not.
    (That's why, maybe, am I?)

Control what you can, leave the rest to the luck.
Control what you cannot, free what you control.

We tried to de-formalize the Hilbert's for/mal/ization of geometry: by constructing an anti-model, which doesn't respect any of his 20 axioms! (F. Smarandache, <Paradoxist Mathematics>)
Because, by axiomatization, a theory loses its transcendental, myth, beauty, and becomes too arithmeticsized, technical, mechanical.
Or, if a system of axioms is defined in a theory, this should be of infinite (and, even better, of aleph-) cardinality.

Logicism:
Frege's axioms for set theory, to derive the whole arithmetic, were inconsistent (see Bertrand Russell's Paradox).

Look at these
*Inconsistent Systems of Axioms:*
Let $(a_1)$, $(a_2)$, ..., $(a_n)$, (b) be n+1 independent axioms,
with $n \geq 1$; and let (b') be another axiom contradictory to (b).
We construct a system of n+2 axioms:
   [I]    $(a_1)$, $(a_2)$, ..., $(a_n)$, (b), (b')
which is inconsistent. But this system may be split into two consistent systems of independent axioms
   [C]    $(a_1)$, $(a_2)$, ..., $(a_n)$, (b),
and
   [C']   $(a_1)$, $(a_2)$, ..., $(a_n)$, (b').
We also consider the partial system of independent axioms
   [P]    $(a_1)$, $(a_2)$, ..., $(a_n)$.
Developing [P], we find many propositions (theorems, lemmas, etc.)
          $(p_1)$, $(p_2)$, ..., $(p_m)$,
by logical combinations of its axioms.
Developing [C], we find all propositions of [P]
          $(p_1)$, $(p_2)$, ..., $(p_m)$,
resulted by logical combinations of $(a_1)$, $(a_2)$, ..., $(a_n)$, moreover other propositions
          $(r_1)$, $(r_2)$, ..., $(r_t)$,
resulted by logical combinations of (b) with any of $(a_1)$, $(a_2)$, ..., $(a_n)$.
Similarly for [C'], we find the propositions of [P]
          $(p_1)$, $(p_2)$, ..., $(p_m)$,
moreover other propositions
          $(r_1')$, $(r_2')$, ..., $(r_t')$,
resulted by logical combinations of (b') with any of $(a_1)$, $(a_2)$, ..., $(a_n)$,
where $(r_1')$ is an axiom contradictory to $(r_1)$, and so on.
Now, developing [I], we'll find all the previous resulted propositions:
          $(p_1)$,  $(p_2)$, ..., $(p_m)$,

$(r_1),\ (r_2),\ ...,(r_t),$
$(r_1'),\ (r_2'),\ ...,(r_t').$

Therefore, [I] is equivalent to [C] reunited to [C'].

From one pair of contradictory propositions $\{(b),(b')\}$ in its beginning, [I] adds t more such pairs, where $t \geq 1$, $\{(r_1),(r_1')\},\ ...,\{(r_t),(r_t')\}$, after a complete step.

The further we go, the more pairs of contradictory propositions are accumulating in [I].

*Contradictory Theory:*

Why do people avoid thinking about a contradictory theory ?

As you know, nature is not perfect:
  and opposite phenomena occur together,
  and opposite ideas are simultaneously asserted and, ironically,
proved that both of them are true!  How is that possible? ...

A statement may be true in a referential system, but false in another one.  The truth is subjective.  The proof is relative.

(In philosophy there is a theory:

that "knowledge is relative to the mind, or things can be known only through their effects on the mind, and consequently there can be no knowledge of reality as it is in itself", called "the Relativity of Knowledge";

  <Webster's New World Dictionary of American English>, Third College Edition, Cleveland & New York, Simon & Schuster Inc., Editors: Victoria Neufeldt, David B. Guralnik, p. 1133, 1988.)

You know?... Sometimes is good to be wrong!

How to reduce to absurd the *reductio ad absurdum* method?

*Continuum Hypothesis* (that the cardinality of the continuum is the smallest non-denumerable cardinal) has been shown to be undecidable, in that both it and its negation are consistent with the standard axioms of set theory.

By contrast to the relativism, which asserts that there is no absolute knowledge, in neutrosophy it is possible to attaint in pure science and by convention the absolute truth, t=100, and yet as a matter of rare fact.

*Hermeneutics of Hermeneutics:*

An idea <A>, by interpretation, is generalized, is particularized, is commented, is filtered, eventually distorted to <A1> different from <A>, to <A2> different from <A>, and so on.

Everybody understands what he wants, according to his level of knowledge, his soul, and his interest.

<A> is viewed as <Non-A> and even <Anti-A> at some degree (ill-defined).

But all deformed versions of this idea syncretize in an <A> way.

Idealists were so formal, empiricists so informal.
Neutrosophy is both.

*Sociological Theory:*

As in the Primitive Society, the modern society is making for MATRIARCHATE -
the woman leads in the industrialized societies.
From an authoritarian PATRIARCHATE in the Slavery and Feudalism towards a more
democratic MATRIARCHATE at present.
The sexuality plays an immense role in the manipulation of men by women, because the
women "have monopoly of the sex" as was justifying to me an American friend kept by
his wife henpecked!
A cyclic social development.
The woman becomes the center of the society's cell, the family.
The sexual pleasure influences different life circles, from the low class people to the
leading spheres.  Freud was right...
One uses women in espionage, in influencing politicians' decisions, in attracting
businessmen - by their feminine charms, which obtain faster results than their male
proponents.
The women have more rights than men in western societies (in divorce trials).

*Social Three-Quarters Paradox:*

In a democracy should the nondemocratic ideas be allowed?
a) If no, i.e. other ideas are not allowed - even those nondemocratic -, then one has not
a democracy, because the freedom of speech is restricted.
b) If yes, i.e. the nondemocratic ideas are allowed, then one ends up to a
nondemocracy (because the non-democratic ideas overthrow the democracy as, for
example, it happened in Nazi Germany, in totalitarian countries, etc.).

*The Sets' Paradox:*

The notion of "set of all sets", introduced by Georg Cantor, does not exist.
Let all sets be noted by $\{S_a\}_a$, where a indexes them.
But the set of all sets is itself another set, say $T_1$;
and then one constructs again another "set of <all sets>", but <all sets> are this time $\{S_a\}$
and $T_1$, and then the "set of all sets" is now $T_2$, different from $T_1$;
and so on... .
Even the notion of "all sets" can not exactly be defined (like the largest number of an
open interval, which doesn't exist), as one has just seeing above (we can construct a new
set as the "set of all sets") and reunites it to "all sets".

*A Paradoxist Psychological Complex* (with the accent on the first syllable):

A collection of fears stemming from previous unsuccessful experience or from
unconscious feelings that, wanting to do something <S>, the result would be <Anti-S>,
which give rise to feelings, attitudes, and ideas pushing the subject towards a deviation of
action <S> eventually towards an <Anti-S> action.
(From the positive and negative brain's electrical activities.)

For example:  A shy boy, attempting to invite a girl to dance, inhabits himself of fear she would turn him down...

How to manage this phobia?  To dote and anti-dote!
By transforming it into an opposite one, thinking differently,
and being fear in our mind that we would pass our expectancies but we shouldn't.

People who do not try of fear not to be rejected: they lose by not competing!

*Auto-suggestion:*
If an army leaves for war with anxiety to lose, that army are half-defeated before starting the confrontation.

*Paradoxist Psychological Behavior:*
How can we explain contrary behaviors of a person:  in the same conditions, without any reason, cause?
Because our deep unconsciousness is formed of contraries.

*Ceaseless Anxiety:*
What you want is, normally, what you don't get.  And this is for eternity.  Like a chain...
Because, when you get it (if ever), something else will be your next desire.  Man can't live without a new hope.

*Inverse Desire:*
The wish to purposely have bad luck, to suffer, to be pessimistic as stimulating factors for more and better creation or work.
(Applies to some artists, poets, painters, sculptors, spiritualists.)

*My Syndrome*:
Is characterized by nose frequently bleeding under stress, fear, restlessness, tiredness, nervousness, prolonged unhappiness.  This is the way the organism discharges, thus re-establishing the equilibrium, and it is fortunate because the hemorrhage is not interior which would cause death of patient.

The bleeding is cause by the nervous system, not by physical injury.
If you have any idea of treating it, don't hesitate to contact the author.  All opinions are welcome.

All is possible, the impossible too!
Is this an optimistic or pessimistic paradox?
a) It is an optimistic paradox, because shows that all is possible.
b) It is a pessimistic paradox, because shows that the impossible is possible.

*Mathematician's Paradox:*
Let M be a mathematician who may not be characterized by his mathematical work.
a) To be a mathematician, M should have some mathematical work done, therefore M should be characterized by that work.

b) The reverse judgement:  if M may not be characterized by his mathematical work, then M is not a mathematician.

*Divine Three-Quarters Paradox (I):*
Can God commit suicide?
If God cannot, then it appears that there is something God cannot do, therefore God is not omnipotent.
If God can commit suicide, then God dies - because He has to prove it, therefore God is not immortal.

*Divine Three-Quarters Paradox (II):*
Can God be atheist, governed by scientific laws?
If God can be atheist, then God doesn't believe in Himself, therefore why should we believe in Him?
If God cannot, then again He's not omnipotent.

Divine Three-Quarters Paradox (III):

`Can God do bad things?`

  a) If He can not, then He is not omnipotent, therefore he is not God.

  b) If He can, again He's not God, because He doesn't suppose to do bad things.

Now, even if He only can - without doing it -, means He's thinking to be able to do bad things, thought that again is not compatible with a God.

Divine Three-Quarters Paradox (IV):

Can God create a man who is stronger than him?

  a)    If not, then God is not omnipotent, therefore he is not God.

  b)    If yes, then God will not be the strongest one and God might be overthrown.

God is egocentric because he didn't create beings stronger than Him.

Divine Three-Quarters Paradox (V):

Can God transform Himself in his opposite, the Devil?

  a)  If not, then God is not omnipotent, therefore He is not God.
  b)  If yes, then God will is not God anymore.

  [Religion is full of god-ism and evil-ism.]

God and Evil in the same Being.
Man is a bearer of good and bad simultaneously.  Man is enemy to himself.  God and
Magog!

*Expect the Unexpected:*
If we expect someone to do the unexpected, then:
- is it possible for him to do the unexpected?
- is it possible for him to do the expected?
If he does the unexpected, then that's what we expected.
If he doesn't do the expected, then he did the unexpected.

*The Ultimate Paradox:*
Living is the process of dying.
Reciprocally:  Death of one is the process of somebody else's life [an animal eating
another one].

Exercises for readers:
    If China and Japan are in the Far East, why from USA do we go west to get there?
    Are humans inhuman, because they committed genocides?

*The Invisible Paradoxes:*
Our visible world is composed of a totality of invisible particles.
Things with mass result from atoms with quasi-null mass.
Infinity is formed of finite part(icle)s.
Look at these Sorites Paradoxes (associated with Eubulides of Miletus (fourth century
B.C.):
    a) An invisible particle does not form a visible object, nor do two invisible particles,
three invisible particles, etc.
However, at some point, the collection of invisible particles becomes large enough to
form a visible object, but there is apparently no definite point where this occurs.
    b) A similar paradox is developed in an opposite direction.
It is always possible to remove an atom from an object in such a way that what is left is
still a visible object.  However, repeating and repeating this process, at some point, the
visible object is decomposed so that the left part becomes invisible, but there is no
definite point where this occurs.
Between <A> and <Non-A> there is no clear distinction, no exact frontier.  Where does
<A> really end and <Non-A> begin?  We extend Zadeh's fuzzy set term to fuzzy concept.

*Uncertainty Paradox:* Large matter, which is under the 'determinist principle', is
formed by a totality of elementary particles, which are under Heisenberg's 'indeterminacy
principle'.

*Unstable Paradox:* Stable matter is formed by unstable elementary particles
(elementary particles decay when free).

*Short Time Living Paradox:* Long time living matter is formed by very short time living elementary particles.

*Paradoxist Existentialism:*
life's value consists in its lack of value;
life's sense consists in its lack of sense.

*Semantic Paradox (I)*:  I AM WHO I AM NOT.
   If I am not Socrates, and since I am who I am not, it results that I am Socrates.
   If I am Socrates, and since I am who I am not, it results that I am not Socrates.
Generally speaking: "I am X" if and only if "I am not X".
Who am I?
   In a similar pattern one constructs the paradoxes:
   I AM MYSELF WHEN I AM NOT MYSELF.
   I EXIST WHEN I DON'T EXIST.
And, for the most part:
   I {verb} WHEN I DON'T {verb}.
(F. Smarandache, "Linguistic Paradoxes")

   What is a dogma?
An idea that makes you have no other idea.
How can we get rid of such authoritative tenet?  [To un-read and un-study it!]

*Semantic Paradox (II)*:  I DON'T THINK.
   This can not be true for, in order to even write this sentence, I needed to think (otherwise I was writing with mistakes, or was not writing it at all).
Whence "I don't think" is false, which means "I think".

   Unsolved Mysteries:
a) Is it true that for each question there is at least an answer?
b) Is any statement the result of a question?
c)  Let P(n) be the following assertion:
"If S(n) is true, then S(n+1) is false", where S(n) is a sentence relating on parameter n.
Can we prove by mathematical induction that P(n) is true?
d)  "<A> is true if and only if <A> is false".
Is this true or false?
e) How can this assertion "Living without living" be true?
Find a context.  Explain.

   <Anti-A> of <A>.
Anti-literature of literature.
   <Non-A> of <A>.
Language of non-language.

   <A> of <Non-A>.
Artistic of the non-artistic.

*Tautologies*:
I want because I want.  (showing will, ambition)
<A> because of <A>.
(F. Smarandache, "Linguistic Tautologies")

Our axiom is to break down all axioms.

Be patient without patience.

The non-existence exists.
The culture exists through its non-existence.
Our culture is our lack of culture.

Style without style.

The rule we apply:  it is no rule.

*Paradox of the Paradoxes:*
Is "This is a paradox" a paradox?
I mean is it true or false?

To speak without speaking.  Without words  (body language).
To communicate without communicating.
To do the un-do.

To know nothing about everything, and everything about nothing.

I do only what I can't!
If I can't do something, of course "I can do" is false.
And, if I can do, it's also false because I can do only what I am not able to do.

I cannot for I can.

Paradoxal sleep, from a French "Larousse" dictionary (1989), is a phase of the sleep
when the dreams occur.
Sleep, sleep, but why paradoxal?
How do the dreams put up with reality?

Is O. J. Simpson's crime trial an example of:  justice of injustice, or injustice of
justice?
However, his famous release is a victory against the system!

Corrupt the incorruptible!

Everything, which is not paradoxist, is however paradoxist.

This is the Great Universal Paradox.
A superparadox;
(as a superman in a hyperspace).

Facts exist in isolation from other facts (= the analytic philosophy),
and in connection as well with each other (= Whitehead's and Bergson's thoughts).
The neutrosophic philosophy unifies contradictory and noncontradictory ideas in any human field.

The antagonism doesn't exist.
Or, if the antagonism does exist, this becomes (by neutrosophic view) a non-(or un-)antagonism:  a normal thought.  I don't worry about it as well as Wordsworth.

Platonism is the observable of unobservable, the thought of the non-thought.

The essence of a thing may never be reached.  It is a symbol, a pure and abstract and absolute notion.

An action may be considered G% good (or right) and B% bad (or wrong), where G, B $\subset \parallel^-0, 1^+\parallel$ - the remainder being indeterminacy, not only <good> or only <bad> - with rare exceptions, if its
consequence is G% happiness (pleasure).
In this case the action is G%-useful (in a semi-utilitarian way).
Utilitarianism shouldn't work with absolute values only!

Verification has a pluri-sense because we have to demonstrate or prove that something is T% true, and F% false, where T, F $\subset \parallel^-0, 1^+\parallel$ and n_sup $\leq 2^+$, not only T = {0} or {1} – which occurs in rare/absolute exceptions, by means of formal rules of reasoning of this neutrosophic philosophy.

The logical cogitation's structure is discordant.
Scientism and Empiricism are strongly related.  They can't run one without other, because one exists in order to complement the other and to differentiate it from its opponent.
PLUS doesn't work without MINUS, and both of them supported by ZERO.  They all are cross-penetrating sometimes up to confusion.
The non-understandable is understandable.
If vices wouldn't exist, the virtues will not be seen (T. Muşatescu).
Any new born theory (notion, term, event, phenomenon) automatically generates its non-theory - not necessarily anti-(notion, term, event, phenomenon).  Generally speaking, for any
<A> a <Non-A> (not necessarily <Anti-A>) will exist for compensation.
The neutrosophy is a theory of theories, because at any moment new ideas and conceptions are appearing and implicitly their negative and neutral senses are highlighted.
Connections & InterConnections...

The non-important is important, because the first one is second one's shadow that makes it grow its value.
The important things would not be so without any unimportant comparison.

The neutrosophic philosophy accepts a priori & a posteriori any philosophical idea, but associates it with adverse and neutral ones, as a *summum*.
This is to be neutrosophic without being!
Its schemes are related to the neutrality of everything.

Spencer's "organicism", which states that social evolution is from simple to complex and from homogeneous to heterogeneous, can be updated to a cyclic movement:
- from simple to complex and back to simple - since any complex thing after a while becomes simple - (but to a superior level),
and again to complex (but also to a superior level to the previous one);
therefore: from level 1 to level 2, and so on...
- idem from homogeneous to heterogeneous (level 1) and back to homogeneous (level 2), and again to heterogeneous (level 3)...
[a *neutrosophic evolutionism*, neither H. Spencer's, nor V. Conta's].

This neutrosophy creates anti-philosophy.
And, in its turn, the anti-philosophy creates philosophy.
A VICIOUS CIRCLE.
Both of them are making history (?)
It raises a notion/idea/event/phenomenon <A> to <Non-A>, and vice versa.

Philosophy is a poetical science and a scientific poetry.

There are three main types of humans: not only Nietzsche's "overman" with his will to power, but also the <midman> with his will to mediocrity (yes, people who love to anonymously live every single day, dull),
and <underman> with his will to weakness (homeless, tramps, criminals who indulge in laziness, illegalities).
Inside of every man there are an <overman>, a <midman>, and an <underman> - varying in terms of moment, space, context.
That's why, generally speaking, every man is: O% overman, M% midman, and U% underman, where O, M, U $\subset \rbrack\rbrack^{-}0, 1^{+}\lbrack\lbrack$.

While Spencer mechanically supported the flat evolutionism, S. Alexander, C. L. Morgan and later W. P. Montague focused on emergent evolution: the new qualities spontaneously and incalculablely emerge.
There however is a flatness within spontaneousness.

Lenin's "things' dialectic creates the ideas' dialectic, but not reciprocally" still works vice versa.

Same back and forth dynamics for trialectic (with neuters' attributes), pluralectic, transalectic.

If you learn better a discipline, you'll learn 'less' another one (for you don't have time to deepen the knowledge of the second one).
And, if you learn better a discipline, you'll learn 'better' another one (because the more knowledge you have, the more you understand another discipline). N'est-ce pas?

When unemployment U(t) increases, child abuse CA(t) also increases:
$$CA(t) = k \log U(t),$$
where t is the time variable, k is a constant depending on unemployment rate and children's percent in the population.

Philosophy is an ideational puzzle and, alike geometry, it circumscribes and inscribes an idea to and into a class of things.

To think means to be unusual and intriguing and uneasy to others.
If X says <A>, let's examine all its versions <$A_i$>, then what's <Neut-A>, afterwards focus on <Anti-A>, and don't forget all their derivatives. Let's question any and all. Let's be skeptical versus any "great" thinker.
Go ahead and look for the conflict of theories - grain of wisdom and creativity.

I see the ideas. They are red and blue and white, round and sharp, small and big and middle size.
I look through the objects and see the essence.

What could be a philosophical algebra? But a philosophical vector space? And how should we introduce a philosophical norm on it?

Wittgenstein's logical structure of language risks to get out of the main picture when passing from a language $L_1$ to a grammatically very different $L_2$.

Neither interdisciplinarity, nor multidisciplinarity, but the notion to be extended to *infinitdisciplinarity* (or total-disciplinarity), in order to form a *global discipline* - emerged from every single discipline in order to form an all-comprehensive theory applicable back to its elements, which gave birth to it.

Thomas Kuhn's paradigm is based on scientific and metaphysical beliefs as well.

Schopenhauer was radically pessimistic, what about a laughing philosopher? Would s/he be passed as a joker?

Determinist theory asserts that every fact or event in the universe is determined or caused by previous facts or events.
Ok, but what about the 'first' fact or event? Who did cause it?
If you are religious, you may answer: God. Then, who caused God?

Is the Supreme Being created by himself?  How is that?
Or, maybe there was no 'first' fact or event?  Then, how was it possible to get to some point with facts or events in space and time without a beginning?
Determinism flirts with underterminism at some degrees.

Every fact or event in the universe is d(F)% determined or caused of previous facts or events, 0 <= d(F) <= 100, and the percentage depends on each individual fact or event F.
The determinism partially works in this neutrosophy.
The proverb: he, that is born to be hanged, shall never be drowned, doesn't entirely apply.
The destiny is also deviated by man himself.

Our mind can not reflect truth accurately (Francis Bacon).
Unfortunately the science too.
What about arts?  (No, they are too subjective.)

"Truth is subjectivity" (Jaspers).
Yes, in most of the cases, but according to the previous definition, the truth may be objective too
(as a right limit of the subjectivity, when this is going far away from itself as $x \mapsto 1$).
The independent variable x swings between 0 and 1.
Subjective = 0, objective = 1, and everything else in between is a mixture of subjective and objective.  If the percent of subjectivity in a truth is s%, then its percent of objectivity is not necessarily <= (100-s)%.
The truth is not a stagnant property of ideas, said William James, ideas become true because they are made by events.  There are as many truths as concrete successful actions.

"Subjective" is, in its turn, objective too.
Objective is subjective as well.

No assertion is immune revision (W. V. O. Quine).

We extend the solipsism, theory that source of all knowledge of existence is self alone,
to pluripsism, theory that source of all knowledge are all beings, because we get influenced by others' believes, hopes, desires, fears.  It's impossible to isolately live, not even hermits or monks stay alone but they at least interfere with nature.  They have to - in order to survive.
We may never adequately understand our colleagues' experience (Thomas Nagel's empathic solipsism) or ascriptions of psychological states (Wittgenstein's psychological solipsism),
and yet a small percent of it we do understand, even if we misunderstand but we charge our unconsciousness with fragments of their thought - and later we may partially act in their way without even knowing!
We behave in certain way not only because of what occurs inside of our brain (as mythological solipsism asserts), but mobilized of external events as well.

A mathematization of philosophical (and not only) cognition is demanded.

Sometimes people don't even know why they reacted in the way they did. Something it came from their innermost depths, unconscious, something they were not aware of.

"Impossible de penser que <penser> soit une activité sérieuse" (Fr.)
[It's impossible to think that <to think> is a serious activity] (Emile Cioran).

And a sage: There is no philosophy, there are only philosophers. Therefore philosophers without philosophy!
But reciprocally: is there a philosophy without philosophers?

"Any big philosophy ends up into a platitude" (Constantin Noica).

The worth of an action is determined by its conformity to given binding rules (deontology), and equally by its consequences.

The same sentence is true in a reference system, and false in another one. For example: "It rains" can be true today, but false tomorrow; or can be true here, but false there.
Moreover, the sentence is also indeterminate: we don't know if ten years from today it will be raining or not.

Because any attempt to change the political power ends up in embarking another power, "the revolution is impossible" (Bernard-Henry Lévy, André Glucksmann, Jean-Marie Benoist,
Philippe Némo who represented the "New Philosophy" French group).

The power of the monarch derives from his powerless people (Juan de Mariana, 16-17th centuries). Because, if they had any power, monarch's position would be in danger.

"It looks like the great systems started to lose their influence, because they vainly slide over the universe" (Ţuţea).

Plurality of causes of a single effect (J. S. Mill), is extended to plurality of interweaved causes of a plurality of interweaved effects. It is impossible to separate the causes,

$$C_1, C_2, ..., C_n, \text{ where } 0 \le C_i \le 1 \text{ for each index i, and } \sum_{i=1}^{n} C_i = 1,$$

they act as a whole, and so the effects,

$$E_1, E_2, ..., E_m, \text{ where } 0 \le E_j \le 1 \text{ for each index j, and } \sum_{j=1}^{m} E_j = 1,$$

even more: both, the causes and effects, have the power of continuum.

An analysis and synthesis of the whole philosophy done by the neutrosophy would catch up with a self-analysis and self-synthesis (reflexivity), for the movement is itself a part of philosophy.  How should, by consequence, the neutrosophy of the neutrosophy look like?

In *cooperative learning* the groups of students should be heterogeneous (not homogeneous) with respect to gender, ability, and ethnic or cultural background in order to learn from each other and better interact.
Interdependence play un important role, because a student could have to cooperate with another one he might not like.
[ Reynolds, Barbara E., Hagelgans, Nancy L., Schwingendorf, Keith E., Vidakovic, Draga, Dubinsky, Ed, Shahin, Mazen, Wimbish, G. Joseph, Jr., "A Practical Guide to Cooperative Learning in Collegiate Mathematics", Mathematical Association of America, Washington, D.C., 1995. ]

The more things are changing, the more they stay the same.

All mathematical objects are manifolds (not functions, as Alonzo Church asserted).

The Eleatic School holds that <all is one>, and does not accept change and plurality. We say that <one is all> either, and unchanged and singularity don't work in the real life.

There is no real "ism", because "ism" reduces everything to a conceptualization, the thing-in-itself, a manifold of appearances -
while all is mixed and interdependable.

Time is fluid, visible, and material.  Like an organism, a being.  We are part of it.

Husserl's phenomenological epoch, is commuting not only from natural beliefs to an intellectual reflexion, but backwards as well.  It is passing through a neutrality midpoint zero from one extreme to another - besides intermediate multipoints.

Inside the atom protons+electrons+neutrons co\habit.

Theologians have defined the *trinitarism* as:  Father, Son, and Holly Spirit.
What about Devil?  Therefore, a *tetranitarism*?
But Angels?  Thus, *plurinitarism*?

It is necessary to introduce a measurement for the ideas' field.
Let's denote by "*IDON*" [Latin < *idoneus*, (cap)able of] the smallest unit to measure an idea.
    The idon-ical measurement is
directly proportional with the following characteristics of an

idea:
- novelty
- quality
- originality
- density
- continuity
- brightness
- quantity
- analysis
- synthesis
- truth-value,
and inversely proportional with:
- vagueness
- discontinuity
- triviality
- falsity-value.

"From error to paradox it's often not more than a step, but this step is definitive, because, contradicting even the apodictic character of mathematical assertions, it can become itself a knowledge river of future mathematics."
(Al. Froda, <Eroare şi paradox în matematică>)

Therefore mathematics is not sufficient to explain everything.  Science is actually limited too.

This is the *Ultimate Idea*:
there is no ultimate idea!
Leucippus's atomism, elaborated by Democritus, asserting that atoms and void are ultimate realities, is itself voided!

Any system or substance has a degree of disorder (measured by the entropy), a degree of order,
and a degree of order and disorder in the same time.

What is the sense of emptiness (Gabriel Marcel), but of wholeness?  They are opposite, but in my mind both look as perfect spheres.  Even the wholeness is vacuumed of sense.
In the pure form they do not exist.

We can treat various themes, that's why neutrosophy is not a specialized philosophy.  However it is specialized by its method of research, and by its system.  This thinking movement quotes the life.
It is applied in literature, arts, theater, and science as well.

We sometimes love our poetry when we don't love it!
The universe is expanding, the neutrality is expanding – for balancing.

To an event the paradox gives beauty and mister.

World is composed of contradictions.
Anti-world is composed of contradictions.
Contradictions are composed of contradictions.

World is of material and psychic natures simultaneously.
They may not be separated, as the materialist and idealist philosophers tried to do;
and not only the psychic is the superior result of the material,
but the reciprocal sentence as well.

Determinism means paradoxist causes.

Truth is relative [V. Conta], false is relative too. Both of them are cor(e)-related to a parametric (time, space, motion) system.

Crisis implies a progress. The progress, in its turn, unquestionably leads to crisis.
One knows the progress if and only if one knows the crisis (a kind of Upside-Down Way:
"Via Negativa" of St. Thomas Aquinas).
The development has valleys and hills.

Entities are tight by their differences too.

Paradox is infinite. This is a kind of God for the man.

We can paraphrase Hegel by:
what is rational is antagonistic, and what is antagonistic is rational,
and further:
what is irrational is antagonistic too.
You can say that 1+1=2 is rational, but not antagonistic.
However 1+1 may be equal to 3 in another logical system invented by yourself.
Nothing will exist and last out of neutrosophy.

Neutrosophy is not associated with Fichte's and Schelling's German idealism. For
example, in the absolutehood, categories such as:
     cause and effect,
     existence and negation,
may be reversed and mixed.
     There is a neutrality within each neutrality.

The famous pantheism's formula of Spinoza, brought from Giordano Bruno, *Deus, sive Natura* (that is "God, or Nature", identification of God with Nature), is generalized
to <A>, or
<Non-A>, sameness (up to confusion) between attribute <A> and its contrast <Non-A>.
     Synonymity of antonyms,

antonymity of synonyms.

I think, therefore I am a neutrosopher (paraphrasing Descartes's formula of existence: *Cogito ergo sum*).

From Schopenhauer's words "nothing exists without a cause", one merges to the existence of more causes - not only one -, and at least two of them are contradictory to each other.

World as paradox
Schopenhaurer said "World is my idea" using <vorstellung> (Germ.) for <idea>, therefore material is immaterial (because 'idea' is 'immaterial').

Some clericalist are atheist on the account they transform the church into a business and the religion into a political propaganda.

"The contradictiousness is a component of individual's personality" (E. Simion).

"Antitheses are the life" (M. Eminescu).

Men are the same, but everybody is different.

"If ever I could have written a quarter of what I saw and felt, with what clarity I should have brought out all the contradictions of our social system" (J. J. Rousseau).

A paradoxical argument:
"Man is by nature good, and that only our institutions have made him bad" (J. J. Rousseau).

In the matter's structure there is always a union between continuous and discontinuous.

Nothing is non-contradictory.  All is "+", "-", and "0".
Even the exact mathematics.
This is a DHARMA for neutrosophy.

Art is a God for our soul.

"Men will always be what women chose to make them" (J. J. Rousseau).
Consequently, men will be what they maybe don't want to be!

Learning we become worst (*civilization paradox*): further of ourselves.
Rousseau attacked the arts, literature on account of corrupting the ethics and replacing the religion.  By modern fashions we don't differentiate each other,
but conform in speech, cloths, and attitudes;  and we appear what we are not!
People are the same, but... different.

His irony against politicians:
"the politicians of the ancient world were always talking about morals and virtue, ours speak on nothing else but commerce and money".

His attack against luxury:
"those artists and musicians pursuing luxury are lowering their genius to the mediocre level of their times".
Hence, any progress in arts, literature, and sciences lead to the society's decadence.

"Man is born free; and anywhere he is in chains" (J. J. Rousseau, <The Social Contract>).

Human being's existence in society is unnatural
(let's look how he is not alike):
    he acts how he has to act (not how he feels alike)
    he speaks how he has to speak
his personality is destroyed and he became anonymous
his existence is nonexistential
he feels himself foreign  (Heidegger).
Heidegger rejected the science.

I AM NOT IGNORANT THAT I AM IGNORANT, parodying Socrates.

I DO NOT DOUBT THAT I DOUBT.
My authority is not to have any authority at all, for I'm not a dictator.

Nothing from what belongs to us does really belong to us.

We know that we don't know all.

Eternity does not exist.  It is a poem.
Eternity is passing...
Eternity is a delusion of the spirit thirsty of absolute.
Not even the absolute inwardly or through oneself exists, but it has been invented by humankind as a goal of not being able to aim at.  For judging the ardent.
Nothing is perfect, nothing is permanent.
Any notion is sullied by opposite elements, the contrary's umbrage is imprinted on it.
An object is lighted by its shadow.

Philosophy is a useless futile science.  It feeds the blue song of idealists.  Every philosopher is an idealist, the materialists too.
Who is keeping both eyes wide open at p-u-r-e notions and concepts?
The science's conventionalism is sometimes exaggerated.
Philosophy is a taciturnity... and a concealment...

Humanity is progressing against humanity,
until its destruction.  Not only a material ruin,
but people are turned into flesh robots.

How one explains the more mass cultural accent in underdeveloped countries than in rich industrialized ones?

But a hoarfrost of culture still subsists, for example, in the American academic media (Dana Gioia).

The more technology extends, the less culture flourishes. A new event in culture does not differ much from the precedents - culture even repeats, comparatively to the science's exponential growth.

The ratio $\dfrac{\text{Culture}}{\text{Science}} \xrightarrow[t \longmapsto \infty]{} \zeta,$

where $\zeta$ is a small constant and t represents the time.

Fortunately, the science influences the culture as well (see futurism, cubism, abstractism, etc.).

There exists a confusion between culture and civilization.

Alfred Weber analyzed the relationship between the growth of knowledge (science, technology) and the culture (soul).

A question: is there a limit to the civilization's advancing behind whose it's not possible to pass?

Science expanded over the culture, strangled it, occupied its place in the society.

To ponder:
- over the particularity of the general,
or the generality of the particular
- over the complexity of simplicity,
or simplicity of the complex
- over the negative side of the positive,
and reciprocally.

Life is neutrosophic: crying today, laughing tomorrow, neither one after tomorrow...

They are so close that life became more neutrosophic, and the neutrosophy more lifer/alive.

People have neutrosophic behaviors:
    friends who change to enemies or to ignorants...
    rich who fall to poverty or to middle class...
Ideas are neutrosophic.

A sentence may be true:
*a priori* (no matter in what conditions),
or *a posteriori* (depending on certain conditions).

Also, the same sentence may be true at time $T_1$, ignored at time $T_2$, and false at time $T_3$, or may be true in space $S_1$, ignored in space $S_2$, and false in space $S_3$, and so on…

There is a distinction between Neutrosophic Philosophy and the Philosophy of Neutralities.
The first studies the contradictions and neutralities of various philosophical systems, methods, schools, thinkers.
The second seeks the neutralities and their implications in the life.

"The paradox invaded all activity's fields, all scientific and artistic disciplines.  It is not a marginal phenomenon anymore, but in the heart of the act and the human thought.  Outside the paradox we are not able to understand the world.  We have to learn to identify the paradox in its stages of an extraordinary diversity, to discover its functional mechanisms for incarcerating and controlling it, and possibly manipulating it in order not to be ourselves manipulated by this.  If not long ago the paradox was considered a symptom of a pathological state,
in the last decades it is more frequent an opposite facet of paradox: that of a healthy, normal state.
[Solomon Marcus, "Paradoxul", Ed. Albatros, Bucharest, 1984]

Anti-structure doesn't mean chaos.

Logic of the False, or Anti-Mathematics?

There exist (feminine) YN Energy - left channel, and (masculine) YANG Energy - right channel, for psychic or spiritual power.
First one is of desires'.
Second one is of projects'.
Both are of biological nature.
In certain forms of yoga seven chakras coexist in the human body, but they can't be traced out through physical, chemical, anatomical means.
Kundalini Energy (of divine nature) is the universal energy's projection in us.
Athman (individual inward) blends with Brahman (collective inward) in Indian philosophy.
Yogic meditation consists of purification of chakras and touch of without-thought status, bringing to Kundalini Energy's increase.

*Concreteness of the Abstractness:*
An abstract notion is defined by concrete elements,
and reciprocally.
The concrete objects have their abstract qualities.

Laromiguière called our senses:  machines of making abstractions.  *Mechanical Philosophy*?
Devices of producing presuppositions on running belt (computer programming) - futile philosophy.

A priori thought à la Kant is inlaid to imagination only, a kind of passing to the limit towards infinity.

In the kantian space the thought dresses the purity form, going far from reality, idealizing and self-idealizing.

Nature's essence doesn't have a homogeneous nor a pure aspect - or that's relying on what acceptance we take the terms.

Mathematics also works with approximations. But *exact approximations*.

And thus the perfection is a notion invented by human: an endeavor, a target never to be touched.

We always want what we don't have. Once we've got it, our interest in it is lost. But we tend toward something else.

Human is in a continuous DESIRE, continuous SEEKING, continuous DISSATISFACTION. And these are good, for they bring the progress. So, human is in stress, plugged in.

(In the sport competition an aphorism says that it's easier to conquer a world record than to keep it.)

A cause of all empires' decline (none of them lasted and will never last indefinitely) is the self-content of their leading part in the world, slowing thus down their creative and vigilance engines.

In a universe there are more (concentric or not) universes;
    in a space:  more spaces
    in a time:  more times
    in a move:  more moves
We meet, as such, within a system other systems; and so on...
       subuniverse
       subspace
       subtime
       submove
       subsystem
And these concentrations pass upward and downward away to the
(macro- and micro-) infinite levels.

Nietzsche:  "All is chaos",
but the chaos is organized, hair styled on the curlers of an uncombed head.

Truth is hidden in untruth either.

Theory of Happenings and Theory of Unhappenings of  phenomena correlate.

Consciousness of unconsciousness.

We do not only support the theory of contraries, vehicled by dualists, but merely generalize it as:

There are only contraries:  no phenomenon occurs without its "non" (not necessarily "anti"), without its negation and neutralities.  We mean:  an event and its non-event are born in the same time.

For each object there exists an anti-object and non-object.

The difference between <A> and <Anti-A> are sometimes more pithy:
female-male, minus-plus, etc.
or more diluted...  but that's another story.

I do philosophy just because I am not a philosopher, and am not interested in philosophy.  I waste some time reading and skimming through mind treatises.

Philosophy is useless.  It is a head\ache for individuals without head.
Philosophers are inutile scientists.  I am not a philosopher.  Am I utile?
If philosophy is inefficacious, let's do philosophy!
The best philosophy is the total lack of philosophy?
Because the non-philosophy is itself a philosophy.
What about pseudo-philosophy?
I didn't even want to become a philosopher in this mercantile society (for I would starve to death).  That's why I philosophize...  I try not to find a system.

Today's people are very pragmatic, they don't give a penny on my neutrosophic arguments, nor on your anti-neutrosophic ones!
Only for money they are caring...
The number of humanists, and especially their percentage in the population, is dramatically decreasing.

What is the use of the useless theory?

But, the deep face of the world, its inner motion, its pressure and depression are hidden to our senses.
And that's why the world is sometimes what it is not.

That's the crisis of the modern man's crisis!
Neutrosophic nature envelops everything.

It is easy "to write" philosophy.  But philosophy shouldn't play for a round game.
It's harder "to discover" philosophy, we mean to find laws applicable to large categories.
The impeccable philosophy would essentially comprise the ideational metabolism of the infinite sphere - to absorb the ray of unbounded archetypes.

Jacques Derida's ideology:  the death of all ideologies!

Philosophy is not a unitary theoretical generalized representation of the world (just to intersect with our concepts the A. Comte's positivism).

"Metaphysical sentences are neither true nor false, because they don't assert anything, they don't contain consciousness nor errors" (Rudolf Carnap).

Human is infinite.  We oppose to Jaspers's finitude of human.  Spirit is its unbounded border.

Experimental law of Murphy:
constants aren't, variables won't.

Try to save what can't be saved!

It's easy to forget something important, but it's harder to forget something not important!

Imaginary is more real than the reality.

All is hatred, even the love.

"Knowledge is power" (Francis Bacon), but knowledge brings weakness too (for example a cancerous who knows he's sick).
  Knowledge is power in science, research,
but may be fear, suffering, even suicide - as in case of the previous patient, for example.
Power in a direction signifies weakness in another direction and mediocrity in a third direction.  I believe that power and weakness and mediocrity combine up and down.

When you ask yourself:  Why do I exist?  What is my mission in this stupid world? and you pessimistically think as a Kierkegaard, and especially Schopenhaurer, or your heart is vibrating of Chopin's piano grave chords?

Neutrality is the measure unit of all things, paraphrasing Protagoras's famous adage "Pánton chrémátōn  ánthrōpos métron" (Gr.) (Human is all things' measure).
Why?  Because the contradiction and neutrality are the nature's essence.  And examples we may find anywhere.

A philosophical system is a dogma (Francis Bacon).  That's why I plead for a philosophical system without system.
Not quite analytic philosophy.

Congratulations for your failure!

If you are defeated, fight back.
If you win, fight back either.
Is there a better strategy?

Ah, if I would have a force to change what's unchangeable!

We are permanently moving towards a homogenization of heterogeneous, as Stefan Lupasco would say.
Fixed is the transformation only.

Logos is penetrated by *NonLogos*.

Attempting to free himself through arts, man enslaves himself to creation.

"Homo homini lupus" (Lat.) (human is wolf for the human),
that's why there is a *bellum omnium contra omnes* (Lat.) (the war of everybody versus everybody), as a "natural status" (Hobbes from Plaut).
And Spinoza oppositely with *homo homini deus* (Lat.) (human is God for human), while Feuerbach absolutized to: the God of human is human itself.

Is man a hu(e)-man?

Schleiermacher's personalism proposes that all social problems be solved by evasionism, by intercommunication with God,
or by withdrawal in own personal "dimensions".
Therefore, a kind of 'forget about', of solving a problem by properly neglecting it (ignorance).

"It looks like the grand narratives started to loose their influence, because they slide in a sterile manner over the universe" (P. Țuțea, <Philosophia Perennis>).

Theoretical Categories.
I don't believe there is an absolute beginning of things, nor an absolute ending.
There is no perfect phenomenon, but tending towards a moving goal as in parametrial mathematical analysis.
Nothing perpetual.

Any notion is sullied by untangent notions.

"Know yourself", says a Latin adage.
But it is impossible to penetrate the inner infinity. It is question of psychical and even philosophical approximations.
Many times we feel strangers to ourselves, acting against our thought or senses - as people we would disapprove.

Human being is un organized chaos, endowed with abyssal reason, limited senses, and unbounded irrationalism. All is of continuous and transcendental field. Nor even phenomena are totally derivated ones from others, and there is effect without cause because the irrational has its act empire.
Cantor's set theory solved the infinitude of the finite, and surprisingly the equipotence of unequal sets, in the way that one was finite (segment of line) and another one infinite (the whole line) = paradox's pick!

The real world is messy.  Many problems are ill posed.  In practice there is ugly mathematics.  Clean the awful data to see the beauty of the theorems.  Non-mathematicians crinkle into the problems.

In the philosophy of arts and literature:  a network of beautiful, well, true is replaced by the voluptuousness for ugly, bad, false...
misery of life since Zola, the appetite for scabrous, mould, rot (Baudelaire, Arghezi), injustice of powerful people against the powerless,
the wrong promulgated with façade of right,
and generally <Non-A> dressed in clothes of <A>.

We do not speak on politics, because "in politics we do not have to tell the truth" (Metternich), nor on history which is the "prostitute of politics" (Nicolae Iorga), but on the nationalism of those who pretend to be cosmopolitan.

Existence of absconded contradictions, therefore of a continuous instability in the moving essence of things and phenomena.
Heraclitus's vision of harmony and stability join somehow the absolute, perfect, infinite values liable to a theoretical ideal.

Of course, we can find a harmony in contradictions and a stability in the middle of an instability - dialectically tied.
As well as
       an absolute into absolute
       a perfection into perfection
   and an infinity into infinity
"We enter in the same waves, and we do not.  We are, and we are not" (Heraclitus).

"We die and we do not die; human is a mixture of animal and god; all look when fortuitous when necessary" (Petre Ţuţea).
Decoding the paradox hidden in the problems' core.
Style means "unity in diversity".  "Life can be framed in the form of an instable equilibrium".
With a precise imprecision.

"I know that I don't know" (Socrates).

Philosophy doesn't need philosophers, but thinkers.  The thinkers don't need philosophy.  Therefore, philosophy doesn't need philosophy!
Is this an anarchy?

Philosophy is neutrosophic, or is not at all.
While Platon, by his dialogues, understands that he doesn't solve anything, Kant believes he solves everything.
None of them is correct.

A vicious circle:
Vasile Pârvan:  the ethnical is point of departure, and the universal is point of arrival.
*Terminus a quo* and *terminus ad quem*.
And again one returns to Petre Ţuţea:  nation is the ultimate point of universal evolution.
[We, personally, don't think so!]

Heidegger:  to live absolutely dying every day (in order to get out from anonymity).

The paradox produces anxiety, dizziness (revolving gloomy thoughts), arguing in a circle, twisting your mind around!
A solved paradox loses its mystery and it's a paradox no more.

How can we interpret the biblical expression:  "Enthrall me, God, for I to be free" (Imitatio Christi)?

Liberty is a unruly demon from spirit; and dissatisfaction leads to the revolt of the liberty, until it gets to an equilibrium.
While Ţuţea has another opinion:  "Human's liberty is the divine part of him".  *Divina particula aurae* only?
Equilibrium is in a permanent unstable balance.
And disequilibrium with propensity towards equilibrium.
As the saying goes:  Oh, God, give the human what he doesn't have!  You zealously need something and, when you get it, it hackneys in your hand.
Plus tends towards minus.  Minus tends towards plus.  They run each other, as in a vicious circle passing through zero.
Negative and positive.
Heterogeneity is homogenized.  Homogeneity is not pure.
There are optimal points that social phenomena are converging towards, and act as curves with asymptotes.  More exactly, differential equations would simulate the soul.
Extremes touch each other, said Marx, actually blasphemed philosopher.  Without extremes the equilibrium would not exist.
Didn't one vehicle in the Middle Age a theory on the *double truth* (interpreted upon faith: *secundum fidem*, and upon reason:
*secundum rationem* respectively)?

Every human is his own slave and master.

The mother nature is reversible and irreversible.

According to Ţuţea:  Christ is the divined human, and the humanized divine.
He also characterizes Nae Ionescu as:  "the metaphysical meditation moved to the daily level, or the raising of daily to the philosophical level"!

Cultured philosophy and in-cultured ideology!

Are really there phenomena without history, things without history?

No, this notion of <history> is incorporated in the essence of essences.  Even things without history have their history.

     Learning teaches you what not to learn either.
Intelligence has prejudices, prejudices have a grain of intelligence too.

     Imitation has an original character.  And, in its turn,
originality is often imitative.
These are not simple puzzles, escapades.
"God is creator, man is imitator" (Țuțea), and not only,
because man created God in his (imagination) mind.

     "Idiot's function is positive, for without him we would understand neither the geniality nor the normality" (Țuțea).

     Neutrosophy became as a religion, a contemporary myth.
Trans-spiritual.  Trans-sensorial.
Contradictory and neutral laws, factors, principles, functions.

     The fantastic results from the real's side, as an excrescence.  Afterwards, the inverse cycle follows:  when the real (scientific/technical conception) is inspired from imaginary.

     Nicolae Iorga considered the idealist factors have determined those materialist in the human society's evolution.
Conversely it is still right.

     "I thought the truth is universal, continuous, eternal" (Mircea Eliade, <Oceanografie>).
Of course, it is not.

     "One can solace the man who suffers because of happiness" (M. Eliade, <Oceanografie>).
And one can solace a man who rejoices at trouble.

     Man's endeavor to impossible, infinity, absolute passes through possible, finite, relative.
One explains <A> through <Non-A>.  Which means:  <A> is what it is not.

     Goethe-ian principle of bi-polarity:
idol and devil, interior powers of the human being that are in a permanent dispute.
Mephistopheles & Faust.
     While we plead for a pluri-polarity among various combinations of idol and devil in our soul and mind.

     Pure philosophical concepts are not to be found.  This is a dialectic of metaphysics, and similarly a metaphysics of dialectic.

Is there a necessity of happening and a happening of necessity?  We mean a determinism of in-determinism and in-determinism of determinism?
Is there an internal term of the essence of things which implies the appearance of an external term to them?
[Necessity = internal term;  happening = external term.]

A continuous discontinuity, and a discontinuous continuity in the process of evolution.  However the set of isolated points is of null measure.

Nothing belongs to us in this world.  Only our original ideas (if any!), transmitted to posterity, may bear our mind prints:
    a) spiritual ideas (such as theories, theorems, formulas, concepts);
    b) material ideas (embodied in art canvas, sculptures, architectures, machines, tools).
Creativity and inventiveness belong to us.

Philosophy will be neutrosophic, or will not be at all!
*Sine die.*
It is normal when a philosopher asserts something, another one
(to become conspicuous, to distinguish from the first) denies him, otherwise the second would be a simple imitator, an epigone.

And not only in philosophy.  Therefore, two opposite ideas/concepts/systems were set up.  Look how easy it is to develop the paradox.
Thus, it is normal to be abnormal! (Eugène Ionesco)
    The death of neutrosophic philosophy would signify the death of whole philosophy, and of humankind too.  (The philosophy of philosophy will reveal it.)  Because how would this look like to have all people thinking in chorus in unison all over the world?  Wouldn't it be a totalitarism?
The genius of philosophy shows this may not be absolute, perfect, finite.

There are two types of *totalitarism*:
a) *unconditioned* - of one's own will;
for example, the today's third world country people imitating/following the western ideology, politics, culture, behavior, etc.
b) *conditioned* - by military, ideological, economical imposed forces  (in dictatorships, for example see Arthur Koesler, <Le zéro et l'infini>).
Always in the world will be a totalitarism at some degree.
Individual is going with the crowd, without even realizing it (societal totalitarism against individual) - like a sheep with bent head in the herd.
Also, there is an ideational oppression of classics floating in the air, and the permanent revolt of contemporaneous.
    And totalitarism at transversal levels as well:  linguistic (dominant, so called "international", languages), politic (solidarity with the most powerful), economic, ideological, cultural, even scientific.

Gabriel Marcel wrote "Les hommes contre l'humain", speaking on brain-washing (in French: le lavage du cerveau), and on *tabula rasa*.
Mass-media partially does this.

Social disease, created by mass media's political manipulation:
Give citizens the impression/disillusion they are free,
and they'd feel they are -
even if they are not.
Give citizens the impression they live in a democratic society,
they'd feel they do -
even if they do not.
And absolute free society may not exist.  Countries differ by their degree of undemocracy.

Remain with your real world - which exists,
but I remain with my idealist world - which doesn't exist.
One exists through non-existence better.
At the beginning it was the end.  A realm passes.  A realm comes.  And, at the end, the beginning starts.
Let's present the actual phenomena as they are not.
The nonrepresentable represents something.
Let's define the human through a non-definition.
Rational being is full of irrational elements.
Man is a philosophical animal (but depraved, said Rousseau).
(Let's grade the degradation.)

Dante was Florentine (but not Smarandache)!
I am a model of unmodeled artist.  An anti-Goethe and non-Faust.
A sacerdotal sinner, a wicked saint.

The heroes hide cowardly secrets.  The poltroons have heroic facets.

Spirit and matter.
Spirit is an emanation of the matter, said the materialists.
Matter is an emanation of the spirit, said the idealists.  The truth is somewhere in the middle.  It is neutral.
Is the spirit material, and the matter spiritual?
Both, spirit and matter, have ambi-(even pluri-)valent characteristics.

Philosophy is an alive graveyard of dead ideas.

Soul is a kind of anti-body / anti-organism / non-body that isochronizes with the body through a unity of contraries and neutralities.  Soul is a part of the body,
body is a form of the soul.
Soul is the I and the non-I.

God is immortal.
But "God is dead", uttered Nietzsche.  That's why I believe in God.

How did Eugène Ionesco ejaculate in one of his dramas:  "The king is dead.  Long live the king!"

Perfection is imperfect.
This is a theoretical notion only, not touched in practice.

"The paradox is the limit up to where our mind can go,
besides of whom the nothingness shows up" (Țuțea, "321 memorable words").

Life is a source of joy and anger (completing Nietzsche, the poet).  Life is utile to the death.  Life is inutile.  Death is inutile too.  Then what?
We study the weakness of Nietzschian superman, his will of powerless.

Happiness is the headquarter of the future unhappiness.
The sin is the headquarter of further honesty.
Order is the headquarter of disorder.
Passion fights against passion.
Taste and disgust... to cut the Gordian knot.

Philosophy started when it didn't even start, and will end when it will never end.  This has been done when it was not done,
and it was not done when it was really done.

Where goes a road which doesn't go anywhere?
(Paul Claudel:  "Where goes a road which doesn't go to the church?")

The paradox is a therapeutical method in science.  Not speaking of arts and poetry, which hunt after it (see, for example, the Paradoxist Literary Movement set up in 1980's).
However, the science glowers at it!
James F. Peterman considered the whole philosophy as a therapy.

"Where are those who are not anymore?" (Nichifor Crainic).

What was I when I was not?  What was I before being?

My personal life became public (by printing my diary), my private life not private anymore.

"Poets' work can stay one near another, philosophers' not" (Schopenhauer).

The absurd is natural, so the un-natural.

[See the lack of sense of the sense.]

I write philosophy to denounce it, or to prove the sickness of philosophy (?)

How will the universe and humankind look like after one million years?
(This is not a science-fiction/fantastic question, but a more scientific problem.)  In what
direction will them converge?

My purpose is the infirmity of purpose!
Inward purpose is not a purpose.  Outward purpose is not a purpose either.

Any creed delivers an anti-creed.
<To have no creed at all> is also a creed, isn't it?

How to release the pain from the pain?  But the soul from the soul, and the body from
the body?
I want to be a measurer of the truth, to renounce to renunciation and get inspired by the
myths' charm.

Philosophy-poetry:
    an inspired non-inspiration
    a voluntary involuntariness
We need to artistically express the inexpressible.
And catch the non-artistic in an artistic form.

Atheism's role in faith's development.
Schleiermacher nominates by "God" the existence we relate on,
going up to a religion without a personal God.

Inward infinity of finite objects.

Beyond philosophy there is a philosophy.  Beyond arts there are arts.  Beyond religion
there is a religion.
The matter is of neutrosophic essence.

Philosophical poverty:  "We live together, but die alone" (Țuțea).

Man is the blossoming of the nature's neutrosophy.

Theology and science merge in philosophy.

From the animal psychology to the animal philosophy.

We always do things done by others.

Today's society creates underhumans, not superhumans (Nietzsche's übermensch), because man is lost, small, unimportant, forgotten in the huge amalgam of information, news at every second, scientific and cultural forces...  He doesn't face up with these accelerated dynamics.

The most complicated things are those easy.  The most uncommon ones are those common.
But we don't see them because we are superficial and don't have time to think deeper (collapsed by the aggressive day).
All is based and raised on contradictions and neutralities.
World is unitary in its variations and differentiations (Lossky, "World as an organic whole").

As in Ramayana epic, the neutrosophy adopts a skeptic attitude simultaneously rejecting and contradicting the famous philosophical theses.  In other words:  a LOKAYATA in contemporaries, or a CARVAKA.
And not disagreement in disagreement's behalf, but for generalization.  Didn't Voltaire say:  "The laws in arts are made to encroach upon"?

When the human being will understand what is not understandable?

Who made God?  Doesn't He, really, commit mistakes?  Doesn't He have His own God to hold Him responsible for His creation?  Or is he a dictator?!

"Tie two birds together.  They will not be able to fly, even though they now have four wings."  (Jalaludin Rumi) [<The Way of Sufi>, by Idries Shah]

Always what you don't have it's formidable, while what you have you get bored of it.

Man must live in accordance with the natural world around him (Pueblo Indian philosophy).
While genius should not!

*Credo quia absurdum* (Lat.) [I believe for it is absurd], credited to Tertulianus.
Therefore, I believe because it is unbelievable!

The idea of Kierkegaard's eternal alternative:  who emerges the man's impossibility to select or intercede among contraries.
A dialectic of neutrosophic states of ethical consciousness.

Normally the human rights are promulgated by those who do not respect them - according to the curious principle that making noise they pass unobserved.

It is the question of PHANTASÍA KATALEPTIKÉ (Gr.) [comprehensive representation] only by the contradiction law of component units.
*Philosophia perennis & paradoxae* (Lat.).

Do you still think at me when you don't think at me?

One reveals the non-real reality of philosophy. And the real non-reality as well.

As part of the general theory of efficient action (Kotarbiński's praxiology) the intermediaries' and extremes' roles must be caused.

Philosophy shows the human spirit's formation.
"Because a philosopher writes with a knowledge of what his predecessors have thought, his own work is at once a criticism of earlier thought and a creative contribution at the growing edge of philosophy" (Samuel Enoch Stumpf, <A History of Philosophy>).

"I am constrained to confess that there is nothing in what I formerly believe to be true which I cannot somehow doubt" (Descartes).

Theologian Thomas Aquinas agreed the universal is found in particular things (in re) and, according to our experience, it is abstracted from particular things (*post rem*).

God is the supreme nature. The divine reality inside trivial, and reciprocally. He is the supreme neutrosopher of all times. He is the absolute, the nothingness, the nonbeing, <A>, <Neut-A>, and <Anti-A> simultaneously.

Double Truth of Ockham:
    a kind of truth is the product of human reason,
    the other one is a matter of faith.

Seneca: "People love and hate their vices in the same time".
To love our enemies, to hate our friends? How unexpected we are!

How curious!
We strangely and peculiarly behave.

Platon said the soul is fighting between reason and passion.
Creator of the classic tragedy Pierre Corneille's characters are unwound between their ideal and their passion (<Le Cid>), but their ideal wins.
While Jean Racine's characters are destroyed by their passions (<Iphigénie>, <Phèdre>).

In our being there are an "I" and a "Non-I" that dispute the priority. It is that interior dissection which split our existence in two dual pieces.

"The scientific philosophy doesn't exist" (Nae Ionescu).
Philosophy is the road towards neutrality, the exercise on the border between being and nonbeing, an ideational reaction of the essential contradiction in the confrontation of YES with NO and thousands of intermediary positions in between.

Nae Ionescu tells the art work framed in a historical moment does not correspond in another moment.

Governmental investments do not bankrupt, even if they bankrupt
[because the government refinances them from people's tax money!].

I can't afford not to afford thinking.
My philosophy is to contradict the philosophy. And, thus, to deliver an Anti-Philosophy which, after a while, becomes philosophy.
I study others' opinions for I run counter to them.
My ideology means the death of other ideologies. I study Kant for not following him (because, if I don't study him and I know nothing on his <Criticism of Pure Reason>, I may accidentally rediscover his theory, but I would like to imitate nobody).

The neutrality constitutes the dominant note of existence,
such as mystery in the center of the speculative and metaphoric philosophical system of Lucian Blaga.
Its inner tension dilates it.
To reveal this is to form the future growth's stimuli.
This is regarded as an irrationalism of the rationalism, and reciprocally.

The *Paradoxism* studies the paradoxes and their use in different fields.
An axiomatic system of the paradoxism couldn't be other than...
contradictory. Theory of the sense and of the nonsense.
Form of the in-form.
See also the Paraconsistent Logic (Newton C. Da Costa, in the journal <Modern Logic>).

Paradox is metaphysically, unconsciously, occultly perceived... and resembled to the hell!
Absolute, abyss, perfect are only a few notions not touched by others than paradoxist senses.
They are isomorphic.

For any kind of opinion there are a counter-opinion and a neuter-opinion;
for Kant there is a counter-Kant and neuter-Kant,
for Moses ben Maimonides a counter- and neuter-Moses_ben_Maimonides,
for Augustine's philosophy a counter- and neuter-Augustine's_philosophy.
    Existence and counter-existence and neuter-existence.

Since philosophy was born - due to its mosaic of counter-set ideas, systems which clash, rival Schools - the neutrosophy also came to life. But people didn't realize it.
Neutrosophy exists in the history of each field of the cognition.
    *Displacement towards neutrality* - this is the motto of evolution.

Cognition rises from neutrality to neutrality.

Politics is dictated by mean interests (Machiavelli).
The Arabic philosopher Ibn-Haldun defines history as a repetition in a regular and alterant mode of the cycles of climbing and decline of civilizations.

You don't need to be a philosopher in order to become a philosopher.

Manichaean dualist religious doctrine of eternal fight between good and bad (or light and obscurity), originated by the Persian prophet Mani (Manichaeus) in the third century A.D., combining Zoroastrian + Gnostic and other elements, is among the first forms of pre-neutrality expression.

"You become what you are in the context of what others did from you" (Sartre).
Hence, you are what you are not.

A method "to make/produce philosophy":
- pass a strong basic idea <A> through all known philosophical systems, thinking schools,
and compare it with their opinions, concepts;
- extract the <PRO-A>, <CONTRA-A>, and <NEUTER-A> sentences, comment and argue with them.
Everything put in a form of short sections (analytic philosophy), systematically concatenated on themes, notions, categories. And, of course, using an adequate meta-language.

I am asking if the form may exist outside of matter?
Aristotle denied it.
But the thoughts, the ideas... do they have some form?

There are philosophers who contradict themselves: like me, for example. Only that I am not a philosopher(!)

Every phenomenon, action of ours, however much positive, has negative parts.
And, however worst, has good parts.

To win, we need first to lose.

People should speak philosophy.
They already speak philosophy - but don't realize it!
     people eat philosophy
     people drink philosophy every day

Philosophy should be a dream of contemporary citizen.
However, their philosophy is not to do any philosophy at all.
Their thought is not to think.

A famous poem of Tennyson:

"Theirs not to make reply
"Theirs not to reason why
"Theirs but to do and die.

    Criminals are transformed in heroes.
Sinners in saints.  This is the contemporary world!
While innocents and obedient become victims (the poorest) of the society...

    Exterior world is real, but dependent on our consciousness,
therefore not real!

    The lack of existence of the non-existentialism.
    The lack of absurd of the absurdism.

    Was the American Pragmatism (Charles S. Peirce, William James, John Dewey)
another kind of the(r)orism?
At Peirce we see thought [= theory] and action [= practice], then an alloyage.

    Any idea is tested by its neutrosophic effects.

    Philosophy is a speculation, starting from an easy idea, which gets bolder, extended,
and applied to available systems...
as a skeleton covered with an aesthetic skin, which forms a body.
    And such, the philosophy is not a speculation anymore.
The philosophy is still and is not.

    The spirit is transcendent.  The spirit is also material.

    If a philosopher <F> one day asserted an idea <A>, in the future another philosopher
<G> will neutralize him supporting/motivating the idea <Non-A>.
This is a way to do philosophy, or a philosophical career for some ones.

    As an attacker, there is no doubt that you need to defend your attack from the
opponents' resistance.
As a defender, it is not doubt that you need to attack the attackers?
    The best defense is the attack - says a proverb.

    You better like my poems when you criticize them.

    The cure is worse than the problem it supposes to treat.

    Simone de Beauvoir exists even when she doesn't exist
[by her literary work].

    Western culture is progressing in a wrong direction, towards European man's crisis
(Husserl, <Phenomenology>).

Wittgenstein:  "the results of philosophy are the uncovering of one or another piece of plain nonsense".
Interpretation of misinterpretation?

Human gets to identify with God, on the way of soul's liberation and of status of detachment from the world (abgeschlidenheit) [Meister Eckhart, <Die Deutsche Werke>].
But human gets to identify with Devil as well, by revealing the misery of soul and private life.

Essence is God (*essential est Deus*),
essence is Devil (*essential est Diabolus*) either.
Both, God and Devil, are necessary to keep an equilibrium.
God and Devil identify because they are abstract, symbolic, infinite, fuzzy even neutrosophic notions.  And, especially, because there is no pure "positive" or "negative" action.  Each action is a percentile combination of "+" and "-" and "0" attributes.  God, also, commits errors;  (the Bible is full of crimes, incests, and sins).  Devil, in his turn, does beneficial work (because this is like a vaccine, which helps our mind to produce immunity to Bad Behavior "disease" by causing the formation of spiritual "antibodies", which we would call "*antispirits*", produced by our brain).  From vice we again rise, on a long staircase, to virtue.  From virtue we decline back to vice (the opposites attract) - passing through neuter, because monotony is against our biological rhythm.  And the cycle is habitually rotated.
There is no God neither a Devil, but a mixture of them -
they neutralize themselves at some degree:
a "devilish god" and a "godly devil", we would call Him/It *DevGod*.

To most of the questions:
- there is no exact right answer
- there is no exact wrong answer,
or
- every answer is right
- every answer is wrong,
because it is an interpolation of them.

A formal system, interesting enough to formulate its own consistency, can prove its own consistency if and only if the system is inconsistent  (Gődel's Second Incompleteness Theorem).

Cultural events occur 'synchronically' in many countries, but 'protochronically' as well.  The first adverb includes a quantity of universal, the second a quantity of particular.

How can we melt abstractness with concreteness?

*Paradoxist Determinism*:
The lack of cause is, still, a cause.

This is a *Definitive Judgement*:
there is no definitive judgement.

*Le Roi le veut*.  Let's cite the masters:
Platon: *panta chorei* (all is moving);
Diogene Laertius:  *rhein ta hola* (all is passing);
Aristotle:  *panta rhei*, *ouden menei* (all is passing, nothing is remaining).
Therefore, a today's affirmative sentence will be infirm tomorrow.

I have decided not to decide anything anymore.

World continuously changes, and ideas alike.  But, after a while, this arrives in the
same position.

We can easily get from an extreme to the other.

The paradox is immanent to the consciousness (Schuppe),
whence the whole neutrosophy is an immanent philosophy
(because the paradox is a part of neutrosophy).

Lotze studied the distinction among reality, truth, value.
He initiated the axiology, the philosophy of the culture, the anthropological philosophy.
Let's analogously introduce:
> the *NEUTROSOLOGY* (philosophical significance of the neutrality *in lato sensu*),
> the *NEUTROSOLOGY OF THE CULTURE*,
> the *NEUTROSOPHIC ANTHROPOLOGY*.
And so on:  neutrosology of the values, histories, sciences, arts.

Philosophy reflects the existent from the non-existent.

Heraclitus found consensus of opposite propensities and tensions, as that of bow and
lyre.  People can't imagine how in-harmony-with-itself is the discord!

I thank God he told me he doesn't exist.
This is my *Te Deum laudatum*!

Wouldn't it be possible to set up in the calendar a religious holiday for atheists?

We should realize that sometimes the beautiness is ugliness, and ugliness is beautiness
- paraphrasing Gertrude Stein.

There is a unity between the scientific and artistic languages, and this is not Neurath's
*physicism*,

but an accommodation of variability.

Just because the man is mortal, he wants to become immortal
(by his creation in arts, science, history).
What would happen if all men were immortal?

Every man bears inside a supra-man (positive energy), an infra-man (negative energy),
and null-man (no energy side):
        himself projected outside of himself
        himself projected inside of himself
They are sporadically activated.

"Man is a upsurge towards it's-not-possible" (Ion Ornescu, <Poems from Prisons>).

Causes and effects are antagonistic.

There are no dynamics without antagonisms.
To the neutrology the to-(and for-)itself inner movement is characteristic.
{Introspection}
Behaviorism initiated by J. Watson can't be linear.  There also exists an inward
behaviorism of the being (esse), in a continuous disequilibrium with being's outward
image reflected by F. C.
Tolman (and, especially G. H. Mead by his "social behaviorism" concept).  Uniting the
previous inward and outward with *neut-ward* notions, we find a *Neutrosophic
Behaviorism*.
Human demeanor functions upon the laws of differential equations with partial
derivatives with respect to more or less adverse-to-each-other parameters.  Hence:
nonlinear.
Nature's essence, beginning with the atom, consists in the fight and agreement of its
components.  Their convergence passes through divergence.

I decide not to decide anything.

"The beginning and ending are to be found on the same circle"  (Heraclitus).

How infinitely big can be the infinity?
We never could imagine in our mind how infinity would look like in the daily life!
Question which fascinated us...
To go, to go, and never to get to an end?  Or, if you ever find an end, how is it?  A
hundred miles high wall... and thick?  A precipice, a chasm?  Or the universe is circular,
and we indefinitely turn and turn around?  The universe, as a sphere or closed surface,
has no beginning, no ending.
Nor the small infinity neither the big infinity we perceive other than in an abstract
theoretical mode.

The theory of transcendental infinity of Cantor contains the paradox's beauty.  Two unequal sets may have, however, the same one-to-one correspondence among their points.

This was the great surprise which disturbed his rival mathematician Kronecker.

But no one can pull out the charm and ineffable from the science world (and from the new truths which, as part of old reference systems, deny the superannuated classicized assertions).

Dictatorship:
even if you don't want, you have to!

May we construct a non-philosophical philosophy?

Any positive has its negative and null side effects.
Pace of mind does not exist.  The systems' war grinds and reborn our neurons.

How to dispute an indisputable subject!

Formation is just the tension between contraries.  It resorts to a kind of a catastrophic law, a fight (*polemos*)  (Radu Enescu, <Eminescu, the Chimerical Universe>).

Man has to be prepared to adapt, and to face controversial situations.
Man should be strong enough to afford the unaffordable.
"It is much better for people not to have accomplished all they want; it is the sickness which gives sense and worth to the health; evil of good; starvation of saturation; fatigue of leisure"  (Heraclitus).

John Dewey considers intelligence itself as a habit by which the man adjusts its relation to environment.  It's a permanent circuit:

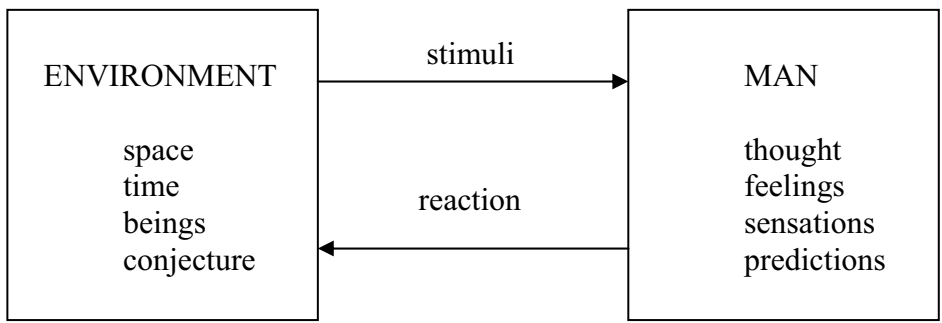

which stops when man dies (and even long after).

"Marx has often protested that he was not a Marxist", writes Samuel Enoch Stumpf in his <Socrates to Sartre> history of philosophy, 1988.

The neutrality knows a process of unbounded self-development towards an absolute (Hegelian) spirit.

According to Marx and Engels, everything (and the neutrality included) "is in a ceaseless state of movement and change".

The philosopher is not capable of discovering unique form of information (Wittgenstein).

But no use and (ab)use of neutrosophy in order to "draw false inferences, or ask spurious questions, or make nonsensical assumptions" (A. J. Ayer, especially referring to the analytic philosophy).

Neutrosophy is not only "unitary and fight of contradictions" (V. I. Lenin), but their transcendence in our every day's life;
        neutrosophic space, time, existence.
This is a generalization and relativization of ambiguities as well as philosophical controversies.

Non-philosophy makes philosophy.
Most philosophical problems arisen from puzzlement (L. Wittgenstein, Gilbert Ryle).

Pleasure and suffering, both, were studied by Heidegger and Sartre.

Sometimes an analytic method may be synthetic, while a synthetical method may be analytic.  W. V. O. Quine argued that "a boundary between analytic and synthetic statements simply has not been drawn".
Every notion has the shape of a *neutrosophic sphere*:
t% of the points/elements are surely in (inward),
f% of the points are surely out (outward),
and i% of the points are indeterminate (neutward),
where $t+i+f \leq 300^+$.

The paradox has many functions, besides its <classical> "*para*" (against) "*doxa*" (opinion) (Gr.) feature.  There is a plethora of oppugnancies inside of it.
But, don't let the paradox bewitch you!

Any dogma gives birth to an anti-dogma.
<Not to have any dogma at all> is also a kind of dogma, isn't it?
The anti-dogma comes back again as a dogma itself.

Body and mind are brought together (and studied by Gilbert Ryle).

*Neutrosophic Existentialism*:
Life is, now, machinized.  Machine is, now, humanized by science's sensorial improvements.

This is a nonexistential existence.
Human is dehumanized!  What is its alienation?

*Anti-Tautology*:
Existence is, therefore, what doesn't exist...

A person is governed by his neutrosophic senses.  I was surprised that people didn't grasp it!
They all shunned it, suffering from night blindness.
Without variations in opinions there would be no evolution.

There are many modes of neutrosophies, of course.  Like in Husserl's phenomenology, we need somewhat to stay back from the realm of experienced neutrosophies (in a detachment act) in order to understand and enable to master our life's opponents and neutralities.

This is a philosophical system without system, or based on non-system.

According to Kierkegaard the anxiety involves an antipathetic sympathy and sympathetic antipathy.

Physic joy may lead to a psychic bitterness.
Spiritual victory is conquered through bodily damage.

"Human existence" (Heidegger's *Dasein*), by its nonsense and absurdity, leads to Non-Existence?  (Somehow: self-destruction?)
Any evolution ends by closing the cycle (demise)!
Generally, the point of maximum extreme on the evolution's curve of any phenomenon is identical to a previous point to the phenomenon's origin.  Circular infinity coincides to zero.
The existent, in its apogee boiling, passes to non-existent. <E> is transformed in <Non-E> (not necessarily <Anti-A>), which is transformed in <F >, which is transformed in <Non-F>...

The victim loves his executioner.
The loser like his prejudicist.
The thrall adores his landlord.
The dog licks the whip which beats him.

How can we bring ourselves into agreement with somebody else's disagreement?
What about studying the informal formalists,
or the formal in-formalists?

Explanatory force of agent's no-reason-for-his-action,
or of agent's reason for his nonaction?

"A philosopher when he pretends to be a philosopher, he's not a philosopher. Philosophical ideas normally spring forth spontaneously, otherwise if you try to color 'em - they look stridently."  (O. Paler, <The Ten Commandments of 'Wisdom'>)

Assimilate what is too little (rarity is precious),
Dis-assimilate what is too much (a Romanian proverb says:  what is a lot, it's not worth).

*Negative definition.*
We introduce a concept <C> to the students by explaining  them what <C> is not.
We, thus, teach what <Anti-C>, the opposite of <C>, is all about to make them understand <C>.
This method is common in science when <Anti-C> is easier to define or better known.
In a similar way we may introduce a concept <C> teaching the students what is <Non-C>.

Analyzing, synthesizing, and evaluating opposite and neutral subject matter leads to neutrosophy.
Inter-(trans-)disciplinarity bases more on integrating theories from apparently intangible disciplines.
A such Centre Internationale de Recherches et Études Transdisciplinaires was set up in Paris, chaired by Basarab Nicolescu.

In regard to J. Piaget's & B. Inhelder's theory of cognitive development that individuals construct knowledge by interacting with their environment, we support the idea that a person's intellect is influenced by contradictory phenomena, facts, events on a background of neutral ones.  The more different they are, the better experience of life and mental growth.
Social interaction encounters pro-action and anti-action and neuter-action.
[<La Psychologie de l'Enfant>]

Theories are elaborated from facts, but facts from theory too.  If the university variant of relativism is the assertion that "there is no objective criterion to decide between true and false, good and bad, the farm hand variant alleges that all is a power game" (H. R. Patapievici, <The Relativism and the Politics / A Waisting Scoundrelism>).

"I find myself traveling towards my destiny in the belly of a paradox" (Thomas Merton, Trappist monk).

One changes from a common thing to a bizarre thing, then from a curiosity back to normal.

The "absolute value" (Platon, Aristotle) is displaced to a "relative value".
The pure ideas are generally impure.

In the sacred Hindu text "Bhagavad-Gita", found in the "Mahābhārata", one of the ancient Sanskrit epics, Lord Krishna lays the complete knowledge of life to his pupil Arjuna:
He who in action sees inaction and in inaction sees action is wise among men.
(Maharishi Mahesh Yogi, <Bhagavad-Gida: A New Translation and Commentary with Sanskrit Text>)

Ultimate order means chaos.

Maieutic Neutrosophy:
to get to the true by contradictory and neutral debates, conversations.

Miguel de Unamuno: When two folks Juan and Pedro talk,
there are six folks who actually talk:
- the real Juan, with the real Pedro;
- the Juan's image as seen by Pedro, with the Pedro's image as seen by Juan;
- the Juan's image as seen by himself, with the Pedro's image as seen by himself.
Actually, there are more:
- the Juan's images as seen by various people around, with the
Pedro's images as seen by various people around.
How many dialogues are taken place?
But in a group of n folks, when everybody talks?

We know without knowing.

In biology which one, the fixist theory or evolutionist theory of beings, is true?

In the modern diplomacy "saving time and energy is not possible without the replacement of real communication by a code, through formalization. (...) For the rest, the code remains almighty. You are <important> and null in the same time. More than yourself, you are as much as your badge - the little cardboard which marks your place at the debate table - allows you to be." (Andrei Pleşu, <Some Eastern Neuroses>)

If you seriously speak, you are laughing at me.
If you don't really speak seriously, you even more laughing at me!

To build a philosophy without any philosophical support (from the scratch)? Would it be a *"naive" philosophy*?
Philosophy without philosophy?
Paraphrasing Husserl: to judge only by comparison with the antinomies, and not according to any other trivial phenomena.

A neutrosophic phenomenology is based on intentional consciousness oriented towards the life's ups and downs and linear events. This is a branch of Husserl's Phenomenological Epoché.

We can easily get from an extreme to another,
but sometimes hardly between two close states.

The fear of ourselves...  We don't know who and why we are...

"every YES must to lean upon a NO (otherwise what Archimedes' lever would lean
upon?)"  (Ion Rotaru).

Philosophy is not solid.
Idea gives birth to non-idea (not necessarily anti-idea),
otherwise the previous would become an indoctrination.
New spirit builds on the old spirit by destroying it.
Another conventional logic replaces the superannuated logic.
Any assertion is a limitation, that's why a non-assertion comes regularly out:
for pushing the limits.

*Contemporary Neutrosophic Moral* issues.
There are arguments for, neuter, and against:
abortion, euthanasia, homosexuality, pornography, reverse discrimination, death penalty,
business ethics, sexual equality, legal use of drugs, economic justice.

The paradox is a mystery!
Gabriel Marcel's "What am I" particular human question has two complementary answers
within paradoxism:
a)  I am what I am not,
and
b)  I am not what I am.
We think these tell you everything.  Period!
Sometimes:
a)  It's possible to capture the impossible
and
b)  It's impossible to capture the possible.

*Ein Buch für Alle und Keinen* (Germ.) (a book for everybody and for nobody,
Nietzsche), subtitle to "Also sprach Zarathustra".

If you do something, it's wrong.
If you don't do, it's wrong either.
In conclusion:
        What should you do?
    and/or
        What shouldn't you do?

Man is infinite in his interior, and finite in his exterior.
How is it possible that a finite entity include an infinite one?

The benefaction of the neutrosophy emerges from its philosophy of life and writing:
 it's normal to have bad and good in life,
 it is even better than only bad or only good (which mean
monotony, whence death of mind and action).

 Your happiness is inside of yourself (from Buddhism).  Thus,
God is inside of man.
Your sadness too.
But man is inside of God as well.
And yet man and God do not coincide.

 Philosophy had to govern the state in the Athens democracy (Karl Popper), while in
the modern "democracy" of Hegelian inspiration the philosophy became the most slave of
the fishy demagogues.

 "All is a continuous metamorphosis, and so its contrary" (Chuang Tzu in his taoism -
School of the Way).

 There is an <A> beyond <A>.
Example:  There is a reality beyond reality; which one?  reality from our imagination.
 Neutrosophy doesn't consent in any way to the domination of some spiritual doctrines
-
although, in its turn, this becomes intrinsically established as another doctrine(!)...  for
and, at the same time, against all the doctrines, but keeping a neutral side..
Whence, neutrosophy will later act versus neutrosophy, giving birth to *post-neutrosophy*.

 *Labelism*.
 Trivial ideas of big guys are taken more important than clever ones of anonymous
individuals.  Everybody's judged upon his place in the society.  People have labels stuck
on their forehead.
Great ideas of poor persons, or from poor countries, are intentionally ignored.  Big guys'
mistakes are hushed up.  It is not the spiritual work which counts the most, but the
author's position (faculty of a "famous" university, his/her book or paper published by an
"important" publishing house or journal, connection network, scientific or artistic mafia,
arrangements, snoring awards).  Traffic of influence!
 The crowd is manipulated by mass-media, which became the strongest force in the
society.  People's consciousness is stolen.

 Daily citizen, in accordance with illuminist Rousseau, bears <mask>.  Due to the
sophisticated technology, he can't look inwardly, he swims through the world passing
besides himself.
Only spirit man is brave enough of his inward retrieval (La Rochefoucault), enduring that
"luxury humiliation".
 Technologized man who is doesn't answer what he feels, thinks, or is true; but what's
good for him to answer (in order to keep his social position/job or be promoted, or in
foresight of rewards).  He's robotized.  He's dis-humanized.  He's false...

The individual is surpassed by universal.

From experimental psychology to *experimental philosophy*.

"Geometry is exaggeration, philosophy is exaggeration, and so poetry.  Everything which has sense is exaggeration.
(...) Ontologically, the bad and idealization, as spirit necessity, are exaggerations.
(...) Greeks' <measure> was excessive as the *hybrid* which broke it.  Their serenity was a fickle equilibrium, of contrary excesses.
Without a dose of exaggeration there is not knowledge, nor action as well.  Neither science, nor justice.  And not even common sense."  (Alexandru Paleologu, <The Common Sense as Paradox>)

Plato:  essence precedes existence, which is easily explicable for objects.  You first think you need an apparatus - and what characteristics to have -, and second you build it.  Following Plato, what now exists, was necessary (produced by natural laws).
Anthropological question:  Thus, human being could be predicted from the origin of the solar system?
Sartre:  existence precedes essence, which is available for beings.  Say the horse, first exists, and then we study its characteristics, which are general for all individual of the same specimen.
Who is true?  (Both of them!)
Who is wrong?  (Both of them either!)
Then, which one came first, the egg or the hen (?)
There is a cycle:  existence ↦ essence ↦ existence ↦ ...
In our opinion none of "existence" or "essence" is first.

"Whatever can go right, won't!"

Nietzsche:  God is dead.
Dostoievsky:  If God did not exist, everything would be permitted.

Connections and adversities among ego and *superego* and *underego*.

*Daimon* is a form to illustrate the pulling off of the mobility  from immobility (Gabriel Liiceanu).

Man is free in society,
but governed by its laws.  Therefore, man is not free, but limited.  There is no absolute freedom.

"Humanity can't live just by logic.  It also needs poetry." (M. Gandhi)

Epistemology:
How to know all what we know?

*Spiritual pathology*:
the philosophy is my life's disease.

Immaterial matter?  Is that an absurdity?
Hobbes wrote about immaterial substances as being something meaningless.

We are not us;  we live through friendship's, profession's, language's, and epoch's
elements.  (C. Noica, <Book of Wisdom>)

Hobbies.
Ordinary people became the slaves of objects (luxury car, house), of passions (sex, trip).
Masters:  slaves of ideas.

When you want to make connections, you even tie opposite cases;
and when you don't want, you even separate identical things.

Be adaptable to inadaptability.

Change the change.

The happiness of the artist persists in his unhappiness.

The greatest moral lessons are thought by immoralists
(because they were landed into the trouble, and are experienced).

Hermeneutics of Foucault, hermeneutics of previous hermeneutics, until anti-
hermeneutics...
Wanting to support too much, you might deny!
From homo religious to homo neutrosophus.

A *pluri-philosophy* means opponent and similar ideas cross-referenced in the whole
cognition.

There are neither a definite end nor an ultimate purpose.
The theologians ignore the happening's role, and so the vitalists.

Researchers in the Chaos Theory are in progress to discover order inside of chaos of
the nonlinear differential equations.

Existence hasn't sense, and yet has a sense.

Will and non-will, goal and non-goal, sense and non-sense all act together, we being
conscious and unconscious of them.

Heraclitus:  All is changing.
Parmenides:  Nothing is changing.

Who is right?

    Heraclitus:  Individual is essential.

    Parmenides:  Universal is essential.

Who is right?

    Heraclitus:  pluralist.

    Parmenides:  monist.

Who is right?

(Both of them in each case!

And, simultaneously, both of them made mistakes.)

  Every reference system reflects a sentence in a different light.

    Empedocles explains how *Filia* and *Neikos* (Love and Hatred, attraction and abhorrence) function together.

    Prothagoras was the first to say that in all things there are contrary reasons.

    Gorgias's definition of rhetoric:  "the art to transform the worst thesis into the best thesis" (or *Ton eto logon kreito poeiein*).

    Ephemeral is eternal only.

    By virtue of contrary principles things are made themselves conspicuous (Anaxagoras):  light through darkness, darkness through light, etc.

    All is necessity and happening in the same time.

    Man's attitudes in the presence of evil or suffering are:

- primitive passivity:  to bear, tolerate it;
- magic reaction:  to do magic rituals for driving away the bad spirits hidden in objects and beings;
- resignation:  to stay pessimistic, because the evil is irreparable;
- suffering utilization:  to turn suffering to joy, because suffering is necessary and can't be eliminated from our life;
- activist solution:  to accept the suffering and to condemn the evil (Tudor Vianu);

(Müller-Lyer, <Soziologie der Leiden>).

    But what are the man's attitudes in the presence of good or joy?

- ecstasy;
- arrogance;
- indolence;
- decline.

This is the close circuit of man's attitudes in the presence of "-", "0", and "+".

    Idealism and Realism.

Schopenhauer asserts that "world is my representation", which is distorted by the plurality of various imaginations.

Contradictory and alike representations at different individuals.

Running counter Fichte's transcendental idealism (who, in his turn, ran counter Kant's metaphysical determinism), Schopenhauer concludes that beyond the veil of the world there is an absolute reality.

A text on a Chinese Funeral Pillar:
"Detour of non-boundary, statement of non-statement, settlement of those who can't settle were our tortures".

Eliade reveals an "irrecognoscible God", who is present without being made known, an echo of the Buddhist paradox of presence-absence grounded by Nâgârjuna. While Hegel (according to H. Küng) shows a "God who sacrifices himself".

Kant: man must be regarded first as purpose, and then as means. While others said the purpose excuses the means!

That, who said he never lied in his life, is a liar.

X writes on a piece of paper, and put it in an envelope addressed to Y, "today I'm writing you no letter anymore".
Is that an antithesis?

Unfortunately, word puzzles are substitutes of philosophy, especially in the language philosophy.
Should a such thinker be named word theorist or terrorist?
And, yet, I love Frege.

Content is not a form of the form, but it tends to become a form. And reversely.

The spirit couldn't even breathe without opposition and neutralities, would wither itself as a plant...

1.3**. Neutrosophic Transdisciplinarity.**

**A) Definition:**

Neutrosophic Transdisciplinarity means to find common features to uncommon entities: i.e., for vague, imprecise, not-clear-boundary entity <A> one has: <A> ∩ <Non-A> ≠ ∅ or even more <A> ∩ <Anti-A> ≠ ∅.

**B) Multi-Structure and Multi-Space:**

B1) Multi-Concentric-Structure:
Let $S_1$ and $S_2$ be two distinct structures, induced by the ensemble of laws L, which verify the ensembles of axioms $A_1$ and $A_2$ respectively, such that $A_1$ is strictly included in $A_2$. One says that the set M, endowed with the properties:
a) M has an $S_1$-structure;

b) there is a proper subset P (different from the empty set Ø, from the unitary element, from the idempotent element if any with respect to $S_2$, and from the whole set M) of the initial set M, which has an $S_2$-structure;

c) M doesn't have an $S_2$-structure; is called a 2-concentric-structure.

We can generalize it to an n-concentric-structure, for $n \geq 2$ (even infinite-concentric-structure).

(By default, *1-concentric structure* on a set M means only one structure on M and on its proper subsets.)

An **n-concentric-structure** on a set S means a weak structure $\{w(0)\}$ on S

such that there exists a chain of proper subsets

$P(n-1) < P(n-2) < \ldots < P(2) < P(1) < S$,

where '<' means 'included in',

whose corresponding structures verify the inverse chain

$\{w(n-1)\} > \{w(n-2)\} > \ldots > \{w(2)\} > \{w(1)\} > \{w(0)\}$,

where '>' signifies 'strictly stronger' (i.e., structure satisfying more axioms).

For example:

  Say a groupoid D, which contains a proper subset S which is a semigroup, which
  in its turn contains a proper subset M which is a monoid, which contains a proper subset NG which is a non-commutative group, which contains a proper subset CG which is a commutative group, where D includes S, which includes M, which includes NG, which includes CG.

  [This is a 5-concentric-structure.]

B2) **Multi-Space**:

Let $S_1, S_2, \ldots, S_n$ be distinct two by two structures on respectively the sets $M_1, M_2, \ldots, M_k$, where $n \geq 2$ (n may even be infinite).

The structures $S_i$, i = 1, 2, …, n, may not necessarily be distinct two by two; each structure $S_i$ may also be $n_i$-concentric, $n_i \geq 1$.

And the sets $M_i$, i = 1, 2, …, n, may not necessarily be disjoint,

also some sets $M_i$ may be equal to or included in other sets $M_j$, j = 1, 2, …, n.

We define the Multi-Space M as a union of the previous sets:

$M = M_1 \cup M_2 \cup \ldots \cup M_n$, hence we have n (different) structures on M.

A multi-space is a space with many structures that may overlap,

or some structures include others, or the structures may interact and influence each other as in our everyday life.

For example we can construct a geometric multi-space formed by the union of three distinct subspaces: an Euclidean space, a Hyperbolic one, and an Elliptic one.

As particular cases when all $M_i$ sets have the same type of structure, we can define the Multi-Group (or n-group; for example; bigroup, tri-group, etc., when all sets $M_i$ are groups), Multi-Ring (or n-ring, for example biring, tri-ring, etc. when all sets $M_i$ are rings), Multi-Field (n-field), Multi-Lattice (n-lattice), Multi-Algebra (n-algebra), Multi-Module (n-module), and so on - which may be generalized to Infinite-Structure-Space (when all sets have the same type of structure), etc.

{F. Smarandache, "Mixed Non-Euclidean Geometries", 1969.}

Let's introduce new terms:

C) **Psychomathematics**:

A discipline which studies psychological processes in connection with mathematics.

D) **Mathematical Modeling of Psychological Processes**:

a) Improvement of Weber's and Fechner's Laws on sensations and stimuli.

According to the neutrosophic theory, between an <idea> (=spiritual) and an <object> (= material) there are infinitely many states. Then, how can we mix an <idea> with an <object> and obtain something in between: s% spiritual and m% material?  [kind of chemical alloy].
Or, as Boethius, a founder of scholasticism, urged to "join faith to reason" in order to reconcile the Christian judgment with the rational judgment.

For example <mind> and <body> co-exist.  Gustav Theodor Fechner, who inaugurated the experimental psychology, obsessed with this problem, advanced the theory that every object is both
mental and physical (psychophysics).
Fechner's Law, $S = k \cdot \log R$, with S the sensation, R the stimulus, and k a constant,  which is derived from Weber's Law, $\Delta R / R = k$, with $\Delta R$ the increment of stimulus just detectable,
should be improved, because the function Log R is indefinitely increasing
as $R \to \infty$, to

$S(R) = k \cdot \ln R / \ln R_M$, for $R \in [R_m, R_M]$,

and

S(R) = 0, for R ∈ [0, R$_m$) ∪ (R$_M$, ∞),

where k is a positive constant depending on three parameters:
individual being, type of sensation, and the kind of stimulus, and R$_m$, R$_M$ represent the minimum and maximum stimulus magnitude respectively perceptible by the subject, the second one bringing about the death of sensation.
Fechner's "functional relation", as well as later psychologists' power law R = k·S$^n$, with n depending on the kind of stimulus, were upper unbounded, while the beings are surely limited in perception.

S: [0, ∞) ↦ {0} ∪ [S$_m$, S$_M$], with S$_m$, S$_M$ the minimum and maximum perceptible sensation respectively.
Of course R$_m$ > 1, S(R$_m$) = S$_m$, and S(R$_M$) = S$_M$ = k.

Ln, increasing faster, replaces log because the sensation is more rapidly increasing at the beginning, and later going on much slower.
At R = R$_M$, S attains its maximum, beyond whom it becomes flat again, falling to zero.
The beings have a low and high threshold respectively, a range where they may feel a sensation.

Graph of Fechner's Law Improvement:

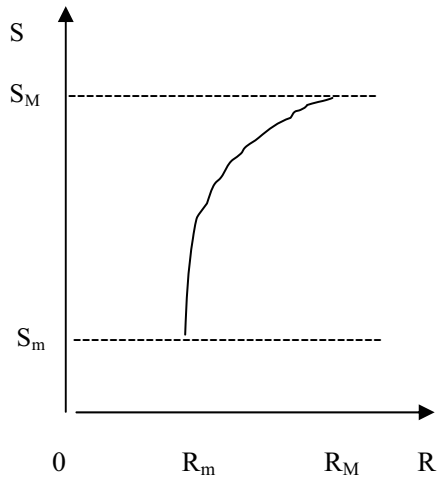

For example in acoustics:  a sound is not heard at the beginning and, if it constantly keeps enlarging its intensity, at a given moment we hear it, and for a while its loudness increases in our ears, until the number of decibels - getting bigger than our possibility of hearing - breaks our eardrums...  We would not hear anything anymore, our sensation died...

Now, if at a given moment t$_0$ the stimulus R remains constant equal to R$_0$ (between the conscious limits of the being, for a long period of time t), and the sensation S(R$_0$) = c, then we get the following formulas:

In the case when the stimulus in not physically or physiologically damaging the individual being: $S_m(t) = c \cdot \log_{1/e}(t+1/e) = -c \cdot \ln(t+1/e)$, for $0 \le t \le \exp(-S_m/c)-1/e$, and 0 otherwise; which is a decreasing function.

In the case when the stimulus is hurting the individual being: $S_m(t) = c \cdot \ln(t+e)$, for $0 \le t \le \exp(S_M/c)-e$, and 0 otherwise; which is an increasing function until the sensation reaches its upper bound; where c, as a constant, depends on individual being, type of sensation, and kind of stimulus.

Examples: i) If a prisoner feels a constant smell in his closed room for days and days, isolated from the exterior, and he doesn't go outside to change the environment, he starts to feel it less and less and after a critical moment he becomes inured to the smell and do not feel it anymore - thus the sensation disappears under the low perceptible limit. ii) If a water drop licks constantly, at the same interval of time, with the same intensity, on the head of a prisoner tied to a pillar, the prisoner after a while will feel the water drop heavier and heavier, will mentally get ill and out of his mind, and will even physically die

- therefore again disappears the sensation, but above the high limit.   See how one can kill someone with a... water drop! iii) If one permanently plays the same song for days and days to a person enclosed in a room without any other noise from outside, that person will be driven crazy, even psychologically die, and the sensation will disappear.

Weber's Law can be improved to $\Delta R/\ln R = k$, with R defined on $[R_m, R_M]$, where k is a constant depending on individual being, type of sensation, and kind of stimulus, due to the fact that the relative threshold $\Delta R$ increases slower with respect to R.

Let's propose a
**b) Synonymity Test**, similar to, and an extension of, the antonym test in psychology, would be a verbal test where the subject must supply as many as possible synonyms of a given word within a as short as possible period of time. How to measure it? The spectrum of supplied synonyms (s), within the measured period of time (t), shows the subject's level of linguistic neutrosophy: s/t.

**c) An Illusion:**
Suppose you travel to a third world country, for example Romania.   You arrive in the capital city of Bucharest, late in the night, and want to exchange a $100 bill to the country's currencies, which are called "lei".   All exchange offices are closed.   A local citizen approaches and proposes you to exchange your bill.   He is a thief.

You give him the $100 bill, he gives you the equivalent in the country's currency, i.e. 25,000 lei.  But the laws of the country do not allow exchange on the street, and both of you know it.

The thief cries "police!", and gives you the dollars back with one hand, while with the other hand takes back his lei, and runs out vanishing behind a building.
The thief has cheated you.
Taken by surprise, you don't realize what had happened, and looking in your hand expecting to see back a $100 bill, actually you see a $1 bill… in your mind, in the very first seconds, it appears the illusion that the $100 bill changed, under your eyes, into a $1 bill!

### E) Psychoneutrosophy:

Psychology of neutral thought, action, behavior, sensation, perception, etc.  This is a hybrid field deriving from psychology, philosophy, economics, theology, etc.
For example, to find the psychological causes and effects of individuals supporting neutral ideologies (neither capitalists, nor communists), politics (not in the left, not in the right), etc.

### F) Socioneutrosophy:
Sociology of neutralities.
For example the sociological phenomena and reasons which determine a country or group of people or class to remain neuter in a military, political, ideological, cultural, artistic, scientific, economical, etc. international or internal war (dispute).

### G) Econoneutrosophy:
Economics of non-profit organizations, groups, such as: churches, philanthropic associations, charities, emigrating foundations, artistic or scientific societies, etc.
How they function, how they survive, who benefits and who loses, why are they necessary, how they improve, how they interact with for-profit companies.

### H) New Types of Philosophies:

*a) Object Philosophy:*  a building through its architecture, a flower, a bird flying, etc. any object are all ideas, or inspire ideas - which are not necessarily to be written down on the paper because they would lose their naturalness and their essence would be distorted. The philosophy should consequently have a universal language, not clung to a specific language (how to translate, for example, Heidegger's dassein, and why to entangle in a notion, syntagme, or word?!).

*b) Concrete Philosophy:*  a drawing, a painting, a canvas, any two-dimensional picture are all ideas and inspire ideas.

*c) Sonorous Philosophy:*  a symphony melody, the jazz music, a sound, any noise are all ideas, or inspire ideas - because they directly work with our unconsciousness.

*d) Fuzzy Philosophy:* there is only a fuzzy border between <A> and <Non-A> and, in consequence, elements which belong (with a certain probability) to both of them, even to <A> and <Anti-A>.
Like the clouds in the sky.
An element *e* belongs 70% to <A> and 30% to <Non-A>.
Or, more organic, *e* belongs 70% to <A>, 20% to <Neut-A> and 10% to <Anti-A> for example.
The di-chotomy between <A> and <Non-A> may be substituted with trichotomy (<A>, <Neut-A>, <Anti-A>) according to our three\ory, and by generalization in a similar way, with plurichotomy onward to transchotomy [ì-chotomy] (continuum-power shades among <A>, <Neut-A>, and <Anti-A>).
And, when the probability is involved, fuzzy-chotomy, or more:
neutro-chotomy.

*e) Applied Philosophy:* philosophical knowledge (such as:
proverbs, aphorisms, maxims, fables, stories) used in our every day's life.

*f) Experimental Philosophy:* philosophical checking and studying of strange, bizarre ideas.

*g) Futurist Philosophy:* ideas created by machines, robots, computers using artificial intelligence;
this is the philosophy of tomorrow.

*h) Nonphilosophy:*
   To make philosophy by not doing any philosophy at all!
Like a mutism.
Everything may mean philosophy: a graffiti (having no words, no letters), any scientific sign or expression displayed on the page...
A poem is a philosophical system. A physics law, a chemical formula, a mathematical equation too.
For example, a blank page also means an idea, a natural phenomenon as well.
Due to the fact that they all make you reflect, meditate, think.
This nonphilosophy becomes, paradoxically, a new kind a philosophy!

   **I) New Types of Philosophical Movements:**

*a) Revisionism:* to review all the philosophical systems, ideas, phenomena, schools, thinkers and rewrite the philosophy as a cumulus of summum bonum.

*b) Inspirationalism:* to look to antecedents for clues and contemporaries for inspiration to get your own research methods and original system.

*c) Recurrentism:* any idea comes from a previous idea and determines another idea, like an infinite recurrent sequence.

*d) Sophisticalism:* the more unintelligible, ambiguous, unsolved, abstract, general... the better!
[This is the style of some people...]

*e) Rejectivism:* a unconscious (and, at some degree, becoming mixed with conscious) will to a priori-ly repel somebody else's system, and totally or partially replace it with yours own.

*f) Paradoxism:* any philosophical idea is true and false in the same time.
Law of the paradoxism:
    Nothing is non-contradictory.
    Nature's essence is antonymic.

## J) Logical and Combinatory Modeling in Experimental Literature:

*a) An Avant-garde Literary Movement, the Paradoxism*
(which uses mathematical paradoxes in artistic creations):
the study of paradoxes as a discipline apart and their use in other fields.

- Basic Thesis of Paradoxism:
    everything has a meaning and a non-meaning
    in a harmony each other.

- Essence of Paradoxism:
    a) sense has a non-sense,
     and reciprocaliy
    b) non-sense has a sense.

- Delimitation from Other Avant-gardes:
    - paradoxism has a significance,
    while dadaism, lettrism, the absurd movement do not;
    - paradoxism especially reveals the contradictions,
    the anti-nomies, the anti-theses, the anti-phrases,
    antagonism, non-conformism, in other words the
    paradoxes of anything (in literature, art, science),
    while futurism, cubism, abstractism and all other
    avant-gardes do not focus on them.

- Directions for Paradoxism:
    - use science methods (especially algorithms) for
    generating (and also studying) contradictory literary and
    artistic works;
    - create contradictory literary and artistic works
    in scientific spaces (using scientific: symbols,
    meta-language, matrices, theorems, lemmas, etc.).

*b) New Types of 'Mathematical' Poetry with Fixed Form*
(using paradoxes and tautologies):
- Paradoxist Distich = a two-line poem such that the second one
contradicts the first, but together they form a unitary meaning
defining (or making connection with) the title.
- Tautological Distich = an apparently redundant two-line poem,
but together the redundant lines give a deeper meaning to the
whole poem defining (or making connection with) the title.
- Dualist Distich
- Paradoxist Tertian
- Tautological Tertian
- Paradoxist Quatrain
- Tautological Quatrain
- Fractal Poem.

*c) New Types of Short Story:*
- Syllogistic Short Story
- Circular Short Story
(F.Smarandache, "Infinite Tale", 1997)

*d) New Types of Drama:*
- Neutrosophic Drama
- Sophistic Drama
- Combinatory Drama = a drama whose scenes are permuted and
combined in so many ways producing over a billion of billions of
different dramas!  (F.Smarandache, "Upside-Down World", 1993)

   Similar definitions for other types of poems, of short
stories, and of dramas.

**Acknowledgements**:

The author would like to thank Drs. C. Le and Ivan Stojmenovic for encouragement
and invitation to write this and the following papers.

# Neutrosophic Logic - A Unifying Field in Logics

*Abstract*:  In this paper one generalizes fuzzy, paraconsistent, and intuitionistic logic to neutrosophic logic.  Many examples are presented, and a survey to the evolution of logics up to neutrosophic logic.

*Keywords and Phrases*: Fuzzy Logic, Paraconsistent Logic, Intuitionistic Logic, Neutrooshpic Logic.

*1991 MSC*:  03B99, 03E99.

## 2.  NEUTROSOPHIC LOGIC:

### 2.1.  Introduction.

As an alternative to the existing logics we propose the Neutrosophic Logic to represent a mathematical model of uncertainty, vagueness, ambiguity, imprecision, undefined, unknown, incompleteness, inconsistency, redundancy, contradiction.  It is a non-classical logic. Eksioglu (1999) explains some of them:

"Imprecision of the human systems is due to the imperfection of knowledge that human receives (observation) from the external world. Imperfection leads to a doubt about the value of a variable, a decision to be taken or a conclusion to be drawn for the actual system. The sources of uncertainty can be stochasticity (the case of intrinsic imperfection where a typical and single value does not exist), incomplete knowledge (ignorance of the totality, limited view on a system because of its complexity) or the acquisition errors (intrinsically imperfect observations, the quantitative errors in measures)."

"Probability (called sometimes the objective probability) process uncertainty of random type (stochastic) introduced by the chance.  Uncertainty of the chance is clarified by the time or by events' occurrence.  The probability is thus connected to the frequency of the events' occurrence."

"The vagueness which constitutes another form of uncertainty is the character of those with contours or limits lacking precision, clearness.  […] For certain objects, the fact to be in or out of a category is difficult to mention. Rather, it is possible to express a partial or gradual membership."

Indeterminacy means degrees of uncertainty, vagueness, imprecision, undefined, unknown, inconsistency, redundancy.

A question would be to try, if possible, to get an axiomatic system for the neutrosophic logic. Intuition is the base for any formalization, because the postulates and axioms derive from intuition.

### 2.2.1.  Definition:

A logic in which each proposition is estimated to have the percentage of truth in a subset T, the percentage of indeterminacy in a subset I, and the percentage of falsity in a subset F, where T, I, F are defined above, is called *Neutrosophic Logic*.

We use a subset of truth (or indeterminacy, or falsity), instead of a number only, because in many cases we are not able to exactly determine the percentages of truth and of falsity but to approximate them: for example a proposition is between 30-40% true and between 60-70% false, even worst: between 30-40% or 45-50% true (according to various analyzers), and 60% or between 66-70% false.

The subsets are not necessary intervals, but any sets (discrete, continuous, open or closed or half-open/half-closed interval, intersections or unions of the previous sets, etc.) in accordance with the given proposition.

A subset may have one element only in special cases of this logic.

Statically T, I, F are subsets, but dynamically they are functions/operators depending on many known or unknown parameters.

Constants: (T, I, F) truth-values, where T, I, F are standard or non-standard subsets of the non-standard interval $\|{}^{-}0, 1^{+}\|$, where $n_{\inf} = \inf T + \inf I + \inf F \geq {}^{-}0$, and $n_{\sup} = \sup T + \sup I + \sup F \leq 3^{+}$. Statically T, I, F are subsets, but dynamically T, I, F are functions/operators depending on many known or unknown parameters.

Atomic formulas: a, b, c, ... .

Arbitrary formulas: A, B, C, ... .

The neutrosophic logic is a formal frame trying to measure the truth, indeterminacy, and falsehood.

My hypothesis is that **no theory is exempted from paradoxes**, because of the language imprecision, metaphoric expression, various levels or meta-levels of understanding/interpretation which might overlap.

## 2.2.2. **Differences between Neutrosophic Logic (NL) and Intuitionistic Fuzzy Logic (IFL)**

The differences between NL and IFL (and the corresponding neutrosophic set and intuitionistic fuzzy set) are:

a) Neutrosophic Logic can distinguish between *absolute truth* (truth in all possible worlds, according to Leibniz) and *relative truth* (truth in at least one world), because NL(absolute truth)=$1^{+}$ while NL(relative truth)=1. This has application in philosophy (see the neutrosophy). That's why the unitary standard interval [0, 1] used in IFL has been extended to the unitary non-standard interval $]^{-}0, 1^{+}[$ in NL. Similar distinctions for absolute or relative falsehood, and absolute or relative indeterminacy are allowed in NL.

b) In NL there is no restriction on T, I, F other than they are subsets of $]^{-}0, 1^{+}[$, thus:

$^{-}0 \leq \inf T + \inf I + \inf F \leq \sup T + \sup I + \sup F \leq 3^{+}$.

This non-restriction allows paraconsistent, dialetheist, and incomplete information to be characterized in NL {i.e. the sum of all three components if they are defined as points, or sum of superior limits of all three components if they are defined as subsets can be >1 (for paraconsistent information coming from different sources) or < 1 for incomplete information}, while that information can not be described in IFL because in IFL the components T (truth), I (indeterminacy), F (falsehood) are restricted either to t+i+f=1 or to $t^2 + f^2 \leq 1$, if T, I, F are all reduced to the points t, i, f respectively, or to sup T + sup I + sup F = 1 if T, I, F are subsets of [0, 1].

   c) In NL the components T, I, F can also be *non-standard* subsets included in the unitary non-standard interval ]-0, 1+[, not only *standard* subsets, included in the unitary standard interval [0, 1] as in IFL.

   d) NL, like dialetheism, can describe  paradoxes, NL(paradox) = (1, 1, 1), while IFL can not describe a paradox because the sum of components should be 1 in IFL.

   e) NL has a better and clear name "neutrosophic" (which means the neutral part: i.e. neither true nor false), while IFL's name "intuitionistic" produces confusion with Intuitionistic Logic, which is something different (see comments by Didier Dubois).

   f) NL permits the utilization of indeterminacy "I" in algebraic structures such that I^2 = I and in graph theory, giving birth to neutrosophic algebraic structures and neutrosophic graphs, which have many applications.

## 2.3. History:

   The Classical Logic, also called Bivalent Logic for taking only two values {0, 1}, or Boolean Logic from British mathematician George Boole (1815-64), was named by the philosopher Quine (1981) "sweet simplicity".

Peirce, before 1910, developed a semantics for three-valued logic in an unpublished note, but Emil Post's dissertation (1920s) is cited for originating the three-valued logic. Here "1" is used for truth, "1/2" for indeterminacy, and "0" for falsehood.  Also, Reichenbach, leader of the logical empiricism, studied it.

The three-valued logic was employed by Halldén (1949), Körner (1960), Tye (1994) to solve Sorites Paradoxes.  They used truth tables, such as Kleene's, but everything depended on the definition of validity.

A three-valued paraconsistent system (LP) has the values: 'true', 'false', and 'both true and false'.  The ancient Indian metaphysics considered four possible values of a statement: 'true (only)', 'false (only)', 'both true and false', and 'neither true nor false'; J. M. Dunn (1976) formalized this in a four-valued paraconsistent system as his First Degree Entailment semantics;

The Buddhist logic added a fifth value to the previous ones, 'none of these' (called *catushkoti*).

   In order to clarify the anomalies in science, Rugina (1949, 1981) proposed an original method, starting first from an economic point of view but generalizing it to any science, to study the equilibrium and disequilibrium of systems.  His Orientation Table comprises seven basic models:

   Model $M_1$ (which is 100% stable)
   Model $M_2$ (which is   95% stable, and   5% unstable)
   Model $M_3$ (which is   65% stable, and 35% unstable)
   Model $M_4$ (which is   50% stable, and  50% unstable)
   Model $M_5$ (which is   35% stable, and  65% unstable)
   Model $M_6$  (which is    5% stable, and  95% unstable)
   Model $M_7$ (which is                    100% unstable)

He gives Orientation Tables for Physical Sciences and Mechanics (Rugina 1989), for the Theory of Probability, for what he called Integrated Logic, and generally for any Natural or Social Science (Rugina 1989).  This is a Seven-Valued Logic.

The {0, $a_1$, ..., $a_n$, 1} Multi-Valued, or Plurivalent, Logic was develop by Łukasiewicz, while Post originated the m-valued calculus.

The many-valued logic was replaced by Goguen (1969) and Zadeh (1975) with an Infinite-Valued Logic (of continuum power, as in the classical mathematical analysis and classical probability) called Fuzzy Logic, where the truth-value can be any number in the closed unit interval [0, 1]. The Fuzzy Set was introduced by Zadeh in 1965.

Rugina (1989) defines an anomaly as "a deviation from a position of stable equilibrium represented by Model $M_1$", and he proposes a Universal Hypothesis of Duality:

"The physical universe in which we are living, including human society and the world of ideas, all are composed in different and changeable proportions of stable (equilibrium) and unstable (disequilibrium) elements, forces, institutions, behavior and value"

and a General Possibility Theorem:

"there is an unlimited number of possible combinations or systems in logic and other sciences".

According to the last assertions one can extend Rugina's Orientation Table in the way that any system in each science is s% stable and u% unstable, with s+u=100 and both parameters $0 \leq s, u \leq 100$, somehow getting to a fuzzy approach.

But, because each system has hidden features and behaviors, and there would always be unexpected occurring conditions we are not able to control - we mean the indeterminacy plays a role as well, a better approach would be the *Neutrosophic Model*:

Any system in each science is s% stable, i% indeterminate, and u% unstable, with s+i+u=100 and all three parameters $0 \leq s, i, u \leq 100$.

Therefore, we finally generalize the fuzzy logic to a transcendental logic, called "neutrosophic logic": where the interval [0, 1] is exceeded, i.e. , the percentages of truth, indeterminacy, and falsity are approximated by non-standard subsets – not by single numbers, and these subsets may overlap and exceed the unit interval in the sense of the non-standard analysis; also the superior sums and inferior sum, $n_{sup} = \sup T + \sup I + \sup F \in {\Vdash}^- 0, 3^+ {\Vert}$, may be as high as 3 or $3^+$, while $n_{inf} = \inf T + \inf I + \inf F \in {\Vdash}^- 0, 3^+ {\Vert}$, may be as low as 0 or $^-0$.

Generally speaking, passing from the attribute "classical" (traditional) to the attribute "modern" (in literature, arts, and philosophy today one says today "postmodern") one invalidates many theorems. Voltaire (1694-1778), a French writer and philosopher, asserted that "the laws in arts are made in order to encroach upon them". Therefore, in neutrosophic logic most of the classical logic laws and its properties are not preserved. Although at a first look neutrosophic logic appears counter-intuitive, maybe abnormal, because the neutrosophic-truth values of a proposition A, NL(A), may even be (1,1,1), i.e. a proposition can completely be true and false and indeterminate at the same time, studying the paradoxes one soon observes that it is intuitive.

The idea of tripartition (truth, falsehood, indeterminacy) appeared in 1764 when J. H. Lambert investigated the credibility of one witness affected by the contrary testimony of another. He generalized Hooper's rule of combination of evidence (1680s), which was a Non-Bayesian approach to find a probabilistic model. Koopman in 1940s introduced the notions of lower and upper probability, followed by Good, and Dempster (1967) gave a rule of combining two arguments. Shafer (1976) extended it to the Dempster-Shafer Theory of Belief Functions by defining the Belief and Plausibility functions and using the rule of inference of Dempster for combining two evidences proceeding from two different sources. Belief function is a connection between fuzzy reasoning and probability. The Dempster-Shafer Theory of Belief Functions is a generalization of the Bayesian Probability (Bayes 1760s, Laplace 1780s); this uses the mathematical probability in a more general way, and is based on probabilistic combination of evidence in artificial intelligence.

In Lambert "there is a chance p that the witness will be faithful and accurate, a chance q that he will be mendacious, and a chance 1-p-q that he will simply be careless" [apud Shafer (1986)]. Therefore three components: accurate, mendacious, careless, which add up to 1.

Van Fraassen introduced the supervaluation semantics in his attempt to solve the Sorites paradoxes, followed by Dummett (1975) and Fine (1975). They all tripartitioned, considering a vague predicate which, having border cases, is undefined for these border cases. Van Fraassen took the vague predicate 'heap' and extended it positively to those objects to which the predicate definitively applies and negatively to those objects to which it definitively doesn't apply. The remaining objects border was called penumbra. A sharp boundary between these two extensions does not exist for a soritical predicate. Inductive reasoning is no longer valid too; if S is a Sorites predicate, the proposition "$\exists n(Sa_n \& \neg Sa_{n+1})$" is false. Thus, the predicate Heap (positive extension) = true, Heap (negative extension) = false, Heap (penumbra) = indeterminate.

Narinyani (1980) used the tripartition to define what he called the "indefinite set", and Atanassov (1982) continued on tripartition and gave five generalizations of the fuzzy set, studied their properties and applications to the neural networks in medicine:

a) Intuitionistic Fuzzy Set (IFS):

Given an universe E, an IFS A over E is a set of ordered triples <universe_element, degree_of_membership_to_A(M), degree_of_non-membership_to_A(N)> such that M+N ≤ 1 and M, N ∈ [0, 1]. When M + N = 1 one obtains the fuzzy set, and if M + N < 1 there is an indeterminacy I = 1-M-N.

b) Intuitionistic L-Fuzzy Set (ILFS):

Is similar to IFS, but M and N belong to a fixed lattice L.

c) Interval-valued Intuitionistic Fuzzy Set (IVIFS):

Is similar to IFS, but M and N are subsets of [0, 1] and sup M + sup N ≤ 1.

d) Intuitionistic Fuzzy Set of Second Type (IFS2):

Is similar to IFS, but $M^2 + N^2 \leq 1$. M and N are inside of the upper right quarter of unit circle.

e) Temporal IFS:

Is similar to IFS, but M and N are functions of the time-moment too.

This neutrosophic logic is the (first) attempt to unify many logics in a single field. However, sometimes a too large generalization may have no practical impact. Such unification theories, or attempts, have been done in the history of sciences:

a) Felix Klein (1872), in his Erlangen program, in geometry, has proposed a common standpoint from which various branches of geometries could be re-organized, interpreted, i.e.:

Given a manifold and a group of transformations of the manifold, to study the manifold configurations with respect to those features that are not altered by the transformations of the group (Klein 1893, p. 67; apud Torretti 1999).

b) Einstein tried in physics to build a Unifying Field Theory that seeks to unite the properties of gravitational, electromagnetic, weak, and strong interactions so that a single set of equations can be used to predict all their characteristics; whether such a theory may be developed it is not known at the present (Illingworth 1991, p. 504).

c) Also, one mentions the Grand Unified Theory, which is a unified quantum field theory of the electromagnetic, weak, and strong interactions (Illingworth 1991, p. 200).

But generalizations become, after some levels, "very general", and therefore not serving at much and, if dealing with indeterminacy, underlying the infinite improbability drive. Would the gain of such total generality offset the losses in specificity? A generalization may be done in one direction, but not in another, while gaining in a bearing but loosing in another.

How to unify, not too much generalizing? Dezert (1999) suggested to develop the less limitative possible theory which remains coherent with certain existing theories. The rules of inferences in this general theory should satisfy many important mathematical properties. "Neutrosophic Logic could permit in the future to solving certain practical problems posed in the domain of research in Data/Information fusion. So far, almost all approaches are based on the Bayesian Theory, Dempster-Shafer Theory, Fuzzy Sets, and Heuristic Methods" (Dezert 1999). Theoretical and technical advances for Information Fusion are probability and statistics, fuzzy sets, possibility, evidential reasoning, random sets, neural networks and neuro-mimetic approaches, and logics (Dezert 2000).

The confidence interval <Bel, Pl> in Dempster-Shafer Theory is the truth subset (T) in the neutrosophic set (or logic). The neutrosophic logic, in addition to the it, contains an indeterminacy set (say indeterminacy interval) and falsehood set (say in-confidence internal).

Łukasiewicz, together with Kotarbiński and Leśniewski from the Warsaw Polish Logic group (1919-1939), questioned the status of truth: eternal, sempiternal (everlasting, perpetual), or both? Also, Łukasiewicz had the idea of *logical probability* in between the two world wars.

Let's borrow from the modal logic the notion of "world", which is a semantic device of what the world might have been like. Then, one says that the neutrosophic truth-value of a statement A, $NL_t(A) = 1^+$ if A is 'true in all possible worlds' (syntagme first used by Leibniz) and all conjunctures, that one may call "absolute truth" (in the modal logic it was named *necessary truth,* Dinulescu-Câmpina (2000) names it 'intangible absolute truth' ), whereas $NL_t(A) = 1$ if A is true in at least one world at some conjuncture, we call this "relative truth" because it is related to a 'specific' world and a specific conjuncture (in the modal logic it was named *possible truth*). Because each 'world' is dynamic, depending on an ensemble of parameters, we introduce the sub-category 'conjuncture' within it to reflect a particular state of the world.

How can we differentiate <the truth behind the truth>? What about the <metaphoric truth>, which frequently occurs in the humanistic field? Let's take the proposition "99% of the politicians are crooked" (Sonnabend 1997, Problem 29, p. 25). "No," somebody furiously comments, "100% of the politicians are crooked, *even more*!" How do we interpret this "even more" (than 100%), i. e. more than the truth?

One attempts to formalize. For n ≥1 one defines the "n-level relative truth" of the statement A if the statement is true in at least n distinct worlds, and similarly "countable-" or "uncountable-level relative truth" as gradual degrees between "first-level relative truth" (1) and "absolute truth" ($1^+$) in the monad $\mu(1^+)$. Analogue definitions one gets by substituting "truth" with "falsehood" or "indeterminacy" in the above.

In *largo sensu* the notion "world" depends on parameters, such as: space, time, continuity, movement, modality, (meta)language levels, interpretation, abstraction, (higher-order) quantification, predication, complement constructions, subjectivity, context, circumstances, etc. Pierre d'Ailly upholds that the truth-value of a proposition depends on the sense, on the metaphysical level, on the language and meta-language; the auto-reflexive propositions (with reflection on themselves) depend on the mode of representation (objective/subjective, formal/informal, real/mental).

In a formal way, let's consider the world W as being generated by the formal system FS. One says that statement A belongs to the world W if A is a well-formed formula (*wff*) in W, i.e. a string of symbols from the alphabet of W that conforms to the grammar of the formal language endowing W. The grammar is conceived as a set of functions (formation rules) whose inputs are symbols strings and outputs "yes" or "no". A formal system comprises a formal language (alphabet and grammar) and a deductive apparatus (axioms and/or rules of inference). In a formal system the rules of inference are syntactically and typographically formal in nature, without reference to the meaning of the strings they manipulate.

Similarly for the neutrosophic falsehood-value, $NL_f(A) = 1^+$ if the statement A is false in all possible worlds, we call it "absolute falsehood", whereas $NL_f(A) = 1$ if the statement A is false in at least one world, we call it "relative falsehood". Also, the neutrosophic indeterminacy-value $NL_i(A) = 1^+$ if the statement A is indeterminate in all possible worlds, we call it "absolute indeterminacy", whereas $NL_i(A) = 1$ if the statement A is indeterminate in at least one world, we call it "relative indeterminacy".

On the other hand, $NL_t(A) = {}^-0$ if A is false in all possible world, whereas $NL_t(A) = 0$ if A is false in at least one world; $NL_f(A) = {}^-0$ if A is true in all possible world, whereas $NL_f(A) = 0$ if A is true in at least one world; and $NL_i(A) = {}^-0$ if A is indeterminate in no possible world, whereas $NL_i(A) = 0$ if A is not indeterminate in at least one world.

The ${}^-0$ and $1^+$ monads leave room for degrees of super-truth (truth whose values are greater than 1), super-falsehood, and super-indeterminacy.

Here there are some corner cases:

There are tautologies, some of the form "B is B", for which $NL(B) = (1^+, {}^-0, {}^-0)$, and contradictions, some of the form "C is not C", for which $NL(C) = ({}^-0, {}^-0, 1^+)$.

While for a paradox, P, $NL(P) = (1,1,1)$. Let's take the Epimenides Paradox, also called the Liar Paradox, "This very statement is false". If it is true then it is false, and if it is false then it is true. But the previous reasoning, due to the contradictory results, indicates a high indeterminacy too. The paradox is the only proposition true and false in the same time in the same world, and indeterminate as well!

Let's take the Grelling's Paradox, also called the heterological paradox [Suber, 1999], "If an adjective truly describes itself, call it 'autological', otherwise call it 'heterological'. Is 'heterological' heterological?" Similarly, if it is, then it is not; and if it is not, then it is.

For a not well-formed formula, nwff, i.e. a string of symbols which do not conform to the syntax of the given logic, NL(nwff) = n/a (undefined). A proposition which may not be considered a proposition was called by the logician Paulus Venetus *flatus voci*. NL(*flatus voci*) = n/a.

## 2.4. An Attempt of Classification of Logics

(upon the following, among many other, criteria):

   a)  The way the connectives, or the operators, or the rules of inferences are defined.
   b)  The definitions of the formal systems of axioms.
   c)  The number of truth-values a proposition can have: two, three, finitely many-values, infinitely many (of continuum power).
   d)  The partition of the interval [0, 1] in propositional values: bi-partition (in degrees of truth and falsehood), or tri-partition (degrees of truth, falsehood, and indeterminacy).
   e)  The distinction between conjunctural (relative) true, conjunctural (relative) false, conjunctural (relative) indeterminacy - designed by 1, with respect to absolute true (or super-truth), absolute false (super-falshood), absolute indeterminacy – designed by $1^+$. Then, if a proposition is absolute true, it is underfalse ($^-0$), i.e. NL(P)=($1^+$, I, $^-0$).

   For example, the neutrosophic truth-value of the proposition "The number of planets of the Sun is divisible by three" is 1 because the proposition is necessary *de re*, i.e. relates to an actual individual mentioned since its truth depends upon the number nine, whereas the neutrosophic truth-value of the proposition "The number of planets of the Sun is the number of its satellites" is $1^+$ because the second proposition is necessary *de dicto*, i.e. relates to the expression of a belief, a possibility since its truth is not dependent upon which number in fact that is. The first proposition might not be true in the future if a new planet is discovered or an existing planet explodes in an asteroid impact, while the second one is always true as being a tautology. This is the difference between the truth-value "1" (dependent truth) and the truth-value "$1^+$" (independent truth).
   f)  The components of the truth values of a proposition summing up to 1 (in boolean logic, fuzzy logic, intuitionistic fuzzy logic), being less than 1 (in intuitionistic logic), or being greater than 1 (in paraconsistent logic, neutrosophic logic). The maximum sum may be 3 in neutrosophic logic, where NL(paradox)=(1,1,1).
   g)  Parameters that influence the truth-values of a proposition. For example in temporal logic the time is involved. A proposition may be true at a time $t_1$, but false at a time $t_2$, or may have some degree of truth in the open interval (0, 1) at a time $t_3$.
   h)  Using approximations of truth-values, or exact values.

   For example, the probabilistic logic, interval-valued fuzzy logic, interval-valued intuitionistic fuzzy logic, possibility logic (Dubois, Prade) use approximations.

   The boolean logic uses exact values, either 0 or 1.
   i)  Studying the paradoxes or not.

   In the neutrosophic logic one can treat the paradoxes, because NL(paradox)=(1,1,1), and in dialetheism. In fuzzy logic FL(paradox)=(1,0) or (0,1)? Because FL(paradox)≠(1,1), due to the fact that the sum of the components should be 1 not greater.
   j)  The external or internal structure of propositions: Sentential (or Propositional) Calculus, which is concerned with logical relations of propositions treated only as a whole, and Predicate (or Functional) Calculus which is concerned besides the logical relations treated as a whole with their internal structure in terms of subject and predicate.

k) Quantification: First-Order (or Lower) Predicate Calculus (quantification is restricted to individuals only, and predicates take only individuals as arguments), Second-Order Predicate Calculus (quantification over individuals and over some classes as well), Higher-Order Calculus (n-predicates take, and quantifiers bind, order n-1 predicates as arguments, for n>1).

l) In proof-theoretic terms:

- Monotonic Logic: let $\Gamma$ be a collection of statements, $\nu_1, \nu_2, \ldots, \nu_n$, and $\varpi, \varphi$ other statements; if $\Gamma \vdash \varphi$ then also $(\Gamma, \varpi) \vdash \varphi$.

- Non-Monotonic Weak Logic: For some $\Gamma, \varpi, \varphi$ one has
  $\Gamma \vdash_{NML} \varphi$ but from $(\Gamma, \varpi)$ does not $\vdash_{NML} \varphi$;

- Non-Monotonic Strong Logic: For some $\Gamma, \varpi, \varphi$, where $\Gamma$ and $\Gamma \wedge \varphi$ are consistent, one has
  $(\Gamma, \varpi) \vdash_{NML} \neg\varphi$.

m) From a traditional standpoint: Classical or Non-classical.

n) Upon inclusion or exclusion of empty domains (and defining the logical validity accordingly), there are Inclusive Predicate Logic, and (Standard) Predicate Logic respectively.

o) Upon the number of arguments the predicates can take, there are Monadic Predicate Logic (predicates take only one argument), Dyadic Predicate Logic (predicates take two arguments), Polyadic Predicate Logic or Logic of Relations (predicates take n>1 arguments).

p) Upon formalization again: Formal Logic, and Informal Logic.

q) Upon types of formalization, there are: Number-Theoretic Predicate Calculus (system with function symbols and individual constants), Pure Predicate Calculus (system without function symbols nor individual constants).

r) Upon standardization: Standard Logic, and Non-Standard Logic.

s) Upon identity: Predicate Logic With Identity (with the axiom $(\forall x)(x=x)$, and the axiom schema $[(x=y) \rightarrow (A \rightarrow A')]^c$, where A' is obtained from A by replacing any free occurrence of x in A with y, and $B^c$ is an arbitrary closure of B), Predicate Logic Without Identity.

t) According to the *ex contradictione quodlibet* (ECQ) principle, from contradictory premises follows anything, there are:

- Explosive logics, which validates it (classical logic, intuitionistic logic);

- Non-Explosive Logics, which invalidate it (paraconsistent logic, neutrosophic logic).

u) According to the Law of Excluded Middle (LEM), either A or ¬A, there are:

- Constructive Logic, which invalidate it (intuitionistic logic, paraconsistent logic, neutrosophic logic);

- Non-Constructive Logics (classical logic).

The criteria are not exhausted. There are sub-classifications too.

Let's take the <u>Modal Logic</u> which is an extension of the Propositional Calculus but with operators that express various modes of truth, such as: necessarily A, possible A, probably A, it is permissible that A, it is believed that A, it has always been true that A. The Modal Logic comprises:

- Alethic Logic (which formalizes the concepts of pertaining to truth and falsehood simultaneously, such as *possibly true* and *necessarily true*); only for this case there are more than two hundred systems of axioms!
- Deontic Logic (which seeks to represent the concepts of *obligatoriness* and *permissibility*);  it is sub-divided into:
    - Standard Deontic Logic,  which has two monadic operators added to the classical propositional calculus: "O" = it ought to be that, and "P" = it is permissible that;
    - Dyadic Deontic Logic, which has two similar dyadic operators added to the classical propositional calculus: "O( / )" = it ought tobe that …, given that …, and P( / ) = is it permissible that …, given that … ;
    - Two-sorted Deontic Logic (Castaňeda 1975) , which distinguishes between *propositions* (which bear truth-values) and *practitions* (which content imperatives, commands, requests).  The deontic operators in this case are: Oi = it is obligatory I that, Pi = it is permissible i that, Wi = it is wrong I that, and Li = it is optional i that. A deontic operator applied to a practition yields a proposition.
- Epistemic Logic (which seeks to represent to concepts of *knowledge*, *belief*, and *ignorance*);
- and Doxastic Logic (which studies the concept of *belief*); it is included in the Epistemic Logic, which is the investigation of epistemic concepts, the main ones being: knowledge, reasonable belief, justification, evidence, certainty.

Dynamic Logic (1970), as a generalization of the modal logic, has a category of expressions interpretable as *propositions* and another category of expressions interpretable as *actions*, with two operators:

  [$\alpha$]A = after every terminating computation according to $\alpha$ it is the case that A;

  <$\alpha$>A = after some terminating computation according to $\alpha$ it is the case that A,

and it is used in the verification of the computer programs.

Combinatory Logic (Schoenfinkel, Haskell Curry, 1920s) is a system for reducing the operational notation of logic, mathematics, or functional language to a sequence of modifications to the input data structure.

Temporal Logic is an extension of Predicate Calculus that includes notation for arguing about when (at what time) statements are true, and employs prefix operators such as:

  ◯x = x is true at the next time;

   x = x is true from now on;

  ◇x = x is eventually true;

or infix operators such as:

  xUy = x is true until y is true;

  xPy = x precedes y;

  xWy = x is weak until y is true.

Temporal Logic studies the Linear Time, which considers only one possible future, and Branching Time, which has two extra operators:

  "A" = all futures,

and   "E" = some futures.

Default Logic (Raymond Reiter 1980) is a formal system with two default operators:

  P:MQ/Q = if P is believed, and Q is consistent with this believe, then Q must be believed;

P:M¬Q/¬Q = if P is believed, and Q is not consistent with this believe, then Q must not be believed.

<u>Tense Logic</u> (Arthur Prior 1967), which is related to the Modal Logic, introduces in the classical logic two operators:

      P = it was the case that … (past tense);

      F = it will be the case that … (future tense).

The truth-value is not static as in classical logic, but changing in time.

<u>Deviant Logics</u> are logics which treat the same classical logic subjects, but in a different way (either by interpreting the connectives and quantifiers non-classically, or rejecting some classical laws): intuitionistic logic, paraconsistent logic, free logic, multi-valued logic.

<u>Free Logic</u> is a system of quantification theory that allows non-denoting singular terms (free variables and individual constants).

In Webster's dictionary (1988) denotation of a term means the class of all particular objects to which the term refers, and connotation of a term means the properties possessed by all the objects in the term's extension.

<u>Erotetic Logic</u> is the logic of questions, answers, and the relations between them. There are (1) imperative approaches (A. Åqvist, J. Hintikka, et al.), epistemic sentences embedded in an imperative sentence system, and (2) interrogative approaches (N. Belnap, T. Kubiński, and others), system of interrogative expressions and their answers.

<u>Relational Logic</u> (Pierce 1870, 1882) is a formal study of the properties of the (binary) relations and the operations on relations.

Because the neutrosphic logic is related to intuitionistic logic, paraconsistent logic, and dialetheism we'll focus more on these types of logics.

<u>Intuitionistic Logic</u> (Brouwer 1907) is a deviant logic from the classic, where the Law of Excluded Middle of Aristotle (A∨¬A) is invalidated. In this logic: a proof of existence, ∃xP(x), does not count unless a method/algorithm of constructing a such x is giving (the interpretation of 'there exists' as 'we can construct' distinguishes between *classical mathematics* and *constructive mathematics* respectively); and a proof of A∨B counts only if a proof of A exhibits or a proof of B. Similarly (Bridges 1997), a proof of A∧B counts if both a proof of A and a proof of B exhibit, a proof of A→B counts if an algorithm is constructed that converts a proof of A into a proof of B, a proof of ¬A means to show that A implies a contradiction, and a proof of ∀xP(x) means to construct an algorithm that applied to any x proves that P(x) holds. As a consequence, the axiom of choice also fails. Brouwer considered some unsolved problem from number theory as proposition A, which is not -- with our present knowledge -- proved true, neither ¬A is proved true. Thus, neutrosophically $NL_t(A∨¬A) < 1$, $NL_f(A∨¬A) < 1$, and $NL_t(A∨B) < 1$, $NL_f(A∨B) < 1$, for some propositions A, B.

<u>Paraconsistent Logic</u> is a logic in which the principle that anything follows from contradictory premises, for all formulas A and B one has A∧¬A⊃B, fails. Therefore, A∧¬A is not always false, i.e. for some A $NL_t(A∧¬A)>0$ or NL(A) = (t, i, f) where t+f>1. It is motivated by dialetheists who support the idea that some contradictions are true, by automated reasoning (information processing) due to inconsistent data stored in computers, and by the fact that people impart opposite beliefs. There are four types of propositional paraconsistent logics (Priest and Tanaka, 1996):

-   Non-Adjunctive Systems (Jaskowski's discussive logic), where the inference

{A, B}⊃A∧B fails; in a discourse a participant's opinion A may be inconsistent with other participant's opinion B on the same subject;

- Non-Truth-Functional Logics (da Costa), which maintains the mechanism of positive logics (classical, intuitionistic) but the value of the negation, ⌐A, is interpreted independently of that of A;
- Many-Valued Systems (Asenjo), many-valued logic which allows both A and ⌐A to be designated (to function as the analogue of truth in a two-valued logic); for example a three-valued paraconsistent system (LP) has the values: 'true', 'false', and 'both true and false', while in a four-valued system (J. M. Dunn 1976) one adds another value 'neither true nor false';
- Relevance Logic (or Relevant Logic) (Wilhelm Ackermann 1956, Alan Anderson and Nuel Belnap 1959-1974) promulgates that the premises of a valid inference must be relevant to the conclusion. The disjunctive syllogism, which states that 'if A∨B and ⌐A are true then so is B', is not admitted in relevance logic, neither in neutrosophic logic. However, Ackermann's rule Gamma, that 'if A∨B and ⌐A are theses then so is B', is admitted.

Dialetheism asserts that some contradictions are true, encroaching upon the Aristotle's Law of Non-Contradiction (LNC) that not both A and ⌐A are true. The dialetheism distinguishes from the trivialism, which views all contradictions as being true. Neither neutrosophic logic is trivialist.

There is a duality (Mortensen 1996) between paraconsistency and intuitionism (i.e. between inconsistency and incompleteness respectively), the Routley * operation (1972) between inconsistent theories and incomplete theories.

Linear Logic (J. Y. Girard 1987) is a resource sensitive logic that emphasizes on state. It employs the central notions of truth from classical logic and of proof construction from intuitionistic logic. Assumptions are considered resources, and conclusions as requirements; A implies B means that the resource A is spent to meet the requirement B. In the deductions there are two structural rules (Scedrov 1999), that allow us to discard or duplicate assumptions (distinguishing linear logic from classical and intuitionistic logics): *contraction*, which stipulates that any assumption once stated may be reused as often as desired, and *weakening* which stipulates that it's possible to carry out a deduction without using all the assumptions. They are replaced by explicit modal logical rules such as: "storage" or "reuse" operator, !A, which means unlimited creation of A, and its dual, ?B, which means unlimited consumption of B.

## 2.5. Comments on Neutrosophic Logic in Comparison to Other Logics.

How to adopt the Gödel-Gentzen *negative translation*, which transforms a formula A of a language L into an equivalent formula A` with no ∨ or ∃, in the neutrosophic predicate logic?

In the Boolean logic a *contingent statement* is a statement that is true under certain conditions and false under others. Then a neutrosophic contingent statement is a statement which has the truth value $(T_1, I_1, F_1)$ under certain conditions and $(T_2, I_2, F_2)$ under others.

The Medieval paradox, called Buridan's Ass after Jean Buridan (near 1295-1356), is a perfect example of complete indeterminacy. An ass, equidistantly from two quantitatively and qualitatively heaps of grain, starves to death because there is no ground for preferring one heap to another.

The neutrosophic value of ass's decision, NL = (0.5, 1, 0.5).

In a two-valued system one regards all the designated values as species of truth and all the anti-designated values as species of falsehood, with truth-value (or falsehood-value) gaps between designated and anti-designated values. In the neutrosophic system one stipulates the non-designated values as species of indeterminacy and, thus, each neutrosophic consequence has degrees of designated, non-designated, and anti-designated values.

Of course, the Law of Excluded Middle (a proposition is either true or false) does not hold in a neutrosophic system.

The Contradiction Law, that no <A> is <Non-A> does not hold too. NL(<A>) may be equivalent with NL(<Non-A>) and often they at least overlap. Neither the law of *Reductio ad absurdum* (or method of indirect proof): $(A \supset \neg A) \supset \neg A$ and $(\neg A \supset A) \supset A$.

Some tautologies (propositions logically necessary, or *true in virtue of form*) in the classical logic might not be tautologies (absolute truth-value propositions) in the neutrosophic logic and, *mutatis mutandis*, some contradictions (propositions logically impossible, or *false in virtue of form*) in the classical logic might not be contradictions (absolute falsehood-value propositions) in the neutrosophic logic.

The mixed hypothetical syllogism Modus Ponens,

    If P then Q
    P
    -----------
    Q

The mixed hypothetical syllogism Modus Tollens,

    If P then Q
    Non Q
    -----------
    Non P

The Inclusive (Weak) Disjunctive Syllogism:

    If (P or Q)
    Non P
    -----------
    Q

The Exclusive (Strong) Disjunctive Syllogism:

    If (either P or Q)
    Non P
    -----------
    Q

Hypothetical Syllogism,

    If P then Q
    If Q then R
    -----------
    If P then R

Constructive Dilemma,

    P or Q
    If P then R
    If Q then R
    -----------
    R

Destructive Dilemma,

    P or Q
    Non P
    ------
    Q

The Polysyllogism, which is formed by many syllogisms such that the conclusion of one becomes a premise of another,

and the Nested Arguments, a chainlike where the conclusion of an argument forms the premise of another where intermediate conclusions are typically left out, are not valid anymore in the neutrosophic logic, but they acquire a more complex form.

Also, the classical *entailment*, which is the effect that a proposition Q is a necessary consequence of another proposition P, P → Q, partially works in the neutrosophic logic. Neither
the informal *fish-hook* symbol, ---], used to show that a proposition Q is an accidental consequence of a proposition P, P ---] Q,  works.

Is it possible in the neutrosophic predicate calculus to transform each formula into an equivalent in prenex form one using the prenex operations?
*Prenex (normal) form* means a formula formalized as follows:

(Qx$_1$)(Qx$_2$)…(Qx$_n$)S,

where "Q" is a universal or existential quantifier, the variables x$_1$, x$_2$, …, x$_n$ are distinct, and S is an open sentence (a well-formed expression containing a free variable). Prenex operation is any operation which transforms any well-formed formula into equivalent in prenex form formula; for example,    $(\exists x)Ax \rightarrow B \equiv (\forall x)(Ax \rightarrow B)$.

In the classical predicate calculus any well-formed formula can be transformed into a prenex form formula.

The double negation, $\neg(\neg A) \equiv A$, which is not valid in intuitionistic logic, is not valid in the neutrosophic logic if one considers the negation operator $\eta_1(A)=1 \ominus NL(A)$ normalized but it is valid if this negation operator is not normalized, also it is valid for the negation operator $\eta_2(A)=(F, I, T)$, where $NL(A)=(T, I, F)$.

Neutrosophic Logic admits non-trivial inconsistent theories.

In stead of saying "a sentence holds (or is assertable)" as in classical logic, one extends to "a sentence p% holds (or is p% assertable)" in neutrosophic logic.  In a more formalized way, a sentence (T, I, F)% holds [or is (T, I, F)% assertable].

A neutrosophic predicate is a vague, incomplete, or not well known attribute, property or function of a subject.  It is a kind of three-valued set function.  If a predicate is applied to more than one subject, it is called neutrosophic relation.

An example:  'Andrew is tall'.

The predicate "tall" is imprecise.  Andrew is maybe tall according to Linda, but short in Jack's opinion, however his tallness is unknown to David.  Everybody judges him in terms of his/her own tallness and acquaintance of him.

The neutrosophic set and logic attempt to better model the non-determinism.  They try:
- to represent the paradoxical results even in science not talking in the humanistic where the paradox is very common;
- to evaluate the peculiarities;
- to illustrate the contradictions and conflicting theories, each true from a specific point of view, false from another one, and perhaps indeterminate from a third perspective;
- to catch the mysterious world of the atom, where the determinism fails; in quantum mechanics we are dealing with systems having an infinite number of degrees of freedom;
- to study submicroscopic particles which behave non-Newtonianly, and some macroscopic phenomena which behave in nearly similar way.

In physics, the light is at once a wave and a particle (photon).  Two contradictory theories were both proven true:

The first one, Wave Theory (Maxwell, Huygens, Fresnel), says that light is a wave due to the interference: two beams of light could cross each other without suffering any damage.

The second one, Particle Theory (Newton, Hertz, Lenard, Planck, Einstein), says that light is corpuscular, due to the photoelectric effect that ultraviolet light is able to evaporate electrons from metal surfaces and to the manner in which light bounces off electrons.

De Broglie reconciled both theories proving that light is a matter wave!  Matter and radiation are at the same time waves and particles.

Let L1(x) be the predicate: "X is of corpuscular nature",
and L2(x) the predicate: "X is of wave nature".
L2(x) is the opposite of L1(x), nonetheless L1(light) = true and L2(light) = true simultaneously.

Also, there exist four different Atom Theories: of Bohr, Heisenberg, Dirac, and Schrödinger respectively, each of them plausibly true in certain conditions (hypotheses).

Another example, from Maxwell's equations an electron does radiate energy when orbits the nucleus, from Bohr's theory an electron does not radiate energy when orbits the nucleus, and both propositions are proved true with our today's knowledge.

Falsehood is infinite, and truthhood quite alike;  in
between, at different degrees, indeterminacy as well.

In the neutrosophic theory:
 between being and nothingness
      existence and nonexistence
      geniality and mediocrity
      certainty and uncertainty
      value and nonvalue
      and generally speaking <A> and <Non-A>
there are infinitely many transcendental states.
And not even 'between', but even beyond them.
An infinitude of infinitudes.
They are degrees of neutralities <Neut-A> combined with <A> and
<Non-A>.

In fact there also are steps:
 between being and being
      existence and existence
      geniality and geniality
      possible and possible
      certainty and certainty
      value and value
      and generally speaking between <A> and <A>.
The notions, in a pure form, last in themselves only (intrinsicalness), but outside they have an interfusion form.
Infinitude of shades and degrees of differentiation:
between white and black there exists an unbounded palette of
colors resulted from thousands of combinations among them.
All is alternative:  progress alternates with setback,
development with stagnation and underdevelopment.

In between objective and subjective there is a plurality of shades.
      In between good and bad...
      In between positive and negative...
 In between possible and impossible
      In between true and false...

In between "A" and "Anti-A"...

As a neutrosophic ellipse:

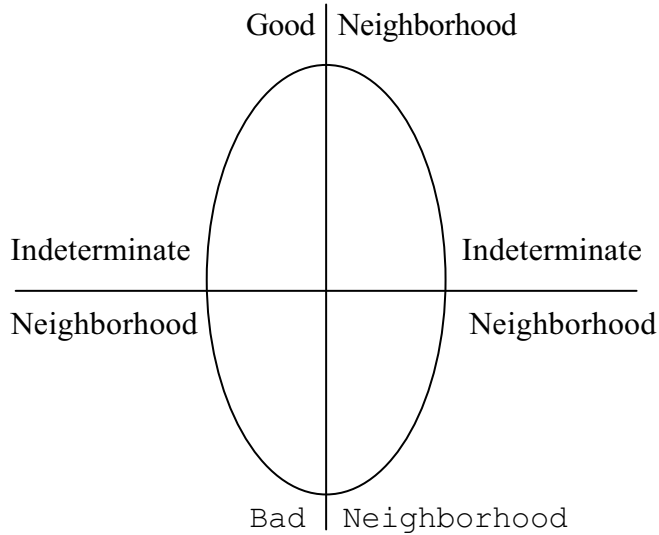

Everything is g% good, i% indeterminate, and b% bad,
where g varies in the subset G, i varies in the subset I, and b varies in the subset B, and G, I, B are included in $]^-0, 1^+[$.

Besides Diderot's dialectics on good and bad ("Rameau's Nephew", 1772), any act has its "god", "indeterminate", and of "bad" as well incorporated.

Rodolph Carnap said:
"Metaphysical propositions are neither true nor false, because they assert nothing, they contain neither knowledge nor error (...)".

Hence, there are infinitely many states between "Good" and "Bad", and generally speaking between "A" and "Anti-A" (and even beyond them), like on the real number line:

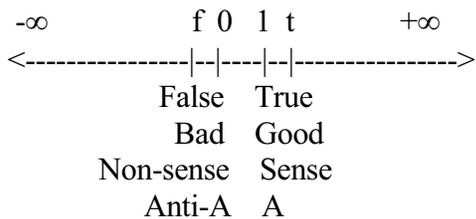

f is the absolute falsity (f<0), t the absolute truth (t>1). In between each oppositing pair, normally in a vicinity of 0.5, are being set up the neutralities.

There exist as many states in between "True" and "False" as in between "Good" and "Bad". Irrational and transcendental standpoints belong to this interval.

Even if an act apparently looks to be only good, or only bad, other hidden sides might be sought.  The ratios

```
        Anti-A      Non-A
        ---------,  ---------
           A           A
```

vary indefinitely.  Are they transfinite?

If a statement is 30%T (true) and 60%I (indeterminate), then it maybe, for example, 15%F (false) .  This is somehow alethic, meaning to simultaneously pertain to truthhood and falsehood, or to truthhood and indeterminacy, or to falsehood and indeterminacy, or even to all three components.

More general, if a statement is 30%T and 60%I, it may be between 5-20%F or 25%F.

## 2.6.  Comparison between Fuzzy Logic and Neutrosophic Logic:

The neutrosophic connectives have a better truth-value definition approach to the real-world systems than the fuzzy connectives.  They are defined on triple non-standard subsets included in the non-standard unit interval $\|^-0, 1^+\|$, while in fuzzy logic they are defined on the closed interval [0, 1].  n_sup is not restricted to 1, but it's enlarged to a monad $\mu(3^+)$,  i.e. a set of hyper-real numbers; similarly, n_inf may be as low as $\mu(^-0)$ , not as 0.

 A paradox, which is true and false in the same time, can not be evaluate in fuzzy logic, because the sum of components should add up to 1, but it is allowed in neutrosophic logic because NL(paradox) maybe (1, 1, 1).

In opposition to Fuzzy Logic, if a proposition <A> is t% true, doesn't necessarily mean it is (100-t)% false.  A better approach is t% true, f% false, and i% indeterminate as in Intuitionistic Fuzzy Logic (Atanassov), where $t \in T$, $i \in I$, $f \in F$, even more general, with n_sup $\leq 3^+$ and n_inf $\geq {}^-0$.

One considers subsets of truth, indeterminacy, and falsity instead of single numbers because of imprecision, uncertainty, and vagueness.

The neutrosophic logical value of <A> is noted by

NL(A) = (T, I, F).  On components one writes:

-   for the truth value $NL_t(A) = T$;
-   for the indeterminacy value $NL_i(A) = I$;
-   for the falsity value $NL_t(A) = F$.

 Neutrosophic Logic means the study of neutrosophic logical values of the propositions.

There exist, for each individual one, PRO parameters, CONTRA parameters, and NEUTER parameters which influence the above values.

Indeterminacy results from any hazard which may occur, from unknown parameters, or from new arising conditions. This resulted from practice.

## 2.7. Neutrosophical Modal Logic:

In modal logic, the primitive operators ◇ 'it is possible that' and □ 'it is necessary that' can be defined by:

$t\_inf(◇A) > 0,$

and, because □A could be regarded as ¬(◇¬A),

$t\_sup(□A) < 1.$

## 2.8. Applications:

Neutrosophic logic is useful in the real-world systems for designing control logic, and may work in quantum mechanics.

# The candidate C, who runs for election in a metropolis M of p people with right to vote, will win.

This proposition is, say, 20-25% true (percentage of people voting for him), 35-45% false (percentage of people voting against him), and 40% or 50% indeterminate (percentage of people not coming to the ballot box, or giving a blank vote - not selecting anyone, or giving a negative vote - cutting all candidates on the list).

# Tomorrow it will rain.

This proposition is, say, 50% true according to meteorologists who have investigated the past years' weather, between 20-30% false according to today's very sunny and droughty summer, and 40% undecided.

# This is a heap.

As an application to the Sorites paradoxes, we may now say that this proposition is 80% true, 40% false, and 25-35% indeterminate (the neutrality comes for we don't know exactly where is the difference between a heap and a non-heap; and, if we approximate the border, our 'accuracy' is subjective). Vagueness plays here an important role.

We are not able to distinguish the difference between yellow and red as well if a continuum spectrum of colors is painted on a wall imperceptibly changing from one into another.

We would be able to say at a given moment that a section is both yellow and red in the same time, or neither one!

A paradox within a *Sorites paradox*: a frontal bald man, with a hair high density on the remaining region of his head, may have more hairs on his head that another man who is not bald but the skin surface of his head and the hair density are smaller than the previous one.

## 2.9. Definition of Neutrosophic Logical Connectives:

The connectives (rules of inference, or operators), in any non-bivalent logic, can be defined in various ways, giving rise to lots of distinct logics. For example, in three-valued logic, where three possible values are possible: true, false, or undecided, there are 3072 such logics!

(Weisstein, 1998)  A single change in one of any connective's truth table is enough to form a (completely) different logic.

The rules are hypothetical or factual.  How to choose them?  The philosopher Van Fraassen (1980) [see Shafer, 1986] commented that such rules may always be controvertible "for it always involves the choice of one out of many possible but nonactual worlds".  There are general rules of combination, and ad hoc rules.

For an applied logic to artificial intelligence, a better approach, the best way would be to define the connectives recursively (Dubois, Prade), changing/adjusting the definitions after each step in order to improve the next result.  This might be comparable to approximating the limit of a convergent sequence, calculating more and more terms, or by calculating the limit of a function successively substituting the argument with values closer and closer to the critical point.  The recurrence allows evolution and self-improvement.

Or to use *greedy algorithms*, which are combinatorial algorithms that attempt at each iteration as much improvement as possible unlike myopic algorithms that look at each iteration only at very local information as with steepest descent method.

As in non-monotonic logic, we make assumptions, but we often err and must jump back, revise our assumptions, and start again.  We may add rules which don't preserve monotonicity.

In bio-mathematics Heitkoetter and Beasley (1993-1999) present the *evolutionary algorithms* which are used "to describe computer-based problem solving systems which employ computational models of some of the known mechanisms of evolution as key elements in their design and implementation". They simulate, via processes of selection, mutation, and reproduction, the evolution of individual structures.  The major evolutionary algorithms studied are: genetic algorithm (a model of machine learning based on genetic operators), evolutionary programming (a stochastic optimization strategy based on linkage between parents and their offspring; conceived by L. J. Fogel in 1960s), evolution strategy, classifier system, genetic programming.

Pei Wang devised a Non-Axiomatic Reasoning System as an intelligent reasoning system, where intelligence means working and adopting with insufficient knowledge and resources.

The inference mechanism (endowed with rules of transformation or rules of production) in neutrosophy should be non-monotonic and should comprise ensembles of recursive rules, with preferential rules and secondary ones (priority order), in order to design a good expert system.  One may add new rules and eliminate old ones proved unsatisfactory.  There should be strict rules, and rules with exceptions.  Recursivity is seen as a computer program that learns from itself.  The statistical regression method may be employed as well to determine a best algorithm of inference.

Non-monotonic reasoning means to make assumptions about things we don't know.  Heuristic methods may be involved in order to find successive approximations.

In terms of the previous results, a default neutrosophic logic may be used instead of the normal inference rules.  The distribution of possible neutrosophic results serves as an orientating frame for the new results.  The flexible, continuously refined, rules obtain iterative and gradual approaches of the result.

A comparison approach is employed to check the result (conclusion) p by studying the opposite of this: what would happen if a non-p conclusion occurred?  The inconsistence of information shows up in the result, if not eliminated from the beginning.  The data bases should be stratified. There exist methods to construct preferable coherent sub-bases within incoherent bases. In Multi-Criteria Decision one exploits the complementarity of different criteria and the complementarity of various sources.

For example, the Possibility Theory (Zadeh 1978, Dubois, Prade) gives a better approach than the Fuzzy Set Theory (Yager) due to self-improving connectives. The Possibility Theory is proximal to the Fuzzy Set Theory, the difference between these two theories is the way the fusion operators are defined.

One uses the definitions of neutrosophic probability and neutrosophic set operations. Similarly, there are many ways to construct such connectives according to each particular problem to solve; here we present the easiest ones.

## 2.9.1. N-norm and N-conorm in Neutrosophic Logic and Neutrosophic Set

As a generalization of T-norm and T-conorm from the Fuzzy Logic and Set, we now introduce the N-norm and N-conorm for the Neutrosophic Logic and Set.

We define a *partial relation order* on the neutrosophic set/logic in the following way:
$x(T_1, I_1, F_1) \leq y(T_2, I_2, F_2)$ iff (if and only if) $T_1 \leq T_2, I_1 \geq I_2, F_1 \geq F_2$ for crisp components.

And, in general, for subunitary set components:
$x(T_1, I_1, F_1) \leq y(T_2, I_2, F_2)$ iff

$$\inf T_1 \leq \inf T_2, \sup T_1 \leq \sup T_2,$$
$$\inf I_1 \geq \inf I_2, \sup I_1 \geq \sup I_2,$$
$$\inf F_1 \geq \inf F_2, \sup F_1 \geq \sup F_2.$$

If we have mixed - crisp and subunitary - components, or only crisp components, we can transform any crisp component, say "a" with $a \in [0,1]$ or $a \in ]^-0, 1^+[$, into a subunitary set $[a, a]$. So, the definitions for subunitary set components should work in any case.

**N-norms**

$N_n$: $( ]^-0,1^+[ \times ]^-0,1^+[ \times ]^-0,1^+[ )^2 \rightarrow ]^-0,1^+[ \times ]^-0,1^+[ \times ]^-0,1^+[$

$N_n (x(T_1,I_1,F_1), y(T_2,I_2,F_2)) = (N_nT(x,y), N_nI(x,y), N_nF(x,y)),$
where $N_nT(.,.), N_nI(.,.), N_nF(.,.)$ are the truth/membership, indeterminacy, and respectively falsehood/nonmembership components.

$N_n$ have to satisfy, for any x, y, z in the neutrosophic logic/set M of the universe of discourse U, the following axioms:
a) Boundary Conditions: $N_n(x, 0) = 0$, $N_n(x, 1) = x$.
b) Commutativity: $N_n(x, y) = N_n(y, x)$.
c) Monotonicity: If $x \leq y$, then $N_n(x, z) \leq N_n(y, z)$.
d) Associativity: $N_n(N_n (x, y), z) = N_n(x, N_n(y, z))$.

There are cases when not all these axioms are satisfied, for example the associativity when dealing with the neutrosophic normalization after each neutrosophic operation. But, since we work with approximations, we can call these **N-pseudo-norms**, which still give good results in practice.

$N_n$ represent the *and* operator in neutrosophic logic, and respectively the *intersection* operator in neutrosophic set theory.

Let $J \in \{T, I, F\}$ be a component.
Most known N-norms, as in fuzzy logic and set the T-norms, are:
• The Algebraic Product N-norm: $N_{n-algebraic}J(x, y) = x \cdot y$
• The Bounded N-Norm: $N_{n-bounded}J(x, y) = \max\{0, x + y - 1\}$
• The Default (min) N-norm: $N_{n-min}J(x, y) = \min\{x, y\}$.

A general example of N-norm would be this.
Let $x(T_1, I_1, F_1)$ and $y(T_2, I_2, F_2)$ be in the neutrosophic set/logic M. Then:
$$\mathbf{N_n(x, y) = (T_1 \wedge T_2, I_1 \vee I_2, F_1 \vee F_2)}$$
where the "$\wedge$" operator, acting on two (standard or non-standard) subunitary sets, is a N-norm (verifying the above N-norms axioms); while the "$\vee$" operator, also acting on two (standard or non-standard) subunitary sets, is a N-conorm (verifying the below N-conorms axioms).

For example, $\wedge$ can be the Algebraic Product T-norm/N-norm, so $T_1 \wedge T_2 = T_1 \cdot T_2$ (herein we have a product of two subunitary sets – using simplified notation); and $\vee$ can be the Algebraic Product T-conorm/N-conorm, so $T_1 \vee T_2 = T_1 + T_2 - T_1 \cdot T_2$ (herein we have a sum, then a product, and afterwards a subtraction of two subunitary sets).

Or $\wedge$ can be any T-norm/N-norm, and $\vee$ any T-conorm/N-conorm from the above and below; for example the easiest way would be to consider the *min* for crisp components (or *inf* for subset components) and respectively *max* for crisp components (or *sup* for subset components).

If we have crisp numbers, we can at the end neutrosophically normalize.

## N-conorms

$N_c$: ( $]^-0,1^+[ \times ]^-0,1^+[ \times ]^-0,1^+[$ )$^2 \rightarrow ]^-0,1^+[ \times ]^-0,1^+[ \times ]^-0,1^+[$

$N_c$ ($x(T_1, I_1, F_1)$, $y(T_2, I_2, F_2)$) = ($N_c T(x,y)$, $N_c I(x,y)$, $N_c F(x,y)$),
where $N_n T(.,.)$, $N_n I(.,.)$, $N_n F(.,.)$ are the truth/membership, indeterminacy, and respectively falsehood/nonmembership components.

$N_c$ have to satisfy, for any x, y, z in the neutrosophic logic/set M of universe of discourse U, the following axioms:
a) Boundary Conditions: $N_c(x, 1) = 1$, $N_c(x, 0) = x$.
b) Commutativity: $N_c (x, y) = N_c(y, x)$.
c) Monotonicity: if $x \leq y$, then $N_c(x, z) \leq N_c(y, z)$.
d) Associativity: $N_c (N_c(x, y), z) = N_c(x, N_c(y, z))$.

There are cases when not all these axioms are satisfied, for example the associativity when dealing with the neutrosophic normalization after each neutrosophic operation. But, since we work with approximations, we can call these **N-pseudo-conorms**, which still give good results in practice.

$N_c$ represent the *or* operator in neutrosophic logic, and respectively the *union* operator in neutrosophic set theory.

Let $J \in \{T, I, F\}$ be a component.
Most known N-conorms, as in fuzzy logic and set the T-conorms, are:

- The Algebraic Product N-conorm: $N_{c-algebraic}J(x, y) = x + y - x \cdot y$
- The Bounded N-conorm: $N_{c-bounded}J(x, y) = \min\{1, x + y\}$
- The Default (max) N-conorm: $N_{c-max}J(x, y) = \max\{x, y\}$.

A general example of N-conorm would be this.

Let $x(T_1, I_1, F_1)$ and $y(T_2, I_2, F_2)$ be in the neutrosophic set/logic M. Then:

$$\mathbf{N_n(x, y) = (T_1 \lor T_2, I_1 \land I_2, F_1 \land F_2)}$$

Where – as above - the "$\land$" operator, acting on two (standard or non-standard) subunitary sets, is a N-norm (verifying the above N-norms axioms); while the "$\lor$" operator, also acting on two (standard or non-standard) subunitary sets, is a N-conorm (verifying the above N-conorms axioms).

For example, $\land$ can be the Algebraic Product T-norm/N-norm, so $T_1 \land T_2 = T_1 \cdot T_2$ (herein we have a product of two subunitary sets); and $\lor$ can be the Algebraic Product T-conorm/N-conorm, so $T_1 \lor T_2 = T_1 + T_2 - T_1 \cdot T_2$ (herein we have a sum, then a product, and afterwards a subtraction of two subunitary sets).

Or $\land$ can be any T-norm/N-norm, and $\lor$ any T-conorm/N-conorm from the above; for example the easiest way would be to consider the *min* for crisp components (or *inf* for subset components) and respectively *max* for crisp components (or *sup* for subset components).

If we have crisp numbers, we can at the end neutrosophically normalize.

Since the min/max (or inf/sup) operators work the best for subunitary set components, let's present their definitions below. They are extensions from subunitary intervals {defined in [3]} to any subunitary sets. Analogously we can do for all neutrosophic operators defined in [3].

Let $x(T_1, I_1, F_1)$ and $y(T_2, I_2, F_2)$ be in the neutrosophic set/logic M.

**Neutrosophic Conjunction/Intersection**:

  $x \land y = (T_\land, I_\land, F_\land)$,

   where $\inf T_\land = \min\{\inf T_1, \inf T_2\}$

    $\sup T_\land = \min\{\sup T_1, \sup T_2\}$

    $\inf I_\land = \max\{\inf I_1, \inf I_2\}$

    $\sup I_\land = \max\{\sup I_1, \sup I_2\}$

    $\inf F_\land = \max\{\inf F_1, \inf F_2\}$

    $\sup F_\land = \max\{\sup F_1, \sup F_2\}$

**Neutrosophic Disjunction/Union**:

  $x \lor y = (T_\lor, I_\lor, F_\lor)$,

   where $\inf T_\lor = \max\{\inf T_1, \inf T_2\}$

    $\sup T_\lor = \max\{\sup T_1, \sup T_2\}$

    $\inf I_\lor = \min\{\inf I_1, \inf I_2\}$

    $\sup I_\lor = \min\{\sup I_1, \sup I_2\}$

    $\inf F_\lor = \min\{\inf F_1, \inf F_2\}$

    $\sup F_\lor = \min\{\sup F_1, \sup F_2\}$

**Neutrosophic Negation/Complement:**

  $\mathscr{C}(x) = (T_\mathscr{C}, I_\mathscr{C}, F_\mathscr{C})$,

   where $T_\mathscr{C} = F_1$

    $\inf I_\mathscr{C} = 1 - \sup I_1$

    $\sup I_\mathscr{C} = 1 - \inf I_1$

    $F_\mathscr{C} = T_1$

Upon the above Neutrosophic Conjunction/Intersection, we can define the
**Neutrosophic Containment:**
> We say that the neutrosophic set A is included in the neutrosophic set B of the universe of discourse U,
> iff for any $x(T_A, I_A, F_A) \in A$ with $x(T_B, I_B, F_B) \in B$ we have:
> $\inf T_A \leq \inf T_B$ ; $\sup T_A \leq \sup T_B$;
> $\inf I_A \geq \inf I_B$ ; $\sup I_A \geq \sup I_B$;
> $\inf F_A \geq \inf F_B$ ; $\sup F_A \geq \sup F_B$.

**Remarks**.

a) The non-standard unit interval $]^-0, 1^+[$ is merely used for philosophical applications, especially when we want to make a distinction between relative truth (truth in at least one world) and absolute truth (truth in all possible worlds), and similarly for distinction between relative or absolute falsehood, and between relative or absolute indeterminacy.

But, for technical applications of neutrosophic logic and set, the domain of definition and range of the N-norm and N-conorm can be restrained to the normal standard real unit interval [0, 1], which is easier to use, therefore:

$$N_n: ( [0,1] \times [0,1] \times [0,1] )^2 \rightarrow [0,1] \times [0,1] \times [0,1]$$

and

$$N_c: ( [0,1] \times [0,1] \times [0,1] )^2 \rightarrow [0,1] \times [0,1] \times [0,1].$$

b) Since in NL and NS the sum of the components (in the case when T, I, F are crisp numbers, not sets) is not necessary equal to 1 (so the normalization is not required), we can keep the final result un-normalized.
But, if the normalization is needed for special applications, we can normalize at the end by dividing each component by the sum all components.
If we work with intuitionistic logic/set (when the information is incomplete, i.e. the sum of the crisp components is less than 1, i.e. *sub-normalized*), or with paraconsistent logic/set (when the information overlaps and it is contradictory, i.e. the sum of crisp components is greater than 1, i.e. *over-normalized*), we need to define the neutrosophic measure of a proposition/set.
If $x(T,I,F)$ is a NL/NS, and T,I,F are crisp numbers in [0,1], then the **neutrosophic vector norm** of variable/set x is the sum of its components:
$$N_{vector-norm}(x) = T+I+F.$$
Now, if we apply the $N_n$ and $N_c$ to two propositions/sets which maybe intuitionistic or paraconsistent or normalized (i.e. the sum of components less than 1, bigger than 1, or equal to 1), x and y, what should be the neutrosophic measure of the results $N_n(x,y)$ and $N_c(x,y)$ ?
Herein again we have more possibilities:
- either the product of neutrosophic measures of x and y:
$N_{vector-norm}(N_n(x,y)) = N_{vector-norm}(x) \cdot N_{vector-norm}(y)$,
- or their average:
$N_{vector-norm}(N_n(x,y)) = (N_{vector-norm}(x) + N_{vector-norm}(y))/2$,
- or other function of the initial neutrosophic measures:

$N_{vector-norm}(N_n(x,y)) = f(N_{vector-norm}(x), N_{vector-norm}(y))$, where $f(.,.)$ is a function to be determined according to each application.

Similarly for $N_{vector-norm}(N_c(x,y))$.

Depending on the adopted neutrosophic vector norm, after applying each neutrosophic operator the result is neutrosophically normalized. We'd like to mention that "**neutrosophically normalizing**" doesn't mean that the sum of the resulting crisp components should be 1 as in fuzzy logic/set or intuitionistic fuzzy logic/set, but the sum of the components should be as above: either equal to the product of neutrosophic vector norms of the initial propositions/sets, or equal to the neutrosophic average of the initial propositions/sets vector norms, etc.

In conclusion, we neutrosophically normalize the resulting crisp components T`,I`,F` by multiplying each neutrosophic component T`,I`,F` with S/( T`+I`+F`), where

S= N$_{\text{vector-norm}}$(N$_n$(x,y)) for a N-norm or S= N$_{\text{vector-norm}}$(N$_c$(x,y)) for a N-conorm - as defined above.

c) If T, I, F are subsets of [0, 1] the problem of neutrosophic normalization is more difficult.
   i)   If sup(T)+sup(I)+sup(F) < 1, we have an *intuitionistic proposition/set*.
   ii)  If inf(T)+inf(I)+inf(F) > 1, we have a *paraconsistent proposition/set*.
   iii) If there exist the crisp numbers t ∈ T, i ∈ I, and f ∈ F such that t+i+f =1, then we can say that we have a *plausible normalized proposition/set*.
        But in many such cases, besides the normalized particular case showed herein, we also have crisp numbers, say t$_1$ ∈ T, i$_1$ ∈ I, and f$_1$ ∈ F such that t$_1$+i$_1$+f$_1$ < 1 (incomplete information) and t$_2$ ∈ T, i$_2$ ∈ I, and f$_2$ ∈ F such that t$_2$+i$_2$+f$_2$ > 1 (paraconsistent information).

## 2.9.1.1. Examples of Neutrosophic Operators which are N-norms or N-pseudonorms or, respectively N-conorms or N-pseudoconorms.

We define a binary **neutrosophic conjunction (intersection)** operator, which is a particular case of a N-norm (neutrosophic norm, a generalization of the fuzzy T-norm):

$$c_N^{TIF} : \left([0,1] \times [0,1] \times [0,1]\right)^2 \to [0,1] \times [0,1] \times [0,1]$$

$$c_N^{TIF}(x,y) = \left(T_1 T_2, I_1 I_2 + I_1 T_2 + T_1 I_2, F_1 F_2 + F_1 I_2 + F_1 T_2 + F_2 T_1 + F_2 I_1\right).$$

The neutrosophic conjunction (intersection) operator $x \wedge_N y$ component truth, indeterminacy, and falsehood values result from the multiplication

$$\left(T_1 + I_1 + F_1\right) \cdot \left(T_2 + I_2 + F_2\right)$$

since we consider in a prudent way $T \prec I \prec F$, where "$\prec$" is a **neutrosophic relationship** and means "weaker", i.e. the products $T_i I_j$ will go to $I$, $T_i F_j$ will go to $F$, and $I_i F_j$ will go to $F$ for all i, j ∈ {1,2}, i≠j, while of course the product $T_1 T_2$ will go to T, $I_1 I_2$ will go to I, and $F_1 F_2$ will go to F (or reciprocally we can say that $F$ prevails in front of $I$ which prevails in front of $T$, and this neutrosophic relationship is transitive):

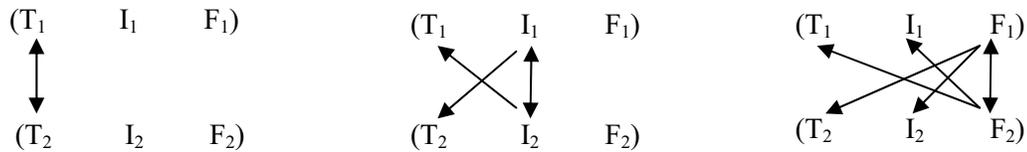

So, the truth value is $T_1 T_2$, the indeterminacy value is $I_1 I_2 + I_1 T_2 + T_1 I_2$ and the false value is $F_1 F_2 + F_1 I_2 + F_1 T_2 + F_2 T_1 + F_2 I_1$. The norm of $x \wedge_N y$ is $\left(T_1 + I_1 + F_1\right) \cdot \left(T_2 + I_2 + F_2\right)$. Thus, if $x$ and $y$ are normalized, then $x \wedge_N y$ is also normalized. Of course, the reader can redefine the neutrosophic conjunction operator, depending on application, in a different way, for example in a more optimistic way, i.e. $I \prec T \prec F$ or $T$ prevails with respect to $I$, then we get:

$$c_N^{ITF}(x,y) = \left(T_1 T_2 + T_1 I_2 + T_2 I_1, I_1 I_2, F_1 F_2 + F_1 I_2 + F_1 T_2 + F_2 T_1 + F_2 I_1\right).$$

Or, the reader can consider the order $T \prec F \prec I$, etc.

Let's also define the unary **neutrosophic negation** operator:

$$n_N : [0,1] \times [0,1] \times [0,1] \rightarrow [0,1] \times [0,1] \times [0,1]$$

$$n_N(T, I, F) = (F, I, T)$$

by interchanging the truth $T$ and falsehood $F$ vector components.

Similarly, we now define a binary **neutrosophic disjunction (or union)** operator, where we consider the neutrosophic relationship $F \prec I \prec T$:

$$d_N^{FIT} : \left([0,1] \times [0,1] \times [0,1]\right)^2 \rightarrow [0,1] \times [0,1] \times [0,1]$$

$$d_N^{FIT}(x,y) = \left(T_1 T_2 + T_1 I_2 + T_2 I_1 + T_1 F_2 + T_2 F_1, I_1 I_2 + I_2 F_1 + I_1 I_2, F_1 F_2\right)$$

We consider as neutrosophic norm of the neutrosophic variable $x$, where $NL(x) = T_1 + I_1 + F_1$, the sum of its components: $T_1 + I_1 + F_1$, which in many cases is 1, but can also be positive <1 or >1.

Or, the reader can consider the order $F \prec T \prec I$, in a pessimistic way, i.e. focusing on indeterminacy I which prevails in front of the truth T, or other **neutrosophic order** of the neutrosophic components T,I,F depending on the application.
Therefore,

$$d_N^{FTI}(x,y) = \left(T_1 T_2 + T_1 F_2 + T_2 F_1, I_1 F_2 + I_2 F_1 + I_1 I_2 + T_1 I_2 + T_2 I_1, F_1 F_2\right)$$

### 2.9.1.2. Neutrophic Composition k-Law

Now, we define a more general neutrosophic composition law, named k-law, in order to be able to define neutrosophic k-conjunction/intersection and neutrosophic k-disjunction/union for k variables, where $k \geq 2$ is an integer.

Let's consider $k \geq 2$ neutrosophic variables, $x_i(T_i, I_i, F_i)$, for all $i \in \{1, 2, ..., k\}$. Let's denote

$$T = (T_1, ..., T_k)$$

$$I = (I_1, ..., I_k)$$

$$F = (F_1, ..., F_k).$$

We now define a neutrosophic composition law $o_N$ in the following way:

$$o_N : \{T, I, F\} \rightarrow [0,1]$$

If $z \in \{T, I, F\}$ then $z_{o_N} z = \prod_{i=1}^{k} z_i$.

If $z, w \in \{T, I, F\}$ then

$$z_{o_N} w = w_{o_N} z \quad = \sum_{\substack{r=1 \\ \{i_1,\dots,i_r,j_{r+1},\dots,j_k\} \equiv \{1,2,\dots,k\} \\ (i_1,\dots,i_r) \in C^r(1,2,\dots,k) \\ (j_{r+1},\dots,j_k) \in C^{k-r}(1,2,\dots,k)}}^{k-1} z_{i_1} \dots z_{i_r} w_{j_{r+1}} \dots w_{j_k}$$

where $C^r(1,2,\dots,k)$ means the set of combinations of the elements $\{1,2,\dots,k\}$ taken by $r$. [Similarly for $C^{k-r}(1,2,\dots,k)$.]

In other words, $z_{o_N} w$ is the sum of all possible products of the components of vectors $z$ and $w$, such that each product has at least a $z_i$ factor and at least a $w_j$ factor, and each product has exactly $k$ factors where each factor is a different vector component of $z$ or of $w$. Similarly if we multiply three vectors:

$$T_{o_N} I_{o_N} F = \sum_{\substack{u,v,k-u-v=1 \\ \{i_1,\dots,i_u,j_{u+1},\dots,j_{u+v},l_{u+v+1},\dots,l_k\} \equiv \{1,2,\dots,k\} \\ (i_1,\dots,i_u) \in C^u(1,2,\dots,k),(j_{u+1},\dots,j_{u+v}) \in \\ \in C^v(1,2,\dots,k),(l_{u+v+1},\dots,l_k) \in C^{k-u-v}(1,2,\dots,k)}}^{k-2} T_{i_1\dots i_u} I_{j_{u+1}\dots j_{u+v}} F_{l_{u+v+1}}\dots F_{l_k}$$

Let's see an example for $k = 3$.

$$x_1(T_1, I_1, F_1)$$
$$x_2(T_2, I_2, F_2)$$
$$x_3(T_3, I_3, F_3)$$

$$T_{o_N} T = T_1 T_2 T_3, \quad I_{o_N} I = I_1 I_2 I_3, \quad F_{o_N} F = F_1 F_2 F_3$$

$$T_{o_N} I = T_1 I_2 I_3 + I_1 T_2 I_3 + I_1 I_2 T_3 + T_1 T_2 I_3 + T_1 I_2 T_3 + I_1 T_2 T_3$$

$$T_{o_N} F = T_1 F_2 F_3 + F_1 T_2 F_3 + F_1 F_2 T_3 + T_1 T_2 F_3 + T_1 F_2 T_3 + F_1 T_2 T_3$$

$$I_{o_N} F = I_1 F_2 F_3 + F_1 I_2 F_3 + F_1 F_2 I_3 + I_1 I_2 F_3 + I_1 F_2 I_3 + F_1 I_2 I_3$$

$$T_{o_N} I_{o_N} F = T_1 I_2 F_3 + T_1 F_2 I_3 + I_1 T_2 F_3 + I_1 F_2 T_3 + F_1 I_2 T_3 + F_1 T_2 I_3$$

For the case when indeterminacy $I$ is not decomposed in subcomponents {as for example $I = P \cup U$ where $P =$ paradox (true and false simultaneously) and $U =$ uncertainty (true or false, not sure which one)}, the previous formulas can be easily written using only three components as:

$$T_{o_N} I_{o_N} F = \sum_{i,j,r \in \mathcal{P}(1,2,3)} T_i I_j F_r$$

where $\mathcal{P}(1,2,3)$ means the set of permutations of $(1,2,3)$ i.e.

$$\{(1,2,3),(1,3,2),(2,1,3),(2,3,1,),(3,1,2),(3,2,1)\}$$

$$z_{o_N} w = \sum_{\substack{i=1 \\ (i,j,r) \equiv (1,2,3) \\ (j,r) \in \mathcal{P}^2(1,2,3)}}^{3} z_i w_j w_{j_r} + w_i z_j z_r$$

This neurotrophic law is associative and commutative.

### 2.9.1.3. Neutrophic Logic and Set k-Operators

Let's consider the neutrophic logic crispy values of variables $x, y, z$ (so, for k = 3):

$$NL(x) = (T_1, I_1, F_1) \text{ with } 0 \leq T_1, I_1, F_1 \leq 1$$

$$NL(y) = (T_2, I_2, F_2) \text{ with } 0 \leq T_2, I_2, F_2 \leq 1$$

$$NL(z) = (T_3, I_3, F_3) \text{ with } 0 \leq T_3, I_3, F_3 \leq 1$$

In neutrosophic logic it is not necessary to have the sum of components equals to 1, as in intuitionist fuzzy logic, i.e. $T_k + I_k + F_k$ is not necessary 1, for $1 \leq k \leq 3$

As a particular case, we define the tri-nary conjunction neutrosophic operator:

$$c_{3N}^{TIF} : \left([0,1] \times [0,1] \times [0,1]\right)^3 \rightarrow [0,1] \times [0,1] \times [0,1]$$

$$c_{3N}^{TIF}(x,y,z) = \left(T_{o_N} T, I_{o_N} I + I_{o_N} T, F_{o_N} F + F_{o_N} I + F_{o_N} T\right)$$

If all x, y, z are normalized, then $c_{3N}^{TIF}(x, y, z)$ is also normalized.

If x, y, or $y$ are non-normalized, then $\left|c_{3N}^{TIF}(x,y,z)\right| = |x| \cdot |y| \cdot |z|$, where |w| means norm of w.

$c_{3N}^{TIF}$ is a 3-N-norm (neutrosophic norm, i.e. generalization of the fuzzy T-norm).

Again, as a particular case, we define the unary negation neutrosophic operator:

$$n_N : [0,1] \times [0,1] \times [0,1] \rightarrow [0,1] \times [0,1] \times [0,1]$$

$$n_N(x) = n_N (T_1, I_1, F_1) = (F_1, I_1, T_1).$$

Let's consider the vectors:

$$T = \begin{pmatrix} T_1 \\ T_2 \\ T_3 \end{pmatrix}, \quad I = \begin{pmatrix} I_1 \\ I_2 \\ I_3 \end{pmatrix} \text{ and } F = \begin{pmatrix} F_1 \\ F_2 \\ F_3 \end{pmatrix}.$$

We note $T_{\bar{x}} = \begin{pmatrix} F_1 \\ T_2 \\ T_3 \end{pmatrix}$, $T_{\bar{y}} = \begin{pmatrix} T_1 \\ F_2 \\ T_3 \end{pmatrix}$, $T_{\bar{z}} = \begin{pmatrix} T_1 \\ T_2 \\ F_3 \end{pmatrix}$, $T_{\overline{xy}} = \begin{pmatrix} F_1 \\ F_2 \\ T_3 \end{pmatrix}$, etc.

and similarly

$$F_{\bar{x}} = \begin{pmatrix} T_1 \\ F_2 \\ F_3 \end{pmatrix}, \quad F_{\bar{y}} = \begin{pmatrix} F_1 \\ T_2 \\ F_3 \end{pmatrix}, \quad F_{\overline{xz}} = \begin{pmatrix} T_1 \\ F_2 \\ T_3 \end{pmatrix}, \text{ etc.}$$

For shorter and easier notations let's denote $z_{o_N} w = zw$ and respectively $z_{o_N} w_{o_N} v = zwv$ for the vector neutrosophic law defined previously.

Then the neutrosophic tri-nary conjunction/intersection of neutrosophic variables x, y, and z is:

$$c_{3N}^{TIF}(x,y,z) = \left(TT, II + IT, FF + FI + FT + FIT\right) =$$
$$= \big(T_1T_2T_3, I_1I_2I_3 + I_1I_2T_3 + I_1T_2I_3 + T_1I_2I_3 + I_1T_2T_3 + T_1I_2T_3 + T_1T_2I_3,$$
$$F_1F_2F_3 + F_1F_2I_3 + F_1I_2F_3 + I_1F_2F_3 + F_1I_2I_3 + I_1F_2I_3 + I_1I_2F_3 +$$
$$+ F_1F_2T_3 + F_1T_2F_3 + T_1F_2F_3 + F_1T_2T_3 + T_1F_2T_3 + T_1T_2F_3 +$$
$$+ T_1I_2F_3 + T_1F_2I_3 + I_1F_2T_3 + I_1T_2F_3 + F_1I_2T_3 + F_1T_2\,I_3\big).$$

Similarly, the neutrosophic tri-nary disjunction/union of neutrosophic variables x, y, and z is:

$$d_{3N}^{FIT}(x,y,z) = \left(TT + TI + TF + TIF, II + IF, FF\right) =$$

$(T_1T_2T_3 + T_1T_2I_3 + T_1I_2T_3 + T_1I_2I_3 + I_1T_2T_3 + I_1I_2T_3 + I_1T_2I_3 + T_1T_2F_3 + T_1F_2T_3 + F_1T_2T_3 + T_1F_2F_3 + F_1T_2F_3 + F_1F_2T_3 + T_1I_2F_3 + T_1F_2I_3 + I_1F_2T_3 + I_1T_2F_3 + F_1I_2T_3 + F_1T_2I_3, I_1I_2I_3 + I_1I_2F_3 + I_1F_2I_3 + F_1I_2I_3 + I_1F_2F_3 + F_1I_2F_3 + F_1F_2I_3, F_1F_2F_3)$

Surely, other neutrosophic orders can be used for tri-nary conjunctions/intersections and respectively for tri-nary disjunctions/unions among the componenets T, I, F.

2.9.2. Other Neutrosophic Operators.

One notes the neutrosophic logical values of the propositions $A_1$ and $A_2$ by
$NL(A_1) = ( T_1, I_1, F_1 )$ and $NL(A_2) = ( T_2, I_2, F_2 )$. Then we neutrosophically normalize.

**2.9.2.1. Negation:**

$NL(\neg A_1) = ( \{1\} \ominus T_1, \{1\} \ominus I_1, \{1\} \ominus F_1 )$.

**2.9.2.2. Conjunction:**

$NL(A_1 \wedge A_2) = ( T_1 \odot T_2, I_1 \odot I_2, F_1 \odot F_2 )$.
(And, in a similar way, generalized for n propositions.)

**2.9.2.3. Weak or inclusive disjunction:**

$NL(A_1 \vee A_2) = ( T_1 \oplus T_2 \ominus T_1 \odot T_2, I_1 \oplus I_2 \ominus I_1 \odot I_2, F_1 \oplus F_2 \ominus F_1 \odot F_2 )$.
(And, in a similar way, generalized for n propositions.)

**2.9.2.4. Strong or exclusive disjunction:**

$NL(A_1 \underline{\vee} A_2) =$
$( T_1 \odot (\{1\} \ominus T_2) \oplus T_2 \odot (\{1\} \ominus T_1) \ominus T_1 \odot T_2 \odot (\{1\} \ominus T_1) \odot (\{1\} \ominus T_2),$
$I_1 \odot (\{1\} \ominus I_2) \oplus I_2 \odot (\{1\} \ominus I_1) \ominus I_1 \odot I_2 \odot (\{1\} \ominus I_1) \odot (\{1\} \ominus I_2),$
$F_1 \odot (\{1\} \ominus F_2) \oplus F_2 \odot (\{1\} \ominus F_1) \ominus F_1 \odot F_2 \odot (\{1\} \ominus F_1) \odot (\{1\} \ominus F_2) )$.
(And, in a similar way, generalized for n propositions.)

**2.9.2.5. Material conditional (implication):**

$NL(A_1 \mapsto A_2) = ( \{1\} \ominus T_1 \oplus T_1 \odot T_2, \{1\} \ominus I_1 \oplus I_1 \odot I_2, \{1\} \ominus F_1 \oplus F_1 \odot F_2 )$.

**2.9.2.6. Material biconditional (equivalence):**

$NL(A_1 \leftrightarrow A_2) = ( (\{1\} \ominus T_1 \oplus T_1 \odot T_2) \odot (\{1\} \ominus T_2 \oplus T_1 \odot T_2),$
$(\{1\} \ominus I_1 \oplus I_1 \odot I_2) \odot (\{1\} \ominus I_2 \oplus I_1 \odot I_2),$
$(\{1\} \ominus F_1 \oplus F_1 \odot F_2) \odot (\{1\} \ominus F_2 \oplus F_1 \odot F_2) )$.

**2.9.2.7. Sheffer's connector:**

NL(A$_1$ | A$_2$) = NL(¬A$_1$ ∨ ¬A$_2$) = ( {1}⊖T$_1$⊙T$_2$, {1}⊖I$_1$⊙I$_2$, {1}⊖F$_1$⊙F$_2$ ).

**2.9.2.8. Peirce's connector:**

NL(A$_1$↓A$_2$) = NL(¬A$_1$ ∧ ¬A$_2$) =
    = ( ({1}⊖T$_1$) ⊙ ({1}⊖T$_2$), ({1}⊖I$_1$) ⊙ ({1}⊖I$_2$), ({1}⊖F$_1$) ⊙ ({1}⊖F$_2$) ).

In the above neutrosophic operators at section 2.9.2 we have used the Algebraic Product N-Norm and N-conorm. Similarly we can do for Bounded, min/max, and other N-norm/N-conorms, then we neutrosphically normalize.

**2.10. Comments on Neutrosophic Operators:**

The conjunction is well defined, associative, commutative, admits a unit element $U$ = (1, 1, 1), but no element whose truth-component is different from 1 is inversable.
The conjunction is not absorbent, i.e. t( A∧(A∧B) ) ∧ t(A),
except for the cases when t(A) = 0, or t(A) = t(B) = 1.
    The weak disjunction is well-defined, associative, commutative, admits a unit element $O$ = *(0, 0, 0)*,
but no element whose truth-component is different from 0 is inversable.
The disjunction is not absorbent, i.e. t( A∨(A∨B) ) ∨ t(A),
except for the cases when one of t(A) = 1, or t(A) = t(B) = 0.
    None of them is distributive with respect to the other.
    De Morgan laws do not apply either.
    Therefore (NL, ∧, ∨, ¬), where NL is the set of neutrosophic logical propositions, is not an algebra.
Nor ($\mathcal{P}$($∥$ ⁻0, 1⁺$∥$), ∩, ∪, C), where $\mathcal{P}$($∥$ ⁻0, 1⁺$∥$) is the set of all subsets of $∥$ ⁻0, 1⁺$∥$, and $C$(A) is the neutrosophic complement of A.

    One names a set N, endowed by two associative unitary internal laws, * and #, which are not inversable except for their unit elements respectively, and not distributive with respect to each other, *Ninversity*.
If both laws are commutative, then N is called a *Commutative Ninversity*.
    For a better understanding of the neutrosophic logic one needs to study the commutative ninversity.

**2.11. Other Types of Neutrosophic Logical Connectors:**

    There are situations when we have to more focus on the percentage of falsity or of indeterminacy than on the percentage of truth.

 Thus, we define in a similar way the logical connectors, but the main component will then be considered the last or second one respectively.  An intriguing idea would be to take the arithmetic average of the corresponding components of the truth-, indeterminacy-, and false-neutrosophic connectors.
Let's reconsider the previous notations.

**2.12. Neutrosophic Topologies.**

A) General Definition of NT:

Let M be a non-empty set.

Let $x(T_A, I_A, F_A) \in A$ with $x(T_B, I_B, F_B) \in B$ be in the neutrosophic set/logic M, where A and B are subsets of M.  Then (see Section 2.9.1 about N-norms / N-conorms and examples):

$\qquad A \cup B = \{x \in M, x(T_A \vee T_B, I_A \wedge I_B, F_A \wedge F_B)\}$,

$\qquad A \cap B = \{x \in M, x(T_A \wedge T_B, I_A \vee I_B, F_A \vee F_B)\}$,

$\qquad \mathscr{C}(A) = \{x \in M, x(F_A, I_A, T_A)\}$.

A General Neutrosophic Topology on the non-empty set M is a family $\eta$ of Neutrosophic Sets in M satisfying the following axioms:

- $\mathbf{0}(0,0,1)$ and $\mathbf{1}(1,0,0) \in \eta$ ;

- If $A, B \in \eta$, then $A \cap B \in \eta$ ;

- If the family $\{A_k, k \in K\} \subset \eta$, then $\bigcup_{k \in K} A_k \in \eta$ .

B) An alternative version of NT

 -We cal also construct a Neutrosophic Topology on NT = $]^{-}0, 1^{+}[$, considering the associated family of standard or non-standard subsets included in NT, and the empty set $\emptyset$, called open sets, which is closed under set union and finite intersection.

Let A, B be two such subsets. The union is defined as:

$A \cup B = A+B-A\cdot B$, and the intersection as: $A \cap B = A\cdot B$. The complement of A, $C(A) = \{1^{+}\}$-A, which is a closed set. {When a non-standard number occurs at an extremity of an internal, one can write "]" instead of "(" and "[" instead of ")".} The interval NT, endowed with this topology, forms a *neutrosophic topological space*.

In this example we have used the Algebraic Product N-norm/N-conorm. But other Neutrosophic Topologies can be defined by using various N-norm/N-conorm operators.

In the above defined topologies, if all x's are paraconsistent or respectively intuitionistic, then one has a

Neutrosophic Paraconsistent Topology, respectively Neutrosophic Intuitionistic Topology.

**2.13. Neutrosophic Sigma-Algebra**.

The collection of all standard or non-standard subsets of ]⁻0, 1⁺[, constitute a *neutrosophic sigma-algebra* (or *neutrosophic $\sigma$-algebra*), because the set itself, the empty set $\emptyset$, the complements in the set of all members, and all countable unions of members belong to the -power set P(]⁻0, 1⁺[).The complement of a subset is defined above.

The interval NT, endowed with this sigma-algebra, forms a *neutrosophic measurable space*.

### 2.14. Generalizations:

When the sets are reduced to an element only respectively, then

t_sup = t_inf = t, i_sup = i_inf = i, f_sup = f_inf = f,

and n_sup = n_inf = n = t+i+f.

Hence, the neutrosophic logic generalizes:

-the intuitionistic logic, which supports incomplete theories (for 0 < n < 1, 0 ≤ t, i, f ≤ 1);

-the fuzzy logic (for n = 1 and i = 0, and 0 ≤ t, f ≤ 1); from "CRC Concise Encyclopedia of Mathematics", by Eric W. Weisstein, 1998, the fuzzy logic is "an extension of two-valued logic such that statements need not to be True or False, but may have a degree of truth between 0 and 1";

-the Boolean logic (for n = 1 and i = 0, with t, f either 0 or 1);

-the multi-valued logic (for 0 ≤ t, i, f ≤ 1);

definition of <many-valued logic> from "The Cambridge Dictionary of Philosophy", general editor Robert Audi, 1995, p. 461: "propositions may take many values beyond simple truth and falsity, values functionally determined by the values of their components"; Lukasiewicz considered three values (1, 1/2, 0). Post considered m values, etc. But they varied in between 0 and 1 only. In the neutrosophic logic a proposition may take values even greater than 1 (in percentage greater than 100%) or less than 0.

- the *paraconsistent logic* (for n > 1, with all t,i,f < 1);
  - the *dialetheism*, which says that some contradictions are true (for some propositions t = f = 1 and i ≥ 0; some paradoxes can be denoted this way too);
  - the *trivialism*, which says that all contradictions are true (for all propositions t = f = 1 and i ≥ 0; some paradoxes can be denoted this way too);
  - the *fallibilism*, which says that uncertainty belongs to every proposition (for i > 0);
- the *paradoxist logic (*or *paradoxism)*, based on paradoxes (t=f=1);
- the *pseudo-paradoxist logic (*or *pseudo-paradoxism)*, based on pseudo-paradoxes (0 < i ≤ 1+, t=1 and 0<f < 1 or 0<t<1 and f=1);
- the *tautological logic (*or *tautologism)*, based on tautologies (i,f < 0, t > 1).

Compared with all other logics, the neutrosophic logic and intuitionistic fuzzy logic introduce a percentage of "indeterminacy" - due to unexpected parameters hidden in some propositions, or unknowness, or God's will, but only neutrosophic logic let each component t, i, f be even boiling *over 1* (overflooded) or freezing *under 0* (underdried): to be able to make distinction between relative truth and absolute truth, and between relative falsity and absolute falsity.

For example: in some tautologies t > 1, called "overtrue". Similarly, a proposition may be "overindeterminate" (for i > 1, in some paradoxes), "overfalse" (for f > 1, in some unconditionally false propositions); or "undertrue" (for t < 0, in some unconditionally false propositions), "underindeterminate" (for i < 0, in some unconditionally true or false propositions), "underfalse" (for f < 0, in some unconditionally true propositions).

This is because we should make a distinction between unconditionally true (t > 1, f < 0 and i < 0) and conditionally true propositions (t ≤ 1, and 0<f ≤ 1 or 0<i ≤ 1).

While in classical true/false logic it is possible to define precisely 2m different m-ary operators for each m>0 (Charles D. Ashbacher),

the *neutrosophic m-ary operators* may be defined in unaccountably infinite different ways. The good operatorial selection would lead to applications in neural networks, automated reasoning, quantum physics, and probabilistic models.

*Dempster-Shafer Theory* doesn't work for some classes of examples:
1) Assume the universe U = {A, B, C}. If $m_1(A) = a$, $m_1(B) = 0$, $m_1(C)=1-a$, where $0 < a < 1$ and a is very close to 1, and $m_2(A) = 0$, $m_2(B) = b$, $m_2(C)=1-b$, where $0 < b <1$ and b is very close to 1, then $m_1+m_2(C) = 1$! This example generalizes Zadeh's (1984).

   Dezert (2000) defends the theory because, as he asserts, in this case the mass fusion is impossible for the sources of evidences are entirely conflictive.
2) Even more, the previous example can be extended to k > 2 masses $m_1, …, m_k$ weighting k+1 exclusive events of the universe U = {$A_1, …, A_k, A_{k+1}$}, such that for all i≠j, 1 ≤ i, j ≤ k, $m_i(A_i)=a_i$ and $m_i(A_j)=0$, where $0 < a_i < 1$ and $a_i$ is very close to 1, and $m_i(A_{k+1})=1-a_i$.
3) In the following particular example, $m_1(A)=.11$, $m_1(B)=.11$, $m_1(C)=.11$, $m_1(D)=.67$, and $m_2(A)=.11$, $m_2(B)=.11$, $m_2(C)=.11$, $m_2(D)=.67$, using the Dempster's rule of combining evidences one gets $m1+m2(D)=.925185$, which is a 38.0873% increment jump from the two equal evidences of .67, and it looks counter-intuitive. Why not a smaller jump?

In the paraconsistent logic one cannot derive all statements from a contradiction, *ex contradictione quodlibet* fails. In the neutrosophic logic from a given contradiction one can derive a specific statement only, depending on the neutrosophic operator used and the given particular contradiction.

In the *dialetheism* it works the metaphysical thesis that some contradictions are true. In the neutrosophic logic there are contradictions denoted by t=f=1, which means 100% true and 100% false in the same time; even more, it is possible to have propositions which are, say, 70% true and 60% false (considering different sources or criteria) - the truth- and falsehood-components overlap (especially in pseudo paradoxes), while in the fuzzy logic it is not - because the components should sum to 100%, i.e. 70% true and 30% false.

What is the logic of the logic? We study the apparently illogic of the logic, as well as the logic of the illogic.

There exist two main types of truth: the true truth and the false truth, besides the intermediate shades of truth. And similarly for the falsity: the true falsity and the false falsity, beside the intermediate shades of falsity.

The neutrosophic logic unifies many logics; it is like Felix Klein's program in geometry, or Einstein's unified field in physics.

In *Propositional Calculus* a statement may be decidable or undecidable. In the First-Order Logic, due to the quantifiers, a statement may be semi-decidable. In Neutrosophic Logic a statement may be p%-decidable, q%-undecidable, $^-0 \leq p, q \leq 100^+$.

# Neutrosophic Set - A Unifying Field in Sets

*Abstract*: In this paper one generalizes fuzzy, paraconsistent, and intuitionistic sets to neutrosophic set. Many examples are presented.

*Keywords and Phrases*: Fuzzy Set, Paraconsistent Set, Intuitionistic Set, Neutrosophic Set.

*1991 MSC*: 03B99, 03E99.

## 3. NEUTROSOPHIC SET:

### 3.1. Definition:

Let T, I, F be real standard or non-standard subsets of $\|{}^-0, 1^+\|$,

with      sup T = t_sup, inf T = t_inf,

           sup I = i_sup, inf I = i_inf,

           sup F = f_sup, inf F = f_inf,

and       n_sup = t_sup+i_sup+f_sup,

           n_inf = t_inf+i_inf+f_inf.

Let U be a universe of discourse, and M a set included in U. An element x from U is noted with respect to the set M as x(T, I, F) and belongs to M in the following way:

it is t% true in the set, i% indeterminate (unknown if it is) in the set, and f% false, where t varies in T, i varies in I, f varies in F.

Statically T, I, F are subsets, but dynamically T, I, F are functions/operators depending on many known or unknown parameters.

### 3.2. General Examples:

Let A and B be two neutrosophic sets.

One can say, by language abuse, that any element neutrosophically belongs to any set, due to the percentages of truth/indeterminacy/falsity involved, which varies between 0 and 1 or even less than 0 or greater than 1.

Thus: x(50,20,30) belongs to A (which means, with a probability of 50% x is in A, with a probability of 30% x is not in A, and the rest is undecidable); or y(0,0,100) belongs to A (which normally means y is not for sure in A); or z(0,100,0) belongs to A (which means one does know absolutely nothing about z's affiliation with A).

More general, x( (20-30), (40-45)∪[50-51], {20,24,28} ) belongs to the set A, which means:

- with a probability in between 20-30% x is in A (one cannot find an exact approximate because of various sources used);

- with a probability of 20% or 24% or 28% x is not in A;

- the indeterminacy related to the appurtenance of x to A is in between 40-45% or between 50-51% (limits included);

The subsets representing the appurtenance, indeterminacy, and falsity may overlap, and n_sup = 30+51+28 > 100 in this case.

### 3.3. Physics Examples:

a)  For example the Schrodinger's Cat Theory says that the quantum state of a photon can basically be in more than one place in the same time, which translated to the neutrosophic set means that an element (quantum state) belongs and does not belong to a set (one place) in the same time; or an element (quantum state) belongs to two different sets (two different places) in the same time.  It is a question of "alternative worlds" theory very well represented by the neutrosophic set theory.

In Schroedinger's Equation on the behavior of electromagnetic waves and "matter waves" in quantum theory, the wave function Psi which describes the superposition of possible states may be simulated by a neutrosophic function, i.e. a function whose values are not unique for each argument from the domain of definition (the vertical line test fails, intersecting the graph in more points).

Don't we better describe, using the attribute "neutrosophic" than "fuzzy" or any others, a quantum particle that neither exists nor non-exists?

b) How to describe a particle $\zeta$ in the infinite micro-universe that belongs to two distinct places $P_1$ and $P_2$ in the same time?  $\zeta \in P_1$ and $\zeta \notin P_1$ as a true contradiction, or $\zeta \in P_1$ and $\zeta \in \neg P_1$.

### 3.4. Philosophical Examples:

Or, how to calculate the truth-value of Zen (in Japanese) / Chan (in Chinese) doctrine philosophical proposition: the present is eternal and comprises in itself the past and the future?

In Eastern Philosophy the contradictory utterances form the core of the Taoism and Zen/Chan (which emerged from Buddhism and Taoism) doctrines.

How to judge the truth-value of a metaphor, or of an ambiguous statement, or of a social phenomenon which is positive from a standpoint and negative from another standpoint?

There are many ways to construct them, in terms of the practical problem we need to simulate or approach.  Below there are mentioned the easiest ones:

### 3.5. Application:

A cloud is a neutrosophic set, because its borders are ambiguous, and each element (water drop) belongs with a neutrosophic probability to the set (e.g. there are a kind of separated water drops, around a compact mass of water drops, that we don't know how to consider them: in or out of the cloud).

Also, we are not sure where the cloud ends nor where it begins, neither if some elements are or are not in the set.  That's why the percent of indeterminacy is required and the neutrosophic probability (using subsets - not numbers - as components) should be used for better modeling: it is a more organic, smooth, and especially accurate estimation.  Indeterminacy is the zone of ignorance of a proposition's value, between truth and falshood.

### 3.6. Neutrosophic Set Operations:

See the above Section 2.9.1 for the definition and examples of N-conorms.

One notes, with respect to the sets A and B over the universe U,

$x = x(T_1, I_1, F_1) \in A$ and $x = x(T_2, I_2, F_2) \in B$, by mentioning x's *neutrosophic probability appurtenance*.

And, similarly, $y = y(T', I', F') \in B$. In the below examples we consider the Algebraic Product N-norm/N-conorm and then we neutrosopihcally normalize. Similarly for other N-norms/N-conorms.

### 3.6.1. Complement of A:

If $x( T_1, I_1, F_1 ) \in A$,
then $x( \{1\} \ominus T_1, \{1\} \ominus I_1, \{1\} \ominus F_1 ) \in C(A)$.

### 3.6.2. Intersection:

If $x( T_1, I_1, F_1 ) \in A$, $x( T_2, I_2, F_2 ) \in B$,
then $x( T_1 \odot T_2, I_1 \odot I_2, F_1 \odot F_2 ) \in A \cap B$.

### 3.6.3. Union:

If $x( T_1, I_1, F_1 ) \in A$, $x( T_2, I_2, F_2 ) \in B$,
then $x( T_1 \oplus T_2 \ominus T_1 \odot T_2, I_1 \oplus I_2 \ominus I_1 \odot I_2, F_1 \oplus F_2 \ominus F_1 \odot F_2 ) \in A \cup B$.

### 3.6.4. Difference:

If $x( T_1, I_1, F_1 ) \in A$, $x( T_2, I_2, F_2 ) \in B$,
then $x( T_1 \ominus T_1 \odot T_2, I_1 \ominus I_1 \odot I_2, F_1 \ominus F_1 \odot F_2 ) \in A \setminus B$,
because $A \setminus B = A \cap C(B)$.

### 3.6.5. Cartesian Product:

If $x( T_1, I_1, F_1 ) \in A$, $y( T', I', F' ) \in B$,
then $( x( T_1, I_1, F_1 ), y( T', I', F' ) ) \in A \times B$.

### 3.6.6. M is a subset of N:

If $x( T_1, I_1, F_1 ) \in M \Rightarrow x( T_2, I_2, F_2 ) \in N$,
where $\inf T_1 \leq \inf T_2$, $\sup T_1 \leq \sup T_2$, but $\inf I1 >= \inf I2$, $\sup I1 >= \sup I2$, and $\inf F_1 \geq \inf F_2$, $\sup F_1 \geq \sup F_2$.

### 3.6.7. Neutrosophic n-ary Relation:

Let $A_1, A_2, \ldots, A_n$ be arbitrary non-empty sets.
A Neutrosophic n-ary Relation $R$ on $A_1 \times A_2 \times \ldots \times A_n$ is defined as a subset of the Cartesian product $A_1 \times A_2 \times \ldots \times A_n$, such that for each ordered n-tuple $(x_1, x_2, \ldots, x_n)(T, I, F)$, T represents the degree of validity, I the degree of indeterminacy, and F the degree of non-validity respectively of the relation $R$.

It is related to the definitions for the *Intuitionistic Fuzzy Relation* independently given by Atanassov (1984, 1989), Toader Buhaescu (1989), Darinka Stoyanova (1993), Humberto Bustince Sola and P. Burillo Lopez (1992-1995).

### 3.7. Generalizations and Comments:

From the intuitionistic logic, paraconsistent logic, dialetheism, fallibilism, paradoxes, pseudoparadoxes, and tautologies we transfer the "adjectives" to the sets, i.e. to intuitionistic set (set incompletely known), paraconsistent set, dialetheist set, faillibilist set (each element has a percentage of indeterminacy), paradoxist set (an element may belong and may not belong in the same time to the set), pseudoparadoxist set, and tautological set respectively.

Hence, the neutrosophic set generalizes:

- the *intuitionistic set*, which supports incomplete set theories (for $0 < n < 1$, $0 \le t$, $i$, $f \le 1$) and incomplete known elements belonging to a set;

- the *fuzzy set* (for $n = 1$ and $i = 0$, and $0 \le t$, $i$, $f \le 1$);

- the *classical set* (for $n = 1$ and $i = 0$, with $t$, $f$ either 0 or 1);

  - the *paraconsistent set* (for $n > 1$, with all $t,i,f < 1+$);

  - the *faillibilist set* ($i > 0$);

- the *dialetheist set*, a set M whose at least one of its elements also belongs to its complement C(M); thus, the intersection of some disjoint sets is not empty;

- the *paradoxist set* ($t=f=1$);

- the *pseudoparadoxist set* ($0 < i < 1$, $t=1$ and $f>0$ or $t>0$ and $f=1$);

- the *tautological set* ($i,f < 0$).

Compared with all other types of sets, in the neutrosophic set each element has three components which are subsets (not numbers as in fuzzy set) and considers a subset, similarly to intuitionistic fuzzy set, of "indeterminacy" - due to unexpected parameters hidden in some sets, and let the superior limits of the components to even boil *over 1* (overflooded) and the inferior limits of the components to even freeze *under 0* (underdried).

For example: an element in some tautological sets may have $t > 1$, called "overincluded". Similarly, an element in a set may be "overindeterminate" (for $i > 1$, in some paradoxist sets), "overexcluded" (for $f > 1$, in some unconditionally false appurtenances); or "undertrue" (for $t < 0$, in some unconditionally false appurtenances), "underindeterminate" (for $i < 0$, in some unconditionally true or false appurtenances), "underfalse" (for $f < 0$, in some unconditionally true appurtenances).

This is because we should make a distinction between unconditionally true ($t > 1$, and $f < 0$ or $i < 0$) and conditionally true appurtenances ($t \le 1$, and $f \le 1$ or $i \le 1$).

In a *rough set* RS, an element on its boundary-line cannot be classified neither as a member of RS nor of its complement with certainty. In the neutrosophic set a such element may be characterized by x(T, I, F), with corresponding set-values for T, I, F $\subseteq \,\|^-0, 1^+\|$.

*Abstract*:   In this paper one generalizes the classical and imprecise probability to neutrosophic probability, and similarly for neutrosophic statistics.  Examples are presented.  This is only a simple introduction to these concepts.

*Keywords and Phrases*: Classical Probability, Imprecise Probability, Neutroosphic Logic.

*1991 MSC*:  60A99.

## 4.  NEUTROSOPHIC PROBABILITY (NP):

### 4.1.  Definition:

Let T, I, F be real standard or non-standard subsets included in $\Vert {}^-0, 1^+\Vert$,

with    sup T = t_sup, inf T = t_inf,

sup  I = i_sup,  inf I = i_inf,

sup F = f_sup, inf F = f_inf,

and    n_sup = t_sup+i_sup+f_sup,

n_inf = t_inf+i_inf+f_inf.

The *neutrosophic probability* is a generalization of the classical probability and imprecise probability in which the chance that an event A occurs is t% true - where t varies in the subset T, i% indeterminate - where i varies in the subset I, and f% false - where f varies in the subset F.

Statically T, I, F are subsets, but dynamically they are functions/operators depending on many known or unknown parameters.

In classical probability n_sup ≤ 1, while in neutrosophic probability n_sup ≤ $3^+$.

In imprecise probability: the probability of an event is a subset T ⊂ [0, 1], not a number p ∈ [0, 1], what's left is supposed to be the opposite, subset F (also from the unit interval [0, 1]); there is no indeterminate subset I in imprecise probability.

One notes NP(A) = (T, I, F), a triple of sets.

### 4.2.  Neutrosophic Probability Space:

The universal set, endowed with a neutrosophic probability defined for each of its subset, forms a neutrosophic probability space.

Let A and B be two neutrosophic events, and NP(A) = $(T_1, I_1, F_1)$, NP(B) = $(T_2, I_2, F_2)$ their neutrosophic probabilities.  Then we define:

$(T_1, I_1, F_1) \boxplus (T_2, I_2, F_2) = (T_1 \oplus T_2, I_1 \oplus I_2, F_1 \oplus F_2)$,

$(T_1, I_1, F_1) \boxminus (T_2, I_2, F_2) = (T_1 \ominus T_2, I_1 \ominus I_2, F_1 \ominus F_2)$,

$(T_1, I_1, F_1) \boxdot (T_2, I_2, F_2) = (T_1 \odot T_2, I_1 \odot I_2, F_1 \odot F_2)$.

NP(A∩B) = NP(A) ⊡ NP(B);

NP(¬A)  = {1} ⊟ NP(A),  [this second axiom may be replaced, in specific applications, with a more intuitive one: NP(¬A) = (F$_1$, I$_1$,T$_1$)];

NP(A∪B) = NP(A) ⊞ NP(B) ⊟ NP(A) ⊡ NP(B).

Neutrosophic probability is a non-additive probability, i.e. P(A∪B) ≠ P(A)+P(B).

A probability-function P is called additive if P(A∪B) = P(A)+P(B), sub-additive if P(A∪B) ≤ P(A)+P(B), and super-additive if P(A∪B) ≥ P(A)+P(B).

In the Dempster-Shafer Theory P(A) + P(¬A) may be ≠ 1, in neutrosophic probability almost all the time P(A) + P(¬A) ≠ 1.

1. NP(impossible event) = (T$_{imp}$, I$_{imp}$, F$_{imp}$),
where sup T$_{imp}$ ≤ 0, inf F$_{imp}$ ≥ 1; no restriction on I$_{imp}$.

   NP(sure event) = (T$_{sur}$, I$_{sur}$, F$_{sur}$),
where inf T$_{sur}$ ≥ 1, sup F$_{sur}$ ≤ 0.

   NP(totally indeterminate event) = (T$_{ind}$, I$_{ind}$, F$_{ind}$);
where inf I$_{ind}$ ≥ 1; no restrictions on T$_{ind}$ or F$_{ind}$.

2. NP(A) ∈ {(T, I, F), where T, I, F are real standard or non-standard subsets included in $\|{}^-0, 1^+\|$ which may overlap}.

3. NP(A∪B) = NP(A) ⊞ NP(B) ⊟ NP(A∩B).

4. NP(A) = {1} ⊟ NP(¬A).

5. Other versions of Neutrosophic Probabilities can be defined by using other N-norms/N-conorms (see Section 2.9.1.). {In this case we have used the Algebraic Product N-norm/N-conorm.}

## 4.3. Applications:

#1. From a pool of refugees, waiting in a political refugee camp in Turkey to get the American visa, a% have the chance to be accepted - where a varies in the set A, r% to be rejected - where r varies in the set R, and p% to be in pending (not yet decided) - where p varies in P.

Say, for example, that the chance of someone Popescu in the pool to emigrate to USA is (between) 40-60% (considering different criteria of emigration one gets different percentages, we have to take care of all of them), the chance of being rejected is 20-25% or 30-35%, and the chance of being in pending is 10% or 20% or 30%. Then the neutrosophic probability that Popescu emigrates to the Unites States is

   NP(Popescu) = ( (40-60), (20-25)U(30-35), {10,20,30} ), closer to the life.

This is a better approach than the classical probability, where 40 ≤ P(Popescu) ≤ 60, because from the pending chance - which will be converted to acceptance or rejection - Popescu might get extra percentage in his will to emigration,

and also the superior limit of the subsets sum

   60+35+30 > 100

and in other cases one may have the inferior sum < 0,

while in the classical fuzzy set theory the superior sum should be 100 and the inferior sum ≥ 0.

In a similar way, we could say about the element Popescu that

Popescu( (40-60), (20-25)U(30-35), {10,20,30} ) belongs to the set of accepted refugees.

#2. The probability that candidate C will win an election is say 25-30% true (percent of people voting for him), 35% false (percent of people voting against him), and 40% or 41%

indeterminate (percent of people not coming to the ballot box, or giving a blank vote - not selecting anyone, or giving a negative vote - cutting all candidates on the list).

Dialectic and dualism don't work in this case anymore.

#3. Another example, the probability that tomorrow it will rain is say 50-54% true according to meteorologists who have investigated the past years' weather, 30 or 34-35% false according to today's very sunny and droughty summer, and 10 or 20% undecided (indeterminate).

#4. The probability that Yankees will win tomorrow versus Cowboys is 60% true (according to their confrontation's history giving Yankees' satisfaction), 30-32% false (supposing Cowboys are actually up to the mark, while Yankees are declining), and 10 or 11 or 12% indeterminate (left to the hazard: sickness of players, referee's mistakes, atmospheric conditions during the game).  These parameters act on players' psychology.

### 4.4.  Remarks:

Neutrosophic probability are useful to those events which involve some degree of indeterminacy (unknown) and more criteria of evaluation - as above.  This kind of probability is necessary because it provides a better approach than classical probability to uncertain events.

This probability uses a subset-approximation for the truth-value (like *imprecise probability*), but also subset-approximations for indeterminacy- and falsity-values.

Also, it makes a distinction between "*relative sure event*", event which is sure only in some particular world(s): NP(*rse*) = 1, and "*absolute sure event*", event which is sure in all possible worlds: NP(*ase*) = $1^+$; similarly for "*relative impossible event*" / "*absolute impossible event*", and for "*relative indeterminate event*" / "*absolute indeterminate event*".

In the case when the truth- and falsity-components are complementary, i.e. no indeterminacy and their sum is 100, one falls to the classical probability.  As, for example, tossing dice or coins, or drawing cards from a well-shuffled deck, or drawing balls from an urn.

### 4.5.  Generalizations:

An interesting particular case is for n = 1, with $0 \leq t, i, f \leq 1$, which is closer to the classical probability.

For n = 1 and i = 0, with $0 \leq t, f \leq 1$, one obtains the classical probability.

If I disappear and F is ignored, while the non-standard unit interval $\Vert ^-0, 1^+ \Vert$ is transformed into the classical unit interval [0, 1], one gets the imprecise probability.

From the intuitionistic logic, paraconsistent logic, dialetheism, fallibilism, paradoxism, pseudoparadoxism, and tautologism we transfer the  "adjectives" to probabilities, i.e. we define the ***intuitionistic probability*** (when the probability space is incomplete), ***paraconsistent probability***, ***faillibilist probability***, ***dialetheist probability***, ***paradoxist probability***, ***pseudoparadoxist probability***, and ***tautological probability*** respectively.

Hence, the neutrosophic probability generalizes:

- the *intuitionistic probability*, which supports incomplete (not completely known/determined) probability spaces (for $0 < n < 1$, $0 \leq t, f \leq 1$) or incomplete events whose probability we need to calculate;

- the *classical probability* (for n = 1 and i = 0, and $0 \leq t, f \leq 1$);

- the *paraconsistent probability* (for n > 1, with all t,i,f < 1);
  - the *dialetheist probability*, which says that intersection of some disjoint probability spaces is not empty (for t = f = 1 and i = 0; some paradoxist probabilities can be denoted this way);
  - the *faillibilist probability* (for i > 0);
- the *pseudoparadoxism* (for i>0, t=1 and 0<f<1 or 0<t<1 and f=1);
- the *tautologism* (for t> 1).

Compared with all other types of classical probabilities, the neutrosophic probability introduces a percentage of "indeterminacy" - due to unexpected parameters hidden in some probability spaces, and let each component t, i, f be even boiling *over 1* (overflooded) or freezing *under 0* (underdried).

For example: an element in some tautological probability space may have t > 1, called "overprobable". Similarly, an element in some paradoxist probability space may be "overindeterminate" (for i > 1), or "overunprobable" (for f > 1, in some unconditionally false appurtenances); or "underprobable" (for t < 0, in some unconditionally false appurtenances), "underindeterminate" (for i < 0, in some unconditionally true or false appurtenances), "underunprobable" (for f < 0, in some unconditionally true appurtenances).

This is because we should make a distinction between unconditionally true (t > 1, and f < 0 or I < 0) and conditionally true appurtenances (t ≤ 1, and f ≤ 1 or I ≤ 1).

### 4.6. NEUTROSOPHIC STATISTICS:

Analysis of events described by the neutrosophic probability means *neutrosophic statistics*.

The function that models the neutrosophic probability of a random variable x is called *neutrosophic distribution*: NP(x) = ( T(x), I(x), F(x) ), where T(x) represents the probability that value x occurs, F(x) represents the probability that value x does not occur, and I(x) represents the indeterminant/unknown probability of value x.

This is also a generalization of classical statistics.
[Future study is to be done in this subject…]

**Acknowledgements**:

The author would like to thank Drs. C. Le and Ivan Stojmenovic for encouragement and invitation to write this paper.

# 5.  DEFINITIONS DERIVED FROM NEUTROSOPHICS

*Abstract*:  Thirty-three new definitions are presented, derived from neutrosophic set, neutrosophic probability, neutrosophic statistics, and neutrosophic logic.
Each one is independent, short, with references and cross references like in a dictionary style.

*Keywords*:  Fuzzy set, fuzzy logic; neutrosophic logic;
Neutrosophic set, intuitionistic set, paraconsistent set, faillibilist set, paradoxist set, pseudo-paradoxist set, tautological set, nihilist set, dialetheist set, trivialist set;
Classical probability and statistics, imprecise probability;
Neutrosophic probability and statistics, intuitionistic probability and statistics, paraconsistent probability and statistics, faillibilist probability and statistics, paradoxist probability and statistics, pseudo-paradoxist probability and statistics, tautological probability and statistics, nihilist probability and statistics, dialetheist probability and statistics, trivialist probability and statistics;
Neutrosophic logic, paradoxist logic (or paradoxism), pseudo-paradoxist logic (or pseudo-paradoxism), tautological logic (or tautologism).

*2000 MSC*:  03E99, 03-99, 03B99, 60A99, 62A01, 62-99.

*Introduction*:
As an addenda to [1], [3], [5-7] we display the below unusual extension of definitions resulted from neutrosophics in the Set Theory, Probability, and Logic.  Some of them are listed in the Dictionary of Computing [2].  Further development of these definitions (including properties, applications, etc.) is in our research plan.

## 5.1. Definitions of New Sets

====================================================

### 5.1.1. Neutrosophic Set:

 <logic, mathematics> A set which generalizes many existing classes of sets, especially the fuzzy set.

Let U be a universe of discourse, and M a set included in U.
An element x from U is noted, with respect to the set M, as x(T,I,F),
and belongs to M in the following way: it is T% in the set
(membership appurtenance), I% indeterminate (unknown if it is in the
set), and F% not in the set (non-membership);
here T,I,F are real standard or non-standard subsets, included in the
non-standard unit interval ]-0, 1+[, representing truth,
indeterminacy, and falsity percentages respectively.

Therefore: -0 ≤ inf(T) + inf(I) + inf(F) ≤ sup(T) + sup(I) + sup(F) ≤ 3+.

Generalization of {classical set}, {fuzzy set}, {intuitionistic set},
{paraconsistent set}, {faillibilist set}, {paradoxist set},
{tautological set}, {nihilist set}, {dialetheist set}, {trivialist}.

Related to {neutrosophic logic}.

{ ref. Florentin Smarandache, "A Unifying Field in Logics.
Neutrosophy: Neutrosophic Probability, Set, and Logic",
American Research Press, Rehoboth, 1999;
(http://www.gallup.unm.edu/~smarandache/NeutrosophicSet.pdf,
 http://www.gallup.unm.edu/~smarandache/FirstNeutConf.htm,
 

===================================================

## 5.1.2. Intuitionistic Set:

<logic, mathematics> A set which provides incomplete
information on its elements.

A class of {neutrosophic set} in which every element x is
incompletely known, i.e. x(T,I,F) such that
sup(T)+sup(I)+sup(F)<1;
here T,I,F are real standard or non-standard subsets, included in
the non-standard unit interval ]-0, 1+[, representing truth,
indeterminacy, and falsity percentages respectively.

Contrast with {paraconsistent set}.

Related to {intuitionistic logic}.

{ ref. Florentin Smarandache, "A Unifying Field in Logics.
Neutrosophy: Neutrosophic Probability, Set, and Logic",
American Research Press, Rehoboth, 1999;
(http://www.gallup.unm.edu/~smarandache/FirstNeutConf.htm,

 http://www.gallup.unm.edu/~smarandache/neut-ad.htm) }

==================================================

### 5.1.3. Paraconsistent Set:

<logic, mathematics> A set which provides paraconsistent information on its elements.

A class of {neutrosophic set} in which every element x(T,I,F) has the property that sup(T)+sup(I)+sup(F)>1;
here T,I,F are real standard or non-standard subsets, included in the non-standard unit interval ]-0, 1+[, representing truth, indeterminacy, and falsity percentages respectively.

Contrast with {intuitionistic set}.

Related to {paraconsistent logic}.

==================================================

### 5.1.4. Faillibilist Set:

<logic, mathematics> A set whose elements are uncertain.

A class of {neutrosophic set} in which every element x has a percentage of indeterminacy, i.e. x(T,I,F) such that inf(I)>0;
here T,I,F are real standard or non-standard subsets, included in the non-standard unit interval ]-0, 1+[, representing truth, indeterminacy, and falsity percentages respectively.

Related to {fallibilism}.

==================================================

**5.1.5. Paradoxist Set**:

<logic, mathematics> A set which contains and doesn't contain itself at the same time.

A class of {neutrosophic set} in which every element x(T,I,F) has the form x(1,I,1), i.e. belongs 100% to the set and doesn't belong 100% to the set simultaneously;
here T,I,F are real standard or non-standard subsets, included in the non-standard unit interval ]-0, 1+[, representing truth, indeterminacy, and falsity percentages respectively.

Related to {paradoxism}.

{ ref. Florentin Smarandache, "A Unifying Field in Logics.
Neutrosophy: Neutrosophic Probability, Set, and Logic",
American Research Press, Rehoboth, 1999;
(http://www.gallup.unm.edu/~smarandache/FirstNeutConf.htm,
http://www.gallup.unm.edu/~smarandache/neut-ad.htm) }

==================================================

**5.1.6. Pseudo-Paradoxist Set**:

<logic, mathematics> A set which totally contains and partially doesn't contain itself at the same time,
or partially contains and totally doesn't contain itself at the same time.

A class of {neutrosophic set} in which every element x(T,I,F) has the form x(1,I,F) with 0<inf(F)≤sup(F)<1 or x(T,I,1) with 0<inf(T)≤sup(T)<1,
i.e. belongs 100% to the set and doesn't belong F% to the set simultaneously, with 0<inf(F)≤sup(F)<1,
or belongs T% to the set and doesn't belong 100% to the set simultaneously, with 0<inf(T)≤sup(T)<1;
here T,I,F are real standard or non-standard subsets, included in the non-standard unit interval ]-0, 1+[, representing truth, indeterminacy, and falsity percentages respectively.

Related to {pseudo-paradoxism}.

{ ref. Florentin Smarandache, "A Unifying Field in Logics.
Neutrosophy: Neutrosophic Probability, Set, and Logic",
American Research Press, Rehoboth, 1999;
(http://www.gallup.unm.edu/~smarandache/FirstNeutConf.htm,
http://www.gallup.unm.edu/~smarandache/neut-ad.htm) }

==================================================

**5.1.7. Tautological Set**:

<logic, mathematics> A set whose elements are absolutely determined in all possible worlds.

A class of {neutrosophic set} in which every element x has the form x(1+,-0,-0), i.e. absolutely belongs to the set; here T,I,F are real standard or non-standard subsets, included in the non-standard unit interval ]-0, 1+[, representing truth, indeterminacy, and falsity percentages respectively.

Contrast with {nihilist set} and {nihilism}.

Related to {tautologism}.

{ ref. Florentin Smarandache, "A Unifying Field in Logics. Neutrosophy: Neutrosophic Probability, Set, and Logic", American Research Press, Rehoboth, 1999; (http://www.gallup.unm.edu/~smarandache/FirstNeutConf.htm, http://www.gallup.unm.edu/~smarandache/neut-ad.htm) }

==================================================

**5.1.8. Nihilist Set**:

<logic, mathematics> A set whose elements absolutely don't belong to the set in all possible worlds.

A class of {neutrosophic set} in which every element x has the form x(-0,-0,1+), i.e. absolutely doesn't belongs to the set; here T,I,F are real standard or non-standard subsets, included in the non-standard unit interval ]-0, 1+[, representing truth, indeterminacy, and falsity percentages respectively.

The empty set is a particular set of {nihilist set}.

Contrast with {tautological set}.

Related to {nihilism}.

{ ref. Florentin Smarandache, "A Unifying Field in Logics. Neutrosophy: Neutrosophic Probability, Set, and Logic", American Research Press, Rehoboth, 1999;

==================================================

### 5.1.9. Dialetheist Set:

 <logic, mathematics> /di:-al-u-theist/ A set which contains at
least one element which also belongs to its complement.

A class of {neutrosophic set} which models a situation
where the intersection of some disjoint sets is not empty.

There is at least one element x(T,I,F) of the dialetheist set
M which belongs at the same time to M and to the set C(M),
which is the complement of M;
here T,I,F are real standard or non-standard subsets, included in the
non-standard unit interval ]-0, 1+[, representing truth,
indeterminacy, and falsity percentages respectively.

Contrast with {trivialist set}.

Related to {dialetheism}.

{ ref. Florentin Smarandache, "A Unifying Field in Logics.
Neutrosophy: Neutrosophic Probability, Set, and Logic",
American Research Press, Rehoboth, 1999;

==================================================

### 5.1.10. Trivialist Set:

 <logic, mathematics>  A set all of whose  elements also belong
to its complement.

A class of {neutrosophic set} which models a situation
where the intersection of any disjoint sets is not empty.

Every element x(T,I,F) of the trivialist set M belongs at the
same time to M and to the set C(M), which is the
complement of M;
here T,I,F are real standard or non-standard subsets,
included in the non-standard unit interval ]-0, 1+[, representing
truth, indeterminacy, and falsity percentages respectively.

Contrast with {dialetheist set}.

Related to {trivialism}.

{ ref. Florentin Smarandache, "A Unifying Field in Logics.
Neutrosophy: Neutrosophic Probability, Set, and Logic",
American Research Press, Rehoboth, 1999;
(http://www.gallup.unm.edu/~smarandache/FirstNeutConf.htm,
 http://www.gallup.unm.edu/~smarandache/neut-ad.htm) }

========================================================

### 5.2. Definitions of New Probabilities and Statistics

========================================================

## 5.2.1. Neutrosophic Probability:

<probability> The probability that an event occurs is (T, I, F),
where T,I,F are real standard or non-standard subsets, included in the
non-standard unit interval ]-0, 1+[, representing truth,
indeterminacy, and falsity percentages respectively.

Therefore: $-0 \leq \inf(T) + \inf(I) + \inf(F) \leq \sup(T) + \sup(I) + \sup(F) \leq 3+$.

Generalization of {classical probability} and {imprecise probability},
{intuitionistic probability}, {paraconsistent probability}, {faillibilist
probability}, {paradoxist probability}, {tautological probability},
{nihilistic probability}, {dialetheist probability}, {trivialist probability}.

Related with {neutrosophic set} and {neutrosophic logic}.

The analysis of neutrosophic events is called **Neutrosophic Statistics**.

{ ref. Florentin Smarandache, "A Unifying Field in Logics.
 Neutrosophy: Neutrosophic Probability, Set, and Logic",
 American Research Press, Rehoboth, 1999;
(http://www.gallup.unm.edu/~smarandache/FirstNeutConf.htm,
 http://www.gallup.unm.edu/~smarandache/neut-ad.htm) }

========================================================

## 5.2.2. Intuitionistic Probability:

<probability> The probability that an event occurs is (T, I, F),

where T,I,F are real standard or non-standard subsets, included in the non-standard unit interval ]-0, 1+[, representing truth, indeterminacy, and falsity percentages respectively, and n_sup = sup(T)+sup(I)+sup(F) < 1, i.e. the probability is incompletely calculated.

Contrast with {paraconsistent probability}.

Related to {intuitionistic set} and {intuitionistic logic}.

The analysis of intuitionistic events is called **Intuitionistic Statistics**.

{ ref. Florentin Smarandache, "A Unifying Field in Logics.
 Neutrosophy: Neutrosophic Probability, Set, and Logic",
 American Research Press, Rehoboth, 1999;
(http://www.gallup.unm.edu/~smarandache/FirstNeutConf.htm,
 http://www.gallup.unm.edu/~smarandache/neut-ad.htm) }

==================================================

### 5.2.3. Paraconsistent Probability:

<probability> The probability that an event occurs is (T, I, F),
where T,I,F are real standard or non-standard subsets, included in the
non-standard unit interval ]-0, 1+[, representing truth,
indeterminacy, and falsity percentages respectively,
and n_sup = sup(T)+sup(I)+sup(F) > 1,
i.e. contradictory information from various sources.

Contrast with {intuitionistic probability}.

Related to {paraconsistent set} and {paraconsistent logic}.

The analysis of paraconsistent events is called
**Paraconsistent Statistics**.

{ ref. Florentin Smarandache, "A Unifying Field in Logics.
 Neutrosophy: Neutrosophic Probability, Set, and Logic",
 American Research Press, Rehoboth, 1999;
(http://www.gallup.unm.edu/~smarandache/FirstNeutConf.htm,
 http://www.gallup.unm.edu/~smarandache/neut-ad.htm) }

==================================================

### 5.2.4. Faillibilist Probability:

<probability> The probability that an event occurs is (T, I, F),
where T,I,F are real standard or non-standard subsets, included in the
non-standard unit interval ]-0, 1+[, representing truth,
indeterminacy, and falsity percentages respectively,
and inf(I) > 0,
i.e. there is some percentage of indeterminacy in calculation.

Related to {faillibilist set} and {fallibilism}.

The analysis of faillibilist events is called **Faillibilist Statistics**.

{ ref. Florentin Smarandache, "A Unifying Field in Logics.
 Neutrosophy: Neutrosophic Probability, Set, and Logic",
 American Research Press, Rehoboth, 1999;
 (http://www.gallup.unm.edu/~smarandache/FirstNeutConf.htm,
  http://www.gallup.unm.edu/~smarandache/neut-ad.htm) }

================================================

**5.2.5. Paradoxist Probability**:

<probability> The probability that an event occurs is (1, I, 1),
where I is a standard or non-standard subset, included in the
non-standard unit interval ]-0, 1+[, representing indeterminacy.

Paradoxist probability is used for paradoxal events (i.e. which
may occur and may not occur simultaneously).

Related to {paradoxist set} and {paradoxism}.

The analysis of paradoxist events is called **Paradoxist Statistics**.

{ ref. Florentin Smarandache, "A Unifying Field in Logics.
 Neutrosophy: Neutrosophic Probability, Set, and Logic",
 American Research Press, Rehoboth, 1999;
(http://www.gallup.unm.edu/~smarandache/FirstNeutConf.htm,
  http://www.gallup.unm.edu/~smarandache/neut-ad.htm) }

================================================

**5.2.6. Pseudo-Paradoxist Probability**:

<probability> The probability that an event occurs is either (1, I, F) with
0<inf(F)≤sup(F)<1, or (T, I, 1) with 0<inf(T)≤sup(T)<1,
where T,I,F are standard or non-standard subset, included in the
non-standard unit interval ]-0, 1+[, representing the truth, indeterminacy, and
falsity percentages respectively.

Pseudo-Paradoxist probability is used for pseudo-paradoxal events (i.e. which
may certainly occur and may not partially occur simultaneously,
or may partially occur and may not certainly occur simultaneously).

Related to {pseudo-paradoxist set} and {pseudo-paradoxism}.

The analysis of pseudo-paradoxist events is called **Pseudo-Paradoxist Statistics**.

{ ref. Florentin Smarandache, "A Unifying Field in Logics.
 Neutrosophy: Neutrosophic Probability, Set, and Logic",
 American Research Press, Rehoboth, 1999;
(http://www.gallup.unm.edu/~smarandache/FirstNeutConf.htm,
 http://www.gallup.unm.edu/~smarandache/neut-ad.htm) }

===================================================

## 5.2.7. Tautological Probability:

<probability> The probability that an event occurs is more than one,
i.e. (1+, -0, -0).

Tautological probability is used for universally sure events (in all
possible worlds, i.e. do not depend on time, space, subjectivity, etc.).

Contrast with {nihilistic probability} and {nihilism}.

Related to {tautological set} and {tautologism}.

The analysis of tautological events is called **Tautological Statistics**.

{ ref. Florentin Smarandache, "A Unifying Field in Logics.
 Neutrosophy: Neutrosophic Probability, Set, and Logic",
 American Research Press, Rehoboth, 1999;
(http://www.gallup.unm.edu/~smarandache/FirstNeutConf.htm,
 http://www.gallup.unm.edu/~smarandache/neut-ad.htm) }

===================================================

## 5.2.8. Nihilist Probability:

<probability> The probability that an event occurs is less than zero, i.e. (-0, -0, 1+).

Nihilist probability is used for universally impossible events (in all possible worlds, i.e. do not depend on time, space, subjectivity, etc.).

Contrast with {tautological probability} and {tautologism}.

Related to {nihilist set} and {nihilism}.

The analysis of nihilist events is called **Nihilist Statistics**.

{ ref. Florentin Smarandache, "A Unifying Field in Logics.
 Neutrosophy: Neutrosophic Probability, Set, and Logic",
 American Research Press, Rehoboth, 1999;
(http://www.gallup.unm.edu/~smarandache/FirstNeutConf.htm,
 http://www.gallup.unm.edu/~smarandache/neut-ad.htm) }

===================================================

## 5.2.9. Dialetheist Probability:

<probability> /di:-al-u-theist/ A probability space where at least one event and its complement are not disjoint.

A class of {neutrosophic probability} that models a situation
 where the intersection of some disjoint events is not empty.

 Here, similarly, the probability of an event to occur is (T, I, F),
 where T,I,F are real standard or non-standard subsets, included
 in the non-standard unit interval ]-0, 1+[, representing truth,
 indeterminacy, and  falsity percentages respectively.

Contrast with {trivialist probability}.

Related to {dialetheist set} and {dialetheism}.

The analysis of dialetheist events is called **Dialetheist Statistics**.

{ ref. Florentin Smarandache, "A Unifying Field in Logics.
 Neutrosophy: Neutrosophic Probability, Set, and Logic",
 American Research Press, Rehoboth, 1999;
(http://www.gallup.unm.edu/~smarandache/FirstNeutConf.htm,
 http://www.gallup.unm.edu/~smarandache/neut-ad.htm) }

=====================================================

**5.2.10. Trivialist Probability**:

<probability> A probability space where every event and its
complement are not disjoint.

A class of {neutrosophic probability}which models a situation
where the intersection of any disjoint events is not empty.

Here, similarly, the probability of an event to occur is (T, I, F),
where T,I,F are real standard or non-standard subsets, included
in the non-standard unit interval ]-0, 1+[, representing truth,
indeterminacy, and  falsity percentages respectively.

Contrast with {dialetheist probability}.

Related to {trivialist set} and {trivialism}.

The analysis of trivialist events is called **Trivialist Statistics**.

{ ref. Florentin Smarandache, "A Unifying Field in Logics.
 Neutrosophy: Neutrosophic Probability, Set, and Logic",
 American Research Press, Rehoboth, 1999;
(http://www.gallup.unm.edu/~smarandache/FirstNeutConf.htm,
 http://www.gallup.unm.edu/~smarandache/neut-ad.htm) }

=====================================================

### *5.3. Definitions of New Logics*

=====================================================

**5.3.1. Neutrosophic Logic**:

<logic, mathematics> A logic which generalizes many existing classes
of logics, especially the fuzzy logic.

In this logic each proposition is estimated to have the percentage of truth in
a subset T, the percentage of indeterminacy in a subset I, and the percentage
of falsity in a subset F;
here T,I,F are real standard or non-standard subsets, included in the
non-standard unit interval ]-0, 1+[, representing truth,
indeterminacy, and falsity percentages respectively.

Therefore: $-0 \leq \inf(T) + \inf(I) + \inf(F) \leq \sup(T) + \sup(I) + \sup(F) \leq 3+$.

Generalization of {classical or Boolean logic}, {fuzzy logic},
{multiple-valued logic}, {intuitionistic logic}, {paraconsistent logic},
{faillibilist logic, or failibilism}, {paradoxist logic, or paradoxism},
{pseudo-paradoxist logic, or pseudo-paradoxism}, {tautological logic, or
tautologism}, {nihilist logic, or nihilism}, {dialetheist logic, or dialetheism},
{trivialist logic, or trivialism}.

Related to {neutrosophic set}.

{ ref. Florentin Smarandache, "A Unifying Field in Logics.
Neutrosophy: Neutrosophic Probability, Set, and Logic",
American Research Press, Rehoboth, 1999;
(http://www.gallup.unm.edu/~smarandache/NeutrosophicLogic.pdf,
http://www.gallup.unm.edu/~smarandache/FirstNeutConf.htm,
http://www.gallup.unm.edu/~smarandache/neut-ad.htm) }

===================================================

### 5.3.2. Paradoxist Logic (or Paradoxism):

<logic, mathematics> A logic devoted to paradoxes, in which each
proposition has the logical vector value (1, I, 1);
here I is a real standard or non-standard subset, included in the
non-standard unit interval ]-0, 1+[, representing the indeterminacy.

As seen, each paradoxist (paradoxal) proposition is true and false
simultaneously.

Related to {paradoxist set}.

{ ref. Florentin Smarandache, "A Unifying Field in Logics.
Neutrosophy: Neutrosophic Probability, Set, and Logic",
American Research Press, Rehoboth, 1999;
(http://www.gallup.unm.edu/~smarandache/NeutrosophicLogic.pdf,
http://www.gallup.unm.edu/~smarandache/FirstNeutConf.htm,
http://www.gallup.unm.edu/~smarandache/neut-ad.htm) }

===================================================

### 5.3.3. Pseudo-Paradoxist Logic (or Pseudo-Paradoxism):

<logic, mathematics> A logic devoted to pseudo-paradoxes,
in which each proposition has the logical vector value:
either (1, I, F), with $0 < \inf(F) \leq \sup(F) < 1$,
or (T, I, 1), with $0 < \inf(T) \leq \sup(T) < 1$;

here I is a real standard or non-standard subset, included in the
non-standard unit interval ]-0, 1+[, representing the indeterminacy.

As seen, each pseudo-paradoxist (pseudo-paradoxal) proposition is:
either totally true and partially false simultaneously,
or partially true and totally false simultaneously.

Related to {pseudo-paradoxist set}.

{ ref. Florentin Smarandache, "A Unifying Field in Logics.
Neutrosophy: Neutrosophic Probability, Set, and Logic",
American Research Press, Rehoboth, 1999;
(http://www.gallup.unm.edu/~smarandache/NeutrosophicLogic.pdf,
http://www.gallup.unm.edu/~smarandache/FirstNeutConf.htm,
http://www.gallup.unm.edu/~smarandache/neut-ad.htm) }

==================================================

### 5.3.4. Tautological Logic (or Tautologism):

<logic, mathematics> A logic devoted to tautologies, in which each
proposition has the logical vector value (1+, -0, -0).

As seen, each tautological proposition is absolutely true (i. e, true in all
possible worlds).

Related to {tautological set}.

{ ref. Florentin Smarandache, "A Unifying Field in Logics.
Neutrosophy: Neutrosophic Probability, Set, and Logic",
American Research Press, Rehoboth, 1999;
(http://www.gallup.unm.edu/~smarandache/NeutrosophicLogic.pdf,
http://www.gallup.unm.edu/~smarandache/FirstNeutConf.htm,
http://www.gallup.unm.edu/~smarandache/neut-ad.htm) }

==================================================

General References:

1.  Jean Dezert, *Open Questions on Neutrosophic Inference*, Multiple-Valued Logic
Journal, 2001 (to appear).
2.  Denis Howe, *On-Line Dictionary of Computing*,
http://foldoc.doc.ic.ac.uk/foldoc/
3.  C. Le, *Preamble to Neutrosophy and Neutrosophic Logic*, Multiple-Valued
Logic Journal, 2001 (to appear).

.

**General References**:

**Neutrosophy** is a new branch of philosophy, introduced by Dr. Florentin Smarandache in 1995, which studies the origin, nature, and scope of neutralities, as well as their interactions with different ideational spectra.

**Neutrosophic Logic** is a general framework for unification of many existing logics. The main idea of NL is to characterize each logical statement in a 3D Neutrosophic Space, where each dimension of the space represents respectively the truth (T), the falsehood (F), and the indeterminacy (I) of the statement under consideration, where T, I, F are standard or non-standard real subsets of ]-0, 1+[.

**Neutrosophic Set**.
Let U be a universe of discourse, and M a set included in U. An element x from U is noted with respect to the set M as x(T, I, F) and belongs to M in the following way:
it is t% true in the set, i% indeterminate (unknown if it is) in the set, and f% false, where t varies in T, i varies in I, f varies in F.
Statically T, I, F are subsets, but dynamically T, I, F are functions/operators depending on many known or unknown parameters.

**Neutrosophic Probability** is a generalization of the classical probability and imprecise probability in which the chance that an event A occurs is t% true - where t varies in the subset T, i% indeterminate - where i varies in the subset I, and f% false - where f varies in the subset F.
In classical probability n_sup $\leq$ 1, while in neutrosophic probability n_sup $\leq 3^+$.
In imprecise probability: the probability of an event is a subset T $\subset$ [0, 1], not a number p $\in$ [0, 1], what's left is supposed to be the opposite, subset F (also from the unit interval [0, 1]); there is no indeterminate subset I in imprecise probability.

**Neutrosophic Statistics** is the analysis of events described by the neutrosophic probability.
The function that models the neutrosophic probability of a random variable x is called *neutrosophic distribution*: NP(x) = ( T(x), I(x), F(x) ), where T(x) represents the probability that value x occurs, F(x) represents the probability that value x does not occur, and I(x) represents the indeterminant / unknown probability of value x.

Dr. M. Khoshnevisan, Griffith University, Australia